\documentclass[10pt]{amsart}
\usepackage[colorlinks=true, urlcolor=blue,bookmarks=true,bookmarksopen=true,
citecolor=blue,hypertex]{hyperref}

\newcommand{\R}{{\mathbb{R}}}

\newcommand{\N}{{\mathbb{N}}}

\newcommand{\Z}{{\mathbb{Z}}}
\newcommand{\C}{{\mathbb{C}}}

\newcommand{\p}{{\partial}}
\newcommand{\al}{{\alpha}}

\newcommand{\be}{{\beta}}

\newcommand{\om}{{\omega}}
\newcommand{\eps}{{\varepsilon}}

\newcommand{\De}{{\Delta}}
\newcommand{\ga}{{\gamma}}
\newcommand{\Ga}{{\Gamma}}

\newcommand{\la}{{\lambda}}
\newcommand{\si}{{\sigma}}

\newcommand{\Aa}{{\mathcal{A}}}
\newcommand{\Bb}{{\mathcal{B}}}
\newcommand{\Cc}{{\mathcal{C}}}
\newcommand{\Cl}{\mathcal{C}\mathcal{\ell} \,}

\newcommand{\Dd}{{\mathcal{D}}}

\newcommand{\Jj}{{\mathcal{J}}}
\newcommand{\Ee}{{\mathcal{E}}}
\newcommand{\Ff}{{\mathcal{F}}}

\newcommand{\Hh}{{\mathcal{H}}}
\newcommand{\Nn}{{\mathcal{N}}}
\newcommand{\Mm}{{\mathcal{M}}}

\newcommand{\Oo}{{\mathcal{O}}}
\newcommand{\Pp}{{\mathcal{P}}}

\newcommand{\Vv}{{\mathcal{V}}}

\newcommand{\Tt}{{\mathcal{T}}}
\newcommand{\Ww}{{\mathcal{W}}}

\newcommand{\Q}{{\mathbb{Q}}}

\newcommand{\La}{{\Lambda}}

\newcommand{\Si}{{\Sigma}}
\newcommand{\n}{{\mathbf{n}}}
\newcommand{\Nu}{{\hspace{.1cm} ^\nu \hspace{-.1cm}}}

\newcommand{\cyl}{{\rm cyl}}

\newcommand{\ind}{{\rm ind}}
\newcommand{\ev}{{\rm ev}}

\newcommand{\Crit}{{\rm Crit}}

\newcommand{\MS}{{\medskip}}
\newcommand{\BS}{{\bigskip}}
\newcommand{\NI}{{\noindent}}

\newcommand{\QED}{\hfill$\Box$\medskip}
\newtheorem{theorem}{Theorem}
\newtheorem{corollary}{Corollary}[section]
\newtheorem{definition}{Definition}
\newtheorem{remark}[corollary]{Remark}
\newtheorem{lemma}[corollary]{Lemma}

\newtheorem{proposition}[corollary]{Proposition}

\newtheorem{example}[corollary]{Example}

\newcommand{\FF}{{ \,  I  \hspace{-.11cm}   { F}  }}

\begin{document}

\title{Cluster homology}

  \author{Octav Cornea and Fran\c{c}ois Lalonde}

\date{August 2005}

\address{\newline\indent O.C. and F.L.\newline\indent University of Montr\'eal\newline
\indent Department of Mathematics and Statistics
\newline \indent C.P. 6128 Succ. Centre-ville
\newline \indent Montr\'eal  H3C 3J7
\newline\indent Qu\'ebec, Canada}

\maketitle

\tableofcontents

\

\section*{Introduction.}\label{sec:introduction}   We introduce in
this paper a new homology attached to a Langrangian submanifold $L$
of a symplectic manifold,
which provides an invariant of the Hamiltonian isotopy class of $L$. This leads 
to various applications concerning the analytical, topological, and dynamical properties of 
Lagrangian submanifolds, and to a new universal Floer homology,
defined without obstruction, for pairs of Lagrangian submanifolds.

Throughout this paper, we shall assume that
all Lagrangian submanifolds are connected, orientable, and relatively spin
(recall that a Lagrangian submanifold
$L\subset (M, \om)$ is relatively spin if the second
Steifel-Withney class of $L$ admits an extension to $H^2(M;
\Z_2)$; a set of Lagrangian submanifolds is relatively spin if their second
Stiefel-Whitney classes admit a common extension to
$H^{2}(M;\Z/2)$). Actually, the notation $L$ for
a Lagrangian submanifold will always implicitely contain the information of a choice of an orientation
and of a relatively spin structure (and the same applies for a set of such submanifolds). 
The ambient symplectic manifold $(M^{2n},\om)$ is supposed to be compact or, if not, it should be either convex at infinity or geometrically bounded in the sense of  \cite{ALP}, so that no sequence of  Riemann surfaces with boundary lying on a set of compact Lagrangian submanifolds  $L_1, \ldots , L_{\ell} \subset M$ can escape to infinity. In fact, the Lagrangian submanifolds need not all be compact, as long as the above control on sequences of Riemann surfaces is ensured.

Our construction starts with the remark that the
usual coefficient rings used to define Floer type complexes are
not rich enough to properly encode the bubbling of disks. Thus,
for a relatively spin Lagrangian submanifold $L$ in an a
symplectic manifold $(M,\omega)$, we first introduce a rational
differential graded algebra, {\em the cluster algebra} of $L$,
which manages algebraically the bubbling of pseudoholomorphic
disks with boundary on $L$. For two relatively spin Lagrangian
submanifolds
$L$, $L'$, we essentially use the tensor product of the cluster
algebras associated to $L$ and to $L'$ as a new coefficient ring
of a Floer type homology which we call the {\em fine Floer}
homology of the pair $(L,L')$.

\

The cluster algebra of $L$, that we shall denote by $\Cl_*(L,J,f)$, depends
on auxiliary data; these data include, in particular, the choice of an almost complex
structure $J$ compatible with the symplectic form $\om$, and the choice of a 
Morse-Smale pair $(f,g)$ consisting of a Morse
function $f:L\to \R$ and of a Riemannian metric $g$ on $L$. As an
algebra over an appropriate (rational) Novikov ring, it is graded
commutative and is freely generated, as a symmetric algebra, 
by the critical points of the
Morse function $f$. Its homology $\Cl H_{\ast}(L)$, called {\em
the cluster homology}, is of interest in itself: we
show that it is independent of $J$ and $(f,g)$ and is invariant with
respect to symplectic diffeomorphisms. In particular, $\Cl H_*(L)$ is an invariant of
the Hamiltonian isotopy class of the Lagrangian embedding $L \subset M$.

The key idea in the definition of the cluster algebra has its origin
in the observation that the main difficulty in modelling
algebraically the bubbling of disks is that this is a
non-localized co-dimension one phenomenon. In other words, the
bubbling of disks produces boundary components - which, from our
point of view, is not problematic - but  in contrast with
the usual breaking of Morse-flow lines, the points where the
bubbling appears are arbitrary. The solution that we propose to this
problem is, in essence, to enlarge the moduli space of
pseudo-holomorphic disks by allowing configurations formed of
disks connected by negative gradient flow lines of an 
a priori fixed Morse function on our Lagrangian submanifold. These new moduli
spaces, called {\em clustered moduli spaces}, are large enough
to transform the bubbling of disks into an internal,
co-dimension one phenomenon. The compactifications of these spaces
certainly have boundaries but they consist only of pairs of
configurations joined at some critical point of the fixed Morse
function.

In short, in this way, the bubbling of disks is absorbed away from
the boundary components of our big moduli spaces -- the boundary
components of the latter are now only associated to the ordinary
breaking mechanism of flow lines of the fixed Morse function at
some critical point. Since these points are known and are actually
the algebra generators of our cluster complex, this enables us to
define a differential $d$ satisfying $d^2 = 0$ as well as the
Leibniz rule. In \S~\ref{sec:cluster_moduli}, we will describe in
detail these moduli spaces.  They are made of smaller ``pieces"
glued together along various codimension 0 strata of the
boundaries, which motivates the name ``cluster'' .

\

Once the cluster algebra introduced and its main properties
identified, it is not difficult to imagine how to define  a
universal Floer theory, the {\it fine Floer homology}, attached to
the pair formed by two relatively spin Lagrangian submanifolds $L,
L' \subset (M, \om)$ in general position: it should be the
homology associated to a chain complex freely generated by
(possibly, some of) the intersection points of $L$ and $L'$ over a
ring roughly (but not exactly) of the form $\Cl_{\ast}(L,f,J)\otimes
\Cl_{\ast}(L',f',J)$. As we shall see, this idea may be
successfully pursued and the putative differential of the fine
Floer chain complex has indeed a vanishing square. This vanishing,
however, is due to a cancellation resulting from of subtle phenomenon
that has appeared before just as a nuisance in standard Floer theory. 
In its simplest form, this phenomenon is
seen in a one dimensional space of pseudoholomorphic strips
resting on $L$ and on $L'$ and joining an intersection point
$x\in L\cap L'$ to itself. It consists in the fact that such a
moduli space might admit as boundary point a pseudoholomorphic
disk with boundary on $L$ and passing through $x$.

When $L'$ is Hamiltonian isotopic to $L$, it is possible to
identify the cluster algebras of $L$ and $L'$ and to use the
resulting symmetry to simplify the construction. We call {\em symmetric fine Floer
homology} the resulting homology which, in general, is different from the fine Floer homology. We establish a natural algebraic relation between the
cluster and the symmetric fine Floer homologies, called the {\it strip-string
symmetrization}. The definition of the fine Floer homology, symmetric or not, makes
it obvious that it vanishes for a disjoint pair $L,L'$. Using this fact, we
prove that, under additional purely algebraic topological restrictions,
a relatively spin, orientable Lagrangian submanifold that can be
displaced from itself by a Hamiltonian isotopy is uniruled. 
This means that through each point $p$ of such a Lagrangian submanifold
$L \subset (M, \om)$, and for any $\om$-compatible almost complex
structure $J$, there is a $J$-holomorphic disk  whose boundary on
$L$ passes through $p$. In fact, there is also a natural
action associated to the fine Floer complex and using it we show
that the disks detected above are of symplectic area at most equal
to the disjunction energy of $L$. An interesting geometric
consequence of this fact, along ideas from Barraud-Cornea \cite{BaCo},  
is that the displacement energy of such a Lagrangian submanifold $L$  
admits as lower bound the square of the real Gromov radius of $L$ (see
\S\ref{subsec:bded_disks}) multiplied by $\pi/2$.
We will also discuss a different application to the detection of
periodic orbits of Hamiltonian flows. This is based on the construction of a
chain morphism relating the cluster homologies of two Lagrangian submanifolds
$L$ and $L'$
by making use of moduli spaces  which integrate the clusters on $L$ and $L'$ with $J$-holomorphic cylinders with their
ends on the two Lagrangian submanifolds, thus pursuing the work in Gatien-Lalonde \cite{GL}.

\

A theory such as this one consists of a certain number of - by now
familiar - levels: 

\NI
\emph{a}. a geometric phenomenon is modeled by
introducing certain moduli spaces; 

\NI
\emph{b}. compactness
properties of these moduli spaces are proved; 

\NI
\emph{c}. regularity
properties of the moduli spaces are established - in general by
proving certain transversality results; 

\NI
\emph{d}. the structure of
the moduli spaces is used to define (hopefully efficient)
algebraic invariants. 

\smallskip
To make the paper as readable and, at the
same time, as complete as possible, we will present here
reasonably complete proofs having to do with the points \emph{a},
\emph{b} and \emph{d} as well as those of a number of applications,
but we will postpone the analytical theorems concerning the point
\emph{c} to a forthcoming article based on Hofer-Wisocki-Zehnder's new
scale-polyfold Fredholm theory. Again, to improve readability, in
the first section we outline many  of our constructions and we
also prove a number of results which follow rapidly from the
general properties of the objects introduced. We come back to the
technical details in later sections.

\

{\bf Acknowledgements.} We thank Katrin Wehrheim for useful discussions and
for insisting on the fact that the constant strips need to be included in the construction
of the fine Floer homology.

\

\section{Overview of the main constructions and results.}

\subsection{Definition and properties of the cluster complex}
\label{subsec:cluster_definition}

The cluster complex is associated to a triple formed by
(1) a Lagrangian embedding  $L^{n}\hookrightarrow (M,\omega)$, equipped with
a choice of an orientation and of a relatively spin structure (which, we recall, is always assumed in the notation $L$), (2) a
 generic almost complex
 structure $J$ on $M$ compatible wth the symplectic form $\om$, and (3) a  pair $(f,g)$ with  $f:L\to \R$ a Morse
function and $g$ a Riemannian metric on $L$ so that  $(f,g)$ is
Morse-Smale. The two
conditions implicit in the notation $L$ (i.e the orientation 
and the relative spin structure) are needed to orient the 
clustered moduli spaces - to be described below - in a coherent way 
(see \S~\ref{subsec:orientation} for the discussion of the orientations).

This complex is denoted by $\Cl_{\ast}(L, J; (f,g))$ and, if
$\Crit(f)$ is the set of critical points of $f$,  we set:
$$
\Cl_{\ast}(L,J;(f,g)) = (S\Q<s^{-1}\ \Crit(f)>
\otimes\ \Lambda)^{\wedge}
$$
where $s^{-1}\ \Crit(f)$ indicates that the natural index grading
of $\Crit(f)$ is decreased  by one unit, $SV$ is the free, graded
commutative algebra over the graded vector space $V$ (as usual,
the sign commutativity rule is $ab=(-1)^{|a||b|}ba$ for any two
elements $a,b\in V$), $^{\wedge}$ indicates a certain completion
described below and $\Lambda$ is an appropriate group ring defined
as follows: on $\pi_2(M,L)$, consider the equivalence relation
$\la \sim \tau$ iff $\om(\la) = \om (\tau)$ and $\mu(\la) = \mu
(\tau)$, where $\om$ and $\tau$ are the area and Maslov classes
respectively. The ring $\La$ is the rational group ring of
$\pi_2(M,L)/ \sim$.  We write the elements of $\La$ under the form
of finite sums $\sum_{i}c_{i}e^{\la_i}$, $c_{i}\in \Q$. Clearly, given the
equivalence relation on $\pi_2(M,L)$, the variable $e^{\la}$ could
as well be replaced by $t_1^{\mu} t_2^{a}$ where $\mu, a$ are the
Maslov and area numbers in the ranges of
$$
\mu: \pi_2(M,L) \to \Z  \qquad \qquad  \om: \pi_2(M,L) \to \R.
$$

The grading in $\La$ is given by $|e^{\la}|=-\mu(\la)$ for $\la\in
\pi_{2}(M,L)/\sim$. With this convention, the grading of the cluster
complex is given by the usual tensor product formula. Thus for
$x_i \in \Crit(f)$, we have $$|x_i|=\ind_{f}(x_i)-1 \ , \
|x_{1}\ldots x_{k}e^{\la}|=\sum_{i=1}^{k} |x_{i}|-\mu(\la)~.~$$

Finally, the completion $\ ^{\wedge}$ is associated to the
filtration given by:

\begin{equation}\label{eq:completion}
L^{k}(S\Q<s^{-1}\ \Crit(f)>\otimes\ \La)= \Q< x_{1}x_{2}\ldots
x_{s} e^{\la}\ : \ s\geq k \ {\rm or}\ \omega(\la)\geq k >~.~
\end{equation}

 \

 In particular, notice that an element $m\in \Cl(L,J;(f,g))$
 can be written as a possibly infinite sum:
 $$m=m_{0}+m_{1}e^{\la_{1}}+\ldots + m_{i}e^{\la_{i}}+\ldots$$
 where $m_{i}$ are monomials in the elements of $\Crit(f)$
 but if this sum is infinite, then any infinite
 subsequence with $\om(\la_i)$ bounded above, must have
  its corresponding word length sequence converging
  to infinity (in the sense that any of its subsequences
  is unbounded). Conversely, any formal sum verifying  this condition
 belongs to the cluster complex.

 In \S~\ref{sec:cluster_moduli}  we will introduce in detail the moduli
  spaces on which is based the definition of the cluster differential. 
We now summarize the construction.

    The generic data  $J, (f,g)$ on $L$ being given,  fix
an order on the critical points of $f$. Choose any integer
$k \ge 0$, any $x \in \Crit(f)$, any non-decreasing sequence of
critical points $x_1, \ldots, x_k$, and $\la \in \pi_{2}(M,L)/\sim$, with the sole 
constraint that the zero class $\la = 0$ is allowed only when $k$ equals $1$.  
Consider then 
the space $\Mm^x_{x_1, \ldots, x_k} (\la)$ consisting of all
configurations of the following kind.  Let $\Tt$ be a finite
connected tree with oriented edges such that one and only one
vertex has no ingoing edge  (it is called the {\it root}) and all
other vertices have exactly one ingoing edge. All vertices can have any number $\ge 0$ of
outgoing ones. To each edge $\al \be$ (the writing means that
the edge joins the vertex $\al$ to the vertex $\be$ in this
order), let $t_{\al \be}$ be a real number in $ [0, \infty)$
and $z_{\al\be}, z_{\be\al}$ be two points on the boundary of the
unit disk $D$ in the $\C$-plane. Consider now the abstract
topological space $\Tt'$ obtained by replacing each vertex $\ga$
by a copy $D_{\ga}$ of the unit disk $D$,  view the point $z_{\al
\be}$ as a point of $\p D_{\al}$ and likewise $z_{\be\al}$ as a
point of $\p D_{\be}$, and replace each edge $\al \be$ by an edge
joining the point $z_{\al \be} \in \p D_{\al}$ to the point
$z_{\be \al} \in \p D_{\be}$. The points $z_{\al \be}, z_{\be
\al}$ are called {\it incidence points}.

    Let $z, z_1, \ldots, z_k$ be $k+1$ additional distinct
{\it marked points}  lying on the
 boundaries of the disks $D_{\ga}$'s  with the
 sole constraints that they are also distinct from the incidence points and that $z$ belongs to the boundary of the root disk.
 Now consider the pair $(\mathbf u, \mathbf z)$ where $\mathbf z$ is
 the sequence $z, z_1, \ldots, z_k \in \Tt'$ and $u: \Tt' \to M$ is
 a map which satisfies the following properties:

\begin{itemize}
\item[1.] for each $\ga$, $u_{\ga} := u |_{D_{\ga}}$ is a $J$-holomorphic map
from $(D_{\ga}, \p D_{\ga})$ to $(M,L)$ whose class is denoted by $\la_{\ga} \in \pi_2(M,L)$;
\item[2.] the restriction of $u$ to each edge $\al \be$ is an
 integral curve of the negative gradient vector field
 of $(f,g)$ when the edge is parametrized by the interval
 $[0, t_{\al\be}]$;
\item[3.] the sum of the classes $\la_{\ga}$ over all vertices is equal to $\la$;
\item[4.] the point $u(z)$ belongs to the unstable manifold of the critical point $x$
while each point $u(z_i), 1 \le i \le k$, belongs to the stable manifold of $x_i$
(i.e one may view this as an extension of the tree by adding one semi-infinite
outgoing flowline from $x$ to the root disc at $z$, and $k$ semi-infinite
outgoing flowlines from the $z_i$'s to the $x_i$'s).
\end{itemize}

There are other more technical constraints related to (i) stability (ghost disks), (ii)
codimension $\ge 2$ phenomena - which we have neglected in the description above,
and (iii) virtual perturbations with coherent orientations,
see \S~\ref{sec:cluster_moduli} for more detail and \S~\ref{subsec:orientation} where
the orientations of these moduli spaces are described. It is worth noticing here that there are two
transversality issues in this theory: the first one is the usual lack of transversality due to multiple
coverings of disks; the second one arises when two special points meet. In this second case, the problem is that the two cluster constraints become identified at that meeting point, so this has the undesirable  effect of increasing the dimension of the moduli space at that point. In both cases, the problem is solved by the Deligne-Mumford setting at the source, which produces ghost disks whenever needed (see 
\S~\ref{sec:cluster_cplx_hom} for details). Note that, in contrast with the absolute case, the second transversality problem arises in generic  one-dimensional families, hence can never be neglected in our theory.

The space $\Mm^x_{x_1, \ldots, x_k} (\la)$ is then the  topological space made of all
pairs $(\mathbf u, \mathbf z)$ quotiented out by the natural automorphisms.
As we shall see, these moduli spaces admit natural compactifications
$\bar{\Mm}^x_{x_1, \ldots, x_k} (\la)$. To express this, it is useful to 
introduce the notation $S$ for a partially ordered set of critical points 
(in which repetitions are allowed) which of course can be identified with 
a unique non-decreasing sequence of points. If $S', S''$ are two such ordered 
subsets, we will denote by $<S' \cup S'' > $ the ordered subset made of the 
elements in $S' \cup S''$. Letting $S$ be the ordered set $\{x_1, \ldots, x_k\}$, the points 
in the top dimensional strata of 
$\bar{\Mm}^x_{S} (\la)\backslash \Mm^x_{S} (\la)$ can 
be  identified with pairs belonging to some $\bar\Mm^x_{<S',\{y\}>}(\la'))\times
\bar\Mm^y_{S''}(\la''))$ where $\la'+\la''=\la$, $y\in \Crit(f)$, and $<S'\cup S''>=S$.
The corresponding moduli spaces obtained by using the perturbed
equation $\bar{\partial}_{J}u=\nu(u)$ are denoted by $\Nu \Mm^x_{x_1, \ldots, x_k} (\la)$.
It is possible  to find perturbations so that the moduli spaces
 $\Nu \bar{\Mm}^x_{x_1, \ldots, x_k} (\la)$
carry the structure of oriented singular manifolds with boundary,  admitting
a fundamental class over the rationals (with boundary) -- see
Definition~\ref{def:orientations} of \S~\ref{sec:cluster_moduli} for the
definition of the orientations of these moduli spaces.
The dimension of   $\Nu\bar{\Mm}^x_{x_1, \ldots, x_k} (\la)$ equals
\begin{equation}\label{eq:dim}
|x|-\sum_{i=1}^{k}|x_{i}|+\mu(\la)-1  ;
\end{equation}
and, moreover, for those moduli spaces of dimension $1$, we have:
\begin{equation}\label{eq:dsquare1}
\begin{array}{l}
\hspace{0.5in}\partial
(^{\nu}\bar{\Mm}^{x}_{S}(\la))=\\[1ex]
= \bigcup_{S=<S'\cup S''>,y, \la'+\la''=\la}\
(^{\nu}\bar\Mm^x_{<S',y>}(\la'))\times
(^{\nu}\bar\Mm^y_{S''}(\la'')).
\end{array}
\end{equation}

Here the summation is taken over all $y \in \Crit(f)$, all partitions of
$S$ into two subsets $S', S''$ and all the
splittings of $\la$ as the sum of two classes $\la',  \la''$.
This equation holds with the orientations refered above, where
the orientation on the left hand side is of course the one induced on the boundary.

Of course, by Gromov's compactification theorem, only a finite number of the moduli
spaces appearing on the right hand side of the last equation are non-empty.

\BS

Given the moduli spaces described above, we define the cluster
differential
$$d^{\nu}:(S\Q<s^{-1}\ \Crit(f)>\otimes\ \Lambda)^{\wedge}_{\ast}\to
((S\Q<s^{-1}\ \Crit(f)>)\otimes\ \Lambda)^{\wedge}_{\ast-1}$$ as
the unique commutative, graded differential algebra extension
of:
\begin{equation}\label{eq:diff}
d^{\nu}x=\sum_{\la, k \ge 0, x_{1},\ldots,
x_{k}}a^{x}_{x_{1},\ldots,x_{k}}(\la)x_{1}\ldots x_{k} e^{\la}
\end{equation}
where $x,x_{1},\ldots,x_{k}\in \Crit(f)$ have the property that
$(x_{1},\ldots,x_{k})$ respects the fixed order on $\Crit(f)$,
$|x|-\sum_{i} |x_{i}|+\mu(\la)=1$ and
$$a^{x}_{x_{1},\ldots,x_{k}}=\#
(^{\nu}\Mm^{x}_{x_{1},\ldots,x_{k}}(\la)) \in \Q.$$ In this
definition we count the elements in
$^{\nu}\Mm^{x}_{x_{1},\ldots,x_{k}}(\la)$ with signs. 

Note that, therefore, with our conventions, we have:
$$d^\nu (xy)=(d^\nu x)\cdot y
+(-1)^{|x|}x\cdot (d^\nu y)$$.

\begin{theorem}\label{theo:d2} The map $d^{\nu}$ satisfies $(d^{\nu})^2 = 0$.
\end{theorem}

This is an immediate consequence of formula (\ref{eq:dsquare1}) above
together with the orientation conventions described in \S \ref{subsec:orientation}.
The verification is left to the reader.
We denote by $\Cl H_*(L, (f,g),J, \nu)$ the homology of this complex,
and call it the {\it cluster homology of} $L$.
We prove in \S~\ref{sec:cluster_invariance} the invariance of
this homology with respect to changes in the choices of
the auxiliary data:  $(f,g),J, \nu$. Notice that this implies that
this homology is also invariant under symplectic diffeomorphism in the
sense that if $\phi:M\to M$ is a symplectic diffeomorphism, then
$\Cl H_{\ast}(L)\simeq \Cl H_{\ast}(\phi(L))$. However, it does depend in general on the choice of the orientation of $L$ and of its relative spin structure.

\subsubsection{First remarks and an example}

We start with a simple observation which will play an important role
later in this article.

\begin{remark}\label{rem:cl_moduli}{\em
 A critical point $y$ of index $0$ can
{\em never} appear as end in a $0$-dimensional, non empty, moduli
space of type
$$\Nu\Mm^{x}_{\ldots,y,\ldots}(\la)$$
(except for usual Morse flow lines). Indeed, consider the obvious
forgetful map

\begin{equation}\label{eq:forget}
^{\nu}\Mm^{x}_{x_{1},\ldots,y\ldots x_{k}}(\la)\to
^{\nu}\Mm^{x}_{x_{1},\ldots,x_{k}}(\la)
\end{equation}
and notice that, due to (\ref{eq:dim}), the domain of this map is of dimension greater by
one unit than its image, a contradiction since this would imply that the right hand side is 
of negative dimension. Similarly, a critical point $y$ of index $n=\dim (L)$ can not appear
as the root in a $0$-dimensional moduli space of the form $\Mm^{y}_{\ldots}(\la)$ (except,
again, for usual Morse flow lines).  Of course, another, more intuitive way of seeing 
this is to realize that if there is a negative flowline from a non-constant loop 
in $L$ to a local minimum, there must be a one-parameter family of these, 
and similarly for the relative maxima.} 
\end{remark}

The next point is quite significant from an algebraic point of view.

\begin{remark}\label{rem:completion}
{\em The completion $^{\wedge}$ is necessary, in
particular because if $x_{0}\in \Crit(f)$ verifies $|x_{0}|=0$ and
if $\dim (^{\nu}\Mm^{x}_{x_{0},x_{1},\ldots,x_{k}}(\la))=0$, then
we also have $$\dim
(^{\nu}\Mm^{x}_{x_{0},\ldots,x_{0},x_{1},\ldots,x_{k}}(\la))=0$$
and so in our definition of $d$ we need to allow formal powers in
the critical points of $f$ which are of index $1$. A critical
point of index at least $2$ cannot appear infinitely many times
as endpoint of a non-vanishing $0$-dimensional moduli space
because of formula (\ref{eq:dim}). However, for algebraic reasons,
it is preferable to complete using the word length filtration
which counts all generators as in (\ref{eq:completion}) and not
only those generators  of degree $1$.}
\end{remark}

An element in the cluster complex, $\tau\in \Cl(L,J;(f,g))$,
is written, in general, as a sum $\tau=(\sum_{\la}a(\la)e^{\la})+m$ where $a(\la)\in
\Q$ and $m$ is a sum of words (in the letters consisting of the critical points of $f$)
of length at least one. We call each of the terms
$a(\la)e^{\la}$  with $a(\la)\not=0$ {\em a free term} of $\tau$ and if  there is a critical
point $x$ of $f$ whose differential contains at least one free term, we say that the {\em complex has free terms}. 
A particular
case will be significant: if  the Morse index of $x$ in the above definition is larger or equal to $1$, 
we will say that  {\em the cluster complex has
high free terms}. Notice that, due to the fact that $\mu(\la)$ is even,
if a critical point $x$ verifies $dx=a_{0}e^{\la}+...$, $a_{0}\not=0$, then
 $\ind_{f}(x)$ is even. Moreover, in view of Remark \ref{rem:cl_moduli},
$\ind_{f}(x)\not=\dim(L)$. To summarize, the complex has high free terms if and only if there is a critical point of even Morse index different from $0$ and $n$ whose cluster differential has a free term.

\begin{proposition}\label{prop:free_terms}
The cluster complex is acyclic ( that is $\Cl H_{\ast}(L)=0$) if
and only if the cluster has free
terms. If the cluster complex $\Cl (L,J,f)$ has high free terms, then
there are, moreover, some $J$-disks of non-positive Maslov class (i.e with Maslov index $\le 0$).
\end{proposition}

\begin{proof} Suppose that the differential of $x\in \Crit(f)$
has a free
term. We write $$dx=\sum_{\la}a^{x}(\la)e^{\la}+m$$ with $m$ a
sum of words of length at least one. For each free term we have
$\mu(\la)-1=|x|$. Therefore, by the definition of the equivalence
relation used to define $\La$, there exists a unique $\la_{0}$ so
that $\omega(\la_{0})$ is minimal among all the free terms in
$dx$. It follows that for $\tau= xe^{-\la_{0}}$ we have
$$d\tau=1+\sum_{\omega(\la)> 0}b(\la)e^{\la}+m'$$
 where $m'$ is again a sum of monomials of length at least $1$.
Set $$b= \sum_{\omega(\la)> 0}b(\la)e^{\la}+m'$$ and note that
 $c=\sum_{i=0}^{\infty}(-1)^{i}b^{i}\in\Cl(L,J;(f,g))$. Clearly,
 $db=0$ and so $dc=0$. Moreover, $d(c\tau)=cd(\tau)=1$, which
 implies that $\Cl H_{\ast}(L)=0$.

 Conversely, if $\Cl H_{\ast}(L)=0$, then $1$ is a
 boundary. But this is possible only if the cluster
 differential of a generator $x\in
\Crit(f)$ has a non-vanishing free term.

 The second part of the statement is obvious: if $dx=a_{0}e^{\la}+\ldots$, then
$\mu(\la)=2-\ind(x)$ and, as $\mu(\la)$ is even, the claim follows.
\end{proof}

\begin{remark}\label{rem:free_Max}{\rm 
It is useful to make explicit the following property of the cluster complex $\Cl(L,J,(f,g))$
where $f$ is a Morse function with a single local minimum $m$ and a single local maximum $M$ (and any number of critical points of intermediate indices).
By the last proposition, if $\Cl H_{\ast}(L)=0$, then
$dm\not=0$ (and there is therefore a $J$-holomorphic disk passing through $m$) or the cluster complex has high free terms.
The proposition above also shows that the existence or not of 
free terms in the cluster complex provides a fundamental dichotomy on the algebraic side of our
theory.}
\end{remark}

\begin{example}{\rm  If $S^{1}$ is a circle in $\C$ we have
 $$\Cl H_{\ast}(S^{1})=0~.~$$
}\end{example}
Indeed, take on $S^{1}$ the perfect Morse function with
one minimum $m$ and one maximum $M$. There exists one pseudoholomorphic
disk passing through $m$, of Maslov index $2$, with class in $\La$
that we will denote by $\la_{0}$. For the maximum, we have $dM=0$. 
The differential of the minimum can be seen to be given
by
$dm= (1+M+M^{2}+M^{3}+....)e^{\la_{0}}$ which, by the proposition, implies the claim.

\begin{example}\label{rem:no_bubbbling}{\rm  In the absence of bubbling (for example
if $\omega|_{\pi_{2}(M,L)}=0$), then $$\Cl H_{\ast}(L,J;(f,g)) \simeq S(s^{-1}H_{\ast}(L;\Q))\otimes\La)^{\wedge}~.~$$
This happens, of course, because in this case the only component of the cluster
differential is provided by the usual Morse differential.
}\end{example}

\subsubsection{Relations to other constructions.} \label{subsec:relation_to_other}

\

A. Formally, at the algebraic level, our construction has some 
similarities with the construction of
contact homology for Legendrian submanifolds which has been
developed by Eckholm, Etnyre and Sullivan \cite{EES}.
It is an interesting open problem at this time to see
whether this analogy is in fact the reflection of a deeper
relationship.

\

B. There are evident relations between our construction and the
 $A^{\infty}$ - machinery of Fukaya, Oh, Ohta and Ono \cite{FOOO}, at least due to the fact that
both theories model algebraically the bubbling of disks. However,
the underlying moduli spaces are quite different in the two theories, which
reflects the fact that the solution that we propose here to lift the obstructions
to the $d^2 = 0$ equation relies on a different, geometric, idea. While, at
some point, some direct relation between the two approaches might be discovered, they
remain, for the moment, complementary.

\

C. Spheres and disks with ``spikes" have been used in various ways
by other authors. Most notable is the announcement made by Oh in
\cite{Oh}: see our \S~\ref{subsec:relation_floer} below in which
we relate the cluster homology and  Oh's version of Floer homology as
described in \cite{Oh} when the minimal Maslov number is larger or
equal to $2$. Other instances are the work of Schwarz
\cite{Schw} as well as Piunikin-Salamon-Schwarz \cite{PSS} and
more recently Bourgeois \cite{Bou}.

\subsection{Fine Floer Homology.} \label{subsec:fine_definition}

The purpose of this paragraph is to introduce a variant of Floer
homology called {\em fine Floer homology} - denoted $\FF
H_{\ast}(-)$ - which is associated to a pair of transversal, orientable
 Lagrangian submanifolds $L_0, L_1$ of a symplectic manifold $(M, \om)$,
under the sole assumption that they be  relatively spin in $M$. Let us recall that this
means by definition that there is a common extension to $M$ of the second
Stiefel-Whitney classes of $L_0$ and $L_1$. This homology is defined,
roughly, using $\Cl(L_0)$ and $\Cl(L_1)$ as coefficients rings via a
Novikov ring that embraces all three terms $L_0, L_1$ and  $L_0 \cap
L_1$.
Note that we do not require $L_0$ to be homeomorphic to $L_1$ nor any
condition on the minimal Maslov numbers -- in this sense, this is a
universal Floer theory. As we shall see, its definition is entirely natural, up to two
delicate points: (1) one must be careful in the definition of the right
Novikov ring, and, more importantly, (2) to verify $d^{2}=0$, we will need to consider
moduli spaces $\Ww^a_{\ldots;b} (\la)$ of strips joining intersection points
$a,b$ to which are attached clusters lying on $L_{0}$ and on $L_{1}$.
These moduli spaces have undesirable boundary components consisting of
of $J$-holomorphic strips joinning $a \in L_0
\cap L_1$ to $b  \in L_0 \cap L_1$ with one gradient flow line in, say, $L_0$
stemming from, say, the point $a$ and leading in $L_0$ to some cluster
configuration. In general, this type of configuration has real
codimension $1$ and cannot therefore be neglected. We will see below how one can cancel them
to get $d^2 = 0$, even in the most general case when $L_0$ and $L_1$ are not even homeomorphic.

\

\subsubsection{Coefficient ring and moduli spaces.}\label{subsubesec:fine_var}
To define the fine Floer complex $$\FF C (L_0, L_1, \eta; J, (f_0,g_0), (f_1,g_1))$$
we first recall that the choices of orientations and of a relative spin structure for $L_0,L_1$
 have been made and are included in the notation $L_0,L_1$. Besides this, we need
auxiliary data  as follows. First, as before, we need an almost complex structure $J$,
Morse-Smale pairs  $(f_i,g_i)$ on $L_i$ and coherent choices of perturbations.
We also assume the $f_{i}$ in generic position with respect to the intersection points
$L_{0}\cap L_{1}$ in the sense that these intersection points are included in the unstable
manifolds of critical points of index $0$ of $f_{i}$. 
We denote by $\Ga = \{ \al:[0,1]
\to M : \al(i) \in L_i, i=0,1 \}$ the  space of continuous  paths from
$L_0$ to $L_1$. Here, $\eta$ is an element in $\Ga$ - its choice means
that we fix a basepoint for this space. We denote by $\Ga_{\eta}$ the
connected component of  $\Gamma$ which contains $\eta$. We denote by $I(L_{0},L_{1})$
the intersection points between $L_{0}$ and $L_{1}$ and we let $I_{\eta}$ be
those intersection points which, viewed as constant paths,
belong to $\Ga_{\eta}$. The generators
of the fine Floer complex will be precisely the elements of $I_{\eta}$.
Up to a shift in degrees,
the resulting fine Floer homology will only depend on $L_{0}$, $L_{1}$,
 the connected component of $\eta$ and the choice of orientations and relative spin structures of $L_0,L_1$.

\

To continue the construction, note that there are two group morphisms
$$\omega :\pi_{1}\Ga_{\eta}\to \R \ , \ \mu:\pi_{1}\Ga_{\eta}\to \Z$$
the first given by integration of $\omega$ and the second
obtained as follows.
Fix a path $l_{t}$ of Lagrangian subspaces of $T_{\eta(t)}M$ so that $l_{i}=T_{\eta(i)}L_{i}$,
$i=0,1$.
An element $[u]\in \pi_{1} \Ga_{\eta}$
is represented by a map $u:[0,1]\times [0,1] \to M$ so that
$u(s,0)\in L_{0}$, $u(s,1)\in L_{1}$, $\forall s$ and $u(0,t)=u(1,t)=\eta(t), \ \forall t$.
Therefore, $u^{\ast}TM$ is trivial and, after fixing one such trivialization, define
$\mu(u)$ as the Maslov index of the loop obtained by the
concatenation of the following paths of Lagrangian subspaces: $u(s,0)^{\ast}\ TL_{0}$,
$u(1,t)^{\ast}\ l_{t}$, $u(1-s,1)^{\ast}\ TL_{1}$ and $u(0,1-t)^{\ast}\ l_{1-t}$.
It is easy to see that $\mu(u)$ only depends on the homotopy class of $u$.
Let $K$ be the kernel of the product morphism $\Psi=\mu\times \omega$ and let $\bar\Ga_{\eta}$ be
the regular covering of $\Ga_{\eta}$ corresponding to $K$. Let $\Pi$ be the image of $\Psi$ and
let $\bar{\La}$ be the rational group ring of $\Pi$. 
For any point $\xi\in\bar\Ga_{\eta}$ we may define $\mu(\xi)=\mu(\xi,\eta)=\mu(\gamma)$
where $\gamma$ is a path in $\Ga_{\eta}$ which joins $\eta$ to $\pi(\xi)$ in the class defined by $\xi$ (see
Robbin-Salamon \cite{RoSa}). This number may be a half-integer.

We now define
\begin{equation}\label{eq:coefficients_fine}
\mathcal{R}=(S\Q<s^{-1}\Crit(f_{i}): \ i=0,1>\otimes \bar\La)^{\wedge}
\end{equation}
where the completion is as in (\ref{eq:completion}) except that we take into consideration
both critical points of $f_{0}$ and of $f_{1}$.

\

Notice that there are injective group morphisms
$\phi_{i}:\pi_{2}(M,L_{i})/\sim\to \Pi$ which are obtained by
first assuming that $\eta$ joins the base points in $L_{0}$ and $L_{1}$ and then viewing a
disk with boundary in, say, $L_{0}$ as a cylinder whose end on $L_{1}$ is constant. Therefore,
if we denote by $\La_{i}$ the group ring of $\pi_{2}(M,L_{i})/\sim$, we have injective ring
morphisms $\phi_{i}:\La_{i}\to \bar{\La}$. Thus,
$\mathcal{R}$ is isomorphic to $$\Cl (L_{0},f_{0},J)\otimes \Cl (L_{1},f_{1},J)\otimes
_{\La_{0}\otimes \La_{1}}\bar\La~.~$$ In other words, $\mathcal{R}$ is obtained
by replacing $\La_{0}\otimes \La_{1}$ in the tensor product $\Cl (L_{0},f_{0},J)\otimes \Cl (L_{1},f_{1},J)$ by $\bar\La$ by means of the ring morphism
$\phi_{0}\otimes\phi_{1}$ (all tensor products here are in the category of
graded, commutative algebras). This also implies that, by the usual Leibniz formula,
the cluster differentials on $\Cl(L_{i},f_{i},J)$ induce a graded, commutative algebra
differential $\delta$ on $\mathcal{R}$. Clearly, the differential algebra $(\mathcal{R},\delta)$
depends of all the choices made till now but to ease notation we will only indicate these
choices when necessary.

\subsubsection{The fine Floer complex}\label{subsubsec:fine_def_first}
 This is the free differential, graded module over $(\mathcal{R},\delta)$
given by
$$\FF C(L_0,L_1,\eta; J, (f_0,g_0), (f_1,g_1))=(\mathcal{R}\otimes \Q<I_{\eta}>, d_{F})~.~$$
The grading of the elements in $I_{\eta}$ is obtained as follows: we consider lifts
$\bar a \in\bar\Ga_{\eta}$ of the points $a\in I_{\eta}\subset \Ga_{\eta}$ (here, as usual, the last inclusion means that we view intersection points as constant paths)
and we define $|a|=\mu(\bar a)$.
Clearly, this grading depends on the choices of the lifts. 

 \begin{remark}\label{rem:fine_coefficient}{\rm It is useful to note that, although our 
 definition of the class $\la\in \Pi$
 depends on the non-canonical choice of the lifts $\bar{a}$, 
 different choices produce isomorphic complexes. Of course, 
 we could as well have defined the complex by taking as generators
 pairs $(a, \al)$ made of a point $a \in I_{\eta}$ and 
 of a lift of the constant path $a$ to an element $\al$ 
 of the covering $\bar{\Ga}_{\eta}$. We would have 
 then obtained a complex on which the group ring $\bar{\La}$ acts in a natural way.
 }
 \end{remark}

We now describe the differential $d_{F}$.
We order the critical points in $\Crit(f_{i})$. The differential $d_{F}$ verifies the Leibniz formula
and for an element $a\in I_{\eta}$ it is of the form:
$$d_{F}a=\sum_{\la, b, k \ge 0, l \ge 0, (x_{1},\ldots,x_{k},y_{1}\ldots,y_{l})}w^{a}_{x_{1},\ldots,
x_{k}, y_{1},\ldots, y_{l};b}(\la) \, x_{1}\ldots x_{k}y_{1}\ldots y_{l}b e^{\la}$$
where the $x_i$'s belong to $ \Crit(f_{0})$, the $y_j$'s to $\Crit(f_{1})$, they
respect the order, $\la \in \Pi$, and finally $b\in I_{\eta}$.

The coefficients $w^{a}_{x_{1},\ldots,x_{k},y_{1},\ldots,y_{l};b}(\la) \in \Q$
count the number of elements in certain $0$-dimensional moduli spaces
$\Ww^{a}_{x_{1},\ldots,x_{k},y_{1},\ldots,y_{l};b}(\la)$ (again after perturbation).
These moduli spaces are defined in a way similar to the $\Mm_{\ldots}(\la)$'s of
\S\ref{subsec:cluster_definition}. We will describe them formally in more detail later
in this paper in  \S\ref{sec:fine_moduli}. We now only indicate the main idea in their construction.
The starting point consists again of  trees as in \S\ref{subsec:cluster_definition} but
the root vertex $v_{0}$ of the tree corresponds to the closed unit disk $D_0$ with {\em two special distinct marked points} $z_0$ and $z_-$ on its boundary. The restriction $u_0 := u |_{D_0}$ maps continuously $D_0$ to $M$ in such a way that it sends $z_{0}$ to $a$,
$z_{-}$ to $b$. Moreover, on the punctured disk $D_{0}\backslash \{z_{0},z_{-}\}$,
$u_0$ verifies the equation of a pseudoholomorphic
strip $u_0:\R\times [0,1]\to M$ with $u(\R,i)\subset L_{i}$, after the reparametrization 
$$
\R\times [0,1] \to D_{0}\backslash \{z_{0},z_{-}\}.
$$
 Except for
codimension two phenomena, each of the
other vertices correspond to pseudoholomorphic disks whith boundaries  on one
of the $L_{i}$'s. Moreover, the gradient flows appearing in the construction correspond
 to one of the two functions $f_{i}$. In short, the elements of these moduli
 spaces are cluster trees on $L_{0}$ and $L_{1}$ that originate, {\em at finite points} from a single strip. 
 Note that there may be (finitely) many such clusters attached to the boundary of 
 the strip (and we still require that the corresponding incidence points on the boundary 
 of the source of the strip $u_0$, i.e the points where the gradient 
 flowlines of $f_0$ or $f_1$ start, be all distinct). We call them
 {\em cluster-strips}.  We also need to associate a class $\la\in\Pi$ to such an
 object. Notice that, such a clustered strip may be viewed, topologically, as a strip
 joining $a$ to $b$.  Thus it lifts to a path in $\bar{\Ga}_{\eta}$ which starts
 at $\bar{a}$ and ends at some element of the form $\bar{b}e^{\la}$. Precisely this $\la$
 is the class associated to the clustered-strip. With this definition we also see that
 $\la$ equals the sum of the classes of the disks included in the clustered strip (these
 classes are well-defined as we have seen before) added to the class - defined as before - of 
 the root $J$-strip that joins $a$ to $b$. The space $\Ww^{a}_{x_{1},\ldots,x_{k},y_{1},\ldots,y_{l};b}(\la)$ is by definition the set of all unparametrized cluster-strips from $\bar a $ to $\bar b$ labelled, in the above sense, by the class $\la$.

 In our formulae, it is important to distinguish the case in which
 the root is a constant strip and we fix the following convention: we will exclude form the moduli spaces $\Ww$'s, that serve to define our differential, all configurations in which the root is the constant strip at some point $a \in I_{\eta}$ and for which there is either no cluster attached to the strip, or more than one. In other words, when the strip is constant, we keep  
 only those configurations  in which {\em exactly one}
 cluster tree is attached to the strip $a$. In particular
 this cluster tree is {\em included in just one of the $L_{i}$'s}.
We obviously  have two types of such elements: the first, grouped
 in the moduli spaces
 $T^{a}_{x_{1},\ldots,x_{k}}(\la)$, consisting of
 configurations of one cluster tree in $L_{0}$ originating
 at $a$, and a similar moduli space  $R^{a}_{y_{1},\ldots,y_{l}}(\la)$ corresponding
 to trees in $L_{1}$. Thus,
  with this convention, the union
  $T^{a}_{x_{1},\ldots,x_{k}}(\la)\cup R^{a}_{y_{1},\ldots,y_{l}}(\la)$ consists of all
  the elements inside some of the $\Ww^{a}_{\ldots;a}(\la)$ for which the
  root is the constant strip. Note finally that, since the source map for the strip $u_0$ is $D_0 - \{z_0,z_-\} \simeq \R \times [0,1]$, our definition implies that a cluster tree in $L_0$ or in $L_1$ always originate from a flowline at a {\em finite } point of the source (i.e distinct from $z_0, z_-$). But, of course, configurations made of a strip with one cluster anchored at infinity, i.e anchored at the endpoints $a$ or $b$ of the strip will appear in the compactifications of our moduli spaces $\Ww^{a}_{x_{1},\ldots,x_{k},y_{1},\ldots,y_{l};b}(\la)$.

   Clearly, the moduli spaces $\Ww$ depend on all of our choices
  and if we need to indicate these choices explicitly we will write:
  $\Ww^{a}_{\ldots;b}((L_{0},f_{0}),(L_{1},f_{1}), J;\la)$.

 \

 For generic choices of $J$, $(f_{i},g_{i})$ and after pertubation,
 the dimension of the moduli space $\Ww^{a}_{x_{1},\ldots,x_{k},y_{1},\ldots,y_{l};b}(\la)$
 is: $$|a|-|b|-\sum |x_{i}|-\sum |y_{j}|-1~.~$$
 For one-dimensional moduli spaces, there is a formula analogue to (\ref{eq:dsquare1}).
 It will be stated and proved in \S\ref{sec:fine_moduli}. As a consequence, we have:

 \begin{corollary}\label{cor:fine-d-square}
 With the notation above we have: $d_{F}^{2}=0$.
 \end{corollary}

Verifying this formula is less immediate than for the cluster
differential because, besides the usual breaking of clusters and
of strips, there is, as we just mentionned above, a third potential way for boundary points to
emerge: they correspond to some cluster tree attached to a strip
at some moving point $p$ which ``slides" along the boundary of the 
strip to one of the ends of the strip.

\

There are two reasons that make these boundary 
components disappear, one is purely algebraic  
and the other one is analytic and consists in
the fact that (as remarked by Oh \cite{Oh1}) 
the usual gluing
argument applies (under generic conditions) 
to a $J$-disk passing
(transversally) through $a$ and 
to $a$ itself viewed as a constant
strip.  To see how these phenomena 
are related, we outline
below the argument proving $d_{F}^{2}=0$.

\

\subsubsection{Verification of $d_{F}^{2}=0$, remarks and an example}
\label{subsubsec:diff_fine1} 


Let $a$ be
a generator of the fine Floer complex, i.e $a \in I_{\eta}$. We have 
$d_F a=(A_{0}-A_{1})a+\sum  w^{a}_{b}b$ here
$A_{0}$ represents the elements in $T^{a}_{\ldots}$, $A_{1}$ the
elements in $R^{a}_{\ldots}$ and the $w^{a}_{b}$ represent the
other moduli spaces which involve only cluster-strips with a
non-constant root (note that, in these non-constant strips, $b$ might be equal to $a$).
We rewrite
this decomposition as $d_F a=  s(a) a +  \delta a$. Note that
$|A_{i}|=-1$ and so $(s(a))^{2}=(A_{0}-A_{1})^{2}=0$, which is due to the 
graded-commutative algebraic setting that we chose (there is no geometric reason that would 
cancel these terms). We now write
$$
(d_F)^{2}a = d(s(a))a + (s(a))(\delta a)+ \sum (dw^{a}_{b})b +\sum w^{a}_{b}(s(b)b)+\sum
w^{a}_{b} (\delta b)~.~
$$ To show that this vanishes, we need to identify, for each
such element in $d_F^2 (a)$, a way to include it as the $0$-boundary end of some one-parameter family
$C_{t \in [0,1]}$ sitting inside the compactification of one of our $\Ww$'s 
where $C_{t \in (0,1)}$ lies in the interior $\Ww$ of $\bar{\Ww}$, 
and where the configuration $C_1$ also corresponds to an element in the above 
expression of $d_F^2(a)$.  It is easy to see
that all the terms in the last four factors appear as boundary
points of $1$-dimensional moduli spaces $\Ww^{a}_{\ldots}$ of
cluster-strips with a non constant root: the terms in $\sum
(dw^{a}_{b})b$ correspond to broken flowlines inside the cluster trees;
the terms in $\sum w^{a}_{b}(\delta b)$ correspond to strip
breaking; those in $ \sum w^{a}_{b}(s(b)b)$ appear  as boundaries
when the  point where some cluster tree is attached to a non-constant strip
slides to $b$; similarly, the terms $(s(a))(\delta a)$ appear
as boundaries when the point  where some cluster tree is attached to a non-constant strip
slides to $a$.

The terms of the fifth type - corresponding to $d(s(a))a$ - appear as top dimensional 
components in the boundary of the compactification of the moduli space $T^{a}_{\ldots}\cup
R^{a}_{\ldots} \subset \Ww^{a}_{\ldots;a}$.  This creates an apparent problem because
these spaces  have an additional type of
top dimensional boundary elements: these appear when the length of the flowline
joining $a$ to the first disk in the cluster tends to $0$.
Thus, this boundary component of $T^{a}_{\ldots}\cup R^{a}_{\ldots}$ 
consists of cluster trees attached to  $a$ so that the first  disk in the cluster goes
through $a$. However, in view of the gluing result mentioned above,
these elements are in fact interior points in the larger moduli space
$\bar{\mathcal{W}}^{a}_{\ldots,a}$. It is clear that
the gluing construction produces a cluster-strip in which 
the cluster is attached at a finite point of the source of the non-constant
strip joining $a$ to $a$, i.e {\em away from $a$} (hence this one-parameter 
family will degenerate to an element of $d_F^2(a)$, not to an element of $d^3_F(a))$.

It is also useful to note that any $0$-dimensional configuration corresponding to an 
element in the above formula for $d_F^2(a)$ which is of a ``sliding form'' may be 
included in a unique one-parameter family: this is because the union of 
a strip $u$ from $a$ to $b$ with a cluster $T(a,x_1, \ldots, x_k,\la)$ attached at $a$ in $L_0$ say, 
is of dimension $0$ iff both $u$ and $T(a,x_1, \ldots, x_k,\la)$ are $0$-dimensional; but, by 
transversality, this means that there is locally in $L_0$ a $n$-dimensional family 
of such $T(p,x_1, \ldots, x_k,\la), p \in L_0$. This implies that 
the pairs $(u, T(p,x_1, \ldots, x_k,\la)), p \in \p_{L_0} u$, form a $1$-dimensional space. 
Conversely, and for the same reason, any one-dimensional family of 
cluster-strips that degenerates to a pair $(u,T(a,x_1, \ldots, x_k,\la))$ by 
sliding must be such that both the class that contains $u$ and the one 
that contains $T(a,x_1, \ldots, x_k,\la)$ are $0$-dimensional.

To summarize, this proves that the 
five types of boundary components exhaust the boundary points of a $1$-dimensional
moduli space -- this concludes the sketch of the proof.  \QED

\

We denote the
resulting homology by $\FF H(L_{0},L_{1};\eta)$. We will show in \S\ref{sec:fine_invariance} that it does
not depend on the choices made in its construction (except the choice of the orientations and relative spin structure of $L_0$ and $L_1$) and that it is invariant with
respect to Hamiltonian diffeomorphisms in the usual sense: if $\phi:M\to M$ is a Hamiltonian
diffeomorphism, then we have isomorphisms
$$\FF H(L_{0},L_{1};\eta)\simeq\FF H(\phi(L_{0}),L_{1};\eta')\simeq \FF H(L_{0},\phi(L_{1});\eta'')$$
where $\eta'$ and $\eta''$ correspond to $\eta$ via the Hamiltonian diffeomorphism.
As a corollary, we obtain that if $L_{0}$ and $L_{1}$ can be disjoined
by a Hamiltonian isotopy, then $\FF H(L_{0},L_{1};\eta)=0$ for all connected components
$\eta$.

\begin{remark}\label{rem:action}{\rm For any $\xi\in \bar\Ga_{\eta}$
we may define its action
$$\mathcal{A}(x)=-\int (u)^{\ast}\omega$$
where $u:[0,1]\to \Ga_{\eta}$ verifies $u(0)=\eta$, $u(1)=\pi(\xi)$ and belongs to the homotopy class defined by $\xi$. As mentioned before, by using an appropriate reparametrization,
each element of the moduli spaces $\Ww^{a}_{\ldots;b}(\la)$ may be viewed as a
strip joining $\bar a$ and $\bar b e^{\la}$ (recall that $\bar a, \bar b\in \bar\Ga_{\eta}$
are the fixed lifts of the elements $a,b\in I_{\eta}$). Clearly, this implies
that, if $v\in\Ww^{a}_{\ldots;b}\not=\emptyset$, then
$$\mathcal{A}(\bar{a})-\mathcal{A}(\bar{b}e^{\la})=\int v^{\ast}\omega$$ and this integral is
non-negative and, obviously, equal to the sum of the symplectic areas of the
strip and the disks appearing in $v$. Thus, as in the usual
Floer theory, the fine Floer complex admits an action filtration.
}\end{remark}

\begin{example}{\rm
Let $L_{0}=S$ be  a circle in $\C$ and $L_{1}=\R\subset \C$  so that
$L_{0}\cap L_{1}=a,b$. We intend to describe the fine Floer complex for
$L_{0},L_{1}$.
For this we fix a Morse function on $S$ with a single minimum, $m$, and
a single
maximum, $M$ both different from $a$ and $b$.
We take on $L_{1}$ a Morse function with no critical point.
The relevant Novikov ring is the group ring of an abelian group with
one
generator $e^{\la}$ with $\mu(\la)=2$ where $\la$ is the homotopy class
of the
disk whose boundary is the circle. This Novikov ring coincides with
that
appearing in the cluster complex of $S$. In this cluster complex we
recall that
$dM=0$ and $dm=(1+M+M^{2}+...)e^{\la}$. There are
two
pseudoholomorphic strips in the picture: one $u$ joining $a$ to $b$
and one $v$ joining $b$ to $a$, both are of Maslov class $1$.

We assume that $M$ belongs to the boundary of $v$
(if this is not the case, the resulting complex will be different even
if, of course,
the resulting homology will be the same).  Here is the
differential $d_F$ of the fine Floer complex. First, $d_Fa=ma+b$  - the
factor $ma$ corresponds
to the flow line along the circle joining the constant strip $a$ to the
minimum; the term
$b$ corresponds to the  strip $u$. Second,
$d_Fb=mb+(dm)a$ - the first term, $mb$, is as above; the term
$(dm)a=(1+M+M^{2}+...)ae^{\la}$ appears because the strip $v$ itself
passes through
$M$ ( $|M|=0$ makes it possible that all the powers of $M$ appear). It
is easy to check that
$(d_F)^{2}=0$ and that the resulting homology is trivial.}
\end{example}

\subsection{Symmetrization.}
\label{ss:cluster_fine_relation}
The fine Floer homology suffers from an algebraic defect:
each generator of $\FF C(L_{0},L_{1}, (f_{0},g_{0}), (f_{1},g_{1}))$, $x\in I_{\eta}$,
 has a differential of the form $dx=(m-m')x+\ldots$
where $m$ and $m'$ are, respectively, minima of $f_{0}$ and $f_{1}$. This makes sometimes the
fine Floer complex difficult to use. Moreover, this property implies that the
fine Floer homology is trivial in many situations. This happens for example, in the absence
of all bubbling, if $L_{0}$ is Hamiltonian isotopic to $L_{1}$ and the Morse functions
$f_{0}$ and $f_{1}$ are both prefect with a trivial Morse differential. In this
case the fine Floer differential is of the form $dx=(m-m')x$ for each generator $x$ and
the complex is therefore acyclic. However, in many interesting cases there is a natural, 
more symmetric version of the fine Floer complex
which is more efficient and easier to use. While various versions of this are possible
we will mostly use here one version which appears when $L_1$ is the image of $L_0$ by a Hamiltonian diffeomorphism.

\

We begin with the remark  that we could have
defined the fine Floer complex using the following, more general, setting: consider a generic time-dependent 
Hamiltonian $H:M\times [0,1] \to \R$, and denote by $\phi_{t \in [0,1]}$ its flow. 
Take a generic  family of almost complex structures $J_{t \in [0,1]}$ 
compatible with $\om$.

In the construction of the fine Floer complex, take as generators, instead 
of the intersection set $I_{\eta}(L_0, L_1)$, 
the set $I_{\eta}(L_0,L_1,H)$ of trajectories 
of $\phi_{t \in [0,1]}$ starting on $L_{0}$ and ending on $L_{1}$ 
in the component class of  some path $\eta$ from $L_{0}$ to $L_{1}$ (when $L_{0}=L_{1}$ and
no such choice is indicated, we work implicitly with $\eta$ equal to the constant path; clearly,
we may assume, generically, for each of these trajectories $\ga$ its $i$-end, $\ga(i)$, 
belongs to the unstable manifolds of the minimum of $f_{i}$, $i=0,1$).
Replace the condition requiring that the root be a $J$-holomorphic strip 
in the construction of the fine Floer complex 
by demanding that the root be a $(J,H)$- {\em semi-cylinder} joining 
two such trajectories $\ga$ and $\ga'$, i.e a map
$$
u: \R\times [0,1] \to M
$$
satisfying the equation
\begin{equation}\label{eq:fl_eq}
\partial_{s}u+J_t(u)\partial_{t}u+\nabla H(u,t)=0
\end{equation}
and the boundary conditions:

    1) $u(s,i) \in L_{i}$ for $i=0,1$  

     2) $\lim_{s \to - \infty} u(s,t) = \ga(t)$ and  $\lim_{s \to  \infty} u(s,t) = \ga'(t)$.

Note that the appropriate gluing of disks and
constant strips used in the proof of $d_{F}^{2}=0$ in
\S\ref{subsubsec:diff_fine1} remains valid in this context by replacing the constant strips by
constant semi-cylinders, i.e by maps $u: \R\times [0,1] \to M
$ that are of the form $u(s,t) = \ga(t)$ or $\ga'(t)$.  Besides these modifications, 
the construction of this  complex is identical to the fine Floer complex
described before: the only change  is that the root does not anymore correspond to a pseudo-holomorphic
strip but rather satisfies the non-homogenous equation (\ref{eq:fl_eq}). 

\   

The notation for the underlying moduli spaces
will be $$\Ww^{a}_{\ldots;b}((L_{0},f_{0}),(L_{1},f_{1}), J_{t \in [0,1]},H_{t \in [0,1]};\la)$$
or, if the context is clear, $\Ww^{a}_{\ldots;b}(\la)$. 

\

We will call the resulting complex and homology the {\em general fine Floer complex (homology)}.
Although we will not show it in this paper (because we will not need it), it is not difficult to see that the 
the proofs concerning $d_F^2= 0$ and the invariance of the fine Floer homology under variations of 
the auxiliairy data, that we give in \S~\ref{sec:fine_invariance}, may be generalized
to the general fine Floer homology (as long as the data remain generic). Now the point 
is that the general fine Floer complex admits two special cases (which, therefore, for generic data, 
 produce isomorphic homologies):
\begin{itemize}
\item[1)]  the first one, as seen before,  is when $H$ is the zero Hamiltonian and 
    this leads to the fine Floer complex.
\item[2)] the second one is when $L_0 = L_1$. 
\end{itemize}    

Let us discuss the second case and put $L=L_{0}=L_{1}$.  The generators $I_{\eta}(L_0,L_1,H)$ 
are now elements in the set $I_{\eta}(L,H)$ that consists of trajectories of $H_t$ starting 
and ending on $L$, while the configuration space that defines the differential is  
made of $(J_{t \in [0,1]}, H_{t \in [0,1]})$-semi-cylinders with the two ends 
coinciding with trajectories $\ga, \ga'$ and the two side boundaries 
lying on $L$. Without loss of transversality, we may choose the generic 
family $J_{t \in [0,1]}$ so that $J_0 = J_1$ (we will denote this almost complex structure
by $J$).  

We can now define a new complex, the {\em symmetric fine Floer complex}, by 
considering this setting in which, additionally, we choose  
the pair $(f_0,g_0)$ equal to the pair $(f_1,g_1)$ (we denote both pairs by
 $(f,g)$). 
Since we have
a differential graded algebra multiplication map:
$$
\Cl (L,(f,g),J)\otimes \Cl (L,(f,g),J)\to \Cl (L,(f,g),J)
$$
and because $J_0 = J_1 = J$, we may replace the ring 
$\mathcal{R} = \Cl(L_0,(f_0,g_0),J_{0})\otimes \Cl(L_1,(f_1,g_1),J_{1}) \otimes \bar{\La}$
that appears naturally in the definition of the fine Floer complex by the 
ring $\hat{\mathcal{R}}=\Cl (L,(f,g),J)\otimes _{\La}\bar\La$.
This leads to a fine Floer complex $(\hat\FF C(L,H,J, (f,g)), d_{\hat F})$ and a type of fine Floer homology
$\hat \FF H(L)$ which we call {\em symmetric}. If no additional notation appears the path $\eta$
used in this case is just the constant path. Regarding the algebraic 
problem mentioned at the beginning of this section,
note that if we choose $f$ will a single local minimum,
then the minima of $f_{0}$ and $f_{1}$ are
identified in $\hat{\mathcal{R}}$ so that the term $(m-m')x$ disappears from the differential
$d_{\hat F}(x)$. It is this algebraic identification of the clusters on the $L_0$-side with 
the clusters on the $L_1$-side that makes the symmetric fine Floer homology different, 
in general,  from the fine Floer homology. 

 This homology has the same type of invariance
properties as the non-symmetric version and, moreover, it is independent on the choice
of $(H, J, (f,g))$ as long as it remains generic.

\ 

In the remainder of the section we will show how to express the symmetric fine Floer
homology $\hat\FF H(L)$ just in terms of  the cluster complex.

\begin{remark}\label{rem:fine_triv} {\rm If $\Cl H_{\ast}(L)=0$, 
then $\hat{\FF} H_{\ast}(L)=0$. Indeed, the Proposition \ref{prop:free_terms} shows that 
$\Cl H_{\ast}(L)=0$ means that $1$ is a boundary in the cluster complex. As the 
fine Floer complex is a module over the cluster complex the claim follows. In other words,
if the cluster differential has free terms, then the symmetric fine Floer homology is trivial.
 This is also true
 of the non-symmetric version of the theory in which it is enough to assume that
just one $\Cl H_{\ast}(L_{i})$ is trivial to deduce that $\FF H_{\ast}(L_{0},L_{1})$ vanishes.
}\end{remark}

\subsubsection{An algebraic preliminary.} \label{subsubsec:alg_preliminary}

Consider a commutative, differential graded algebra of the form
$\mathcal{A}=(SV,d)$ with $V$ a rational vector space.

Consider the $SV$-module $SV\otimes V$ and write its elements
under the form $$x_{1}\ldots x_{k}\otimes v=x_{1}\ldots x_{k}\bar
v$$ where $v,x_{i}\in V$. Define the following linear map $$\alpha
: SV\to SV\otimes V$$ by letting $\alpha(v)=\bar{v}$, $\alpha
(1)=0$ where $1$ is the unit in $SV$ and extending this map by the
formula
$$\alpha(ab)=a\alpha(b) +(-1)^{|a||b|}b\alpha(a)~.~$$

It is easy to see by induction on the length of words that this
map is well defined and that the formula above is verified for any
homogenous monomials $a,b$. Explicitely, we have
$$
\al (x_{1}\ldots x_{k}) =  \sum_{i} (-1)^{\si_i} x_{1}\ldots\hat{x}_{i}\ldots
x_{k} \bar{x}_i.
$$
where $\si_i$ is the product of the degree of $x_i$ with the sum of the degrees of the $x_j, j > i$. 
We now define a map $d$ on $SV\otimes V$ as the unique
$(SV,d)$-module extension  of
$$ d\bar{v}=\alpha(dv)$$
which verifies the standard graded Leibniz rule. 

One easily checks (by induction on word length) that, for any
monomial $m$,  this map verifies $\alpha(dm)=d(\alpha(m))$. This
implies that the map $d$ so defined is a differential and that
$\alpha$ is a chain map. Denote by $$\widetilde{\mathcal{A}}=
(SV\otimes V, d)$$ the $\mathcal{A}$-differential module 
constructed in this way.

\subsubsection{String-strip symmetrization} \label{subsubsec:relation_cluster_fine}

We return to our geometric setting.  The algebraic construction
above appears in the next proposition.

\begin{proposition} \label{prop:cluster_fine}
The symmetric fine Floer homology verifies:

$$s^{-1}\hat\FF H_{\ast}(L)\simeq H_{\ast}(\widetilde{\Cl} (L,J_{0};(f,g)))~.~$$
\end{proposition}

See \S~\ref{se:relation_cluster_fine} for the proof of this
proposition.  It is essentially based on the following idea: suppose that 
in the definition of the symmetric fine Floer complex, instead of the 
Hamiltonian $H_{t \in [0,1]}$ and the loop $J_{t \in [0,1]}$, one takes 
the pull back to the Weinstein neighbourhood of $L$ of a generic Morse 
function $h$ on $L$ (extended conveniently to all of $M$), and a single generic $J$. 
Then the generators, as $\hat{\mathcal{R}}$-module, of this new symmetric 
fine Floer complex are the critical points of $h$. This leads to a 
complex $\hat{\mathcal{C}}(L,J,f,h,g)$, defined in detail 
in  \S~\ref{se:relation_cluster_fine} where we also show that 
its homology coincides with the symmetric fine Floer homology 
defined above. Finally, by letting $h$ be equal to $f$ 
in $\hat{\mathcal{C}}(L,J,f,h,g)$, one gets a complex 
that can be easily related, algebraically, to the 
cluster complex -- this is what produces the 
isomorphism of the above proposition (it is therefore 
produced in two steps: the first one is a PSS type of 
comparison map that relates cluster-strips to string-strips, 
i.e replace the semi-cylinder by  $(J,h)$-linear clusters; 
the second one is based on the fact that the identification 
of $h$ with $f$ leads to a symmetrization of the kind described 
in the above algebraic construction -- this is why it is natural 
to refer to this isomorphism as the {\em string-strip symmetrization}).

An immediate corollary of this proposition gives the description of the symmetric
Floer homology if no bubbling is present.
\begin{corollary} If $\omega|_{\pi_{2}(M,L)}=0$ then
$$\hat\FF H_{\ast}(L)\simeq  (S(s^{-1}H_{\ast}(L;\Q))\otimes\La)^{\wedge}\otimes H_{\ast}(L)~.~$$
\end{corollary}

Clearly, if a Lagrangian is displaceable, then both the fine Floer
homology and its symmetric version vanish.

\subsection{Applications.} \label{subsec:application}

 We will describe three consequences of our machinery.   

  \subsubsection{The Gromov-Sikorav problem.} \label{subsubsec:GS-problem}
  As a first consequence, we will study a plausible conjecture going
 back to Gromov's original paper \cite{Gr}
 on pseudo-holomorphic curves, stated orally by
 Sikorav in the late eighties in the following way:
 given any compact Lagrangian submanifold of $\C^n$,
 there is a holomorphic disc passing through each
 point of $L$ \footnote{The word ``conjecture'' should not be understood in the sense that 
 Sikorav claimed that this is true, but in the milder sense of asking the question and 
 noting that it is not a consequence of Gromov's or Floer's theories. Obviously, what 
 is at stake here, is the question of the degree of the evaluation maps on the boundaries of $J$-holomorphic disks.} 

    To state our first result in this direction, have in mind the notion of {\em free terms} and
  the statement of Proposition \ref{prop:free_terms}. 

\begin{corollary}\label{cor:uni_ruled}  Let $L \subset M$ be a compact,
orientable, relatively spin Lagrangian submanifold of any
symplectic manifold $M$. Assume that $\hat{\FF}
H_{\ast}(L)=0$ (for example if $L$ is displaceable by Hamiltonian isotopy). Then, for any generic
almost complex structure $J$ compatible with the symplectic form, one of the following holds:
\begin{itemize}
\item[i.]there are
$J$-holomorphic disks with boundary on $L$ passing through a dense subset  of points of $L$. 
\item[ii.] the
cluster complex $\Cl (L,J,f)$ has high free terms for some Morse
function $f$ with a single local minimum and a single local maximum (in particular, there are
$J$-disks of non-positive Maslov index).
\end{itemize}
\end{corollary}

\begin{proof} 
Assume that we are not in the case ii, i.e assume that for any
 Morse function with a single local minimum 
and local maximum, the associated cluster complex does not have high free terms.
Fix such a function $f$ and denote its minimum by $m$. 
By Proposition \ref{prop:cluster_fine},
$s^{-1}\hat\FF H_{\ast}(L)$ is isomorphic to $H_{\ast}(\widetilde{\Cl}(L,J;(f,g)))$, which means that 
$H_{\ast}(\widetilde{\Cl}(L,J;(f,g)))$ vanishes. Notice that if $d\bar{m}\not=0$, then
we also have $dm\not=0$ in the cluster complex. 
Assume now that $$\bar{m}\in\widetilde{\Cl}(L,J;(f,g))$$
is a cycle. We want to show that this also leads to $dm\not=0$.
Using Remark \ref{rem:cl_moduli}, we see that  $\bar{m}$ can be 
a boundary only if $\Cl(-)$ has free terms: indeed, by this remark, the only possible primitive of $\bar m$ must have the form $\tau \bar m$ 
where $\tau$ is a primitive of the unit $1$ in the cluster complex.
This means that there is a free term in some $dx$ for some $x\in\Crit(f)$. Once again by  Remark \ref{rem:cl_moduli}, the index of this $x$
cannot be $n$  and it cannot be stricly between $0$ 
and $n$ by our assumption. Therefore, $x=m$ and $dm\not=0$.

The fact that $dm\not=0$ means that there exists a moduli space
$^{\nu}\Mm^{m}_{x_{1},\ldots,x_{k}}(\la)$ of dimension $0$ and
non empty. But for a cluster tree to originate at the minimum $m$, the 
root disk must go through $m$. As we may use a different function $f$ to 
place $m$ in any generic point in $L$, this implies the claim.
\end{proof}

The dichotomy in the statement of the previous corollary can be sometimes
resolved by homological restrictions.

\begin{corollary} \label{cor:cond_uni_ruled} 
Suppose that $L$ is orientable, relatively spin and that $H_{2k}(L;\Q)=0$ for
$2k\not\in\{0,\dim(L)\}$.
If $\hat{\FF}H_{\ast}(L)=0$, then
$L$ verifies i. of the Corollary \ref{cor:uni_ruled} above.
\end{corollary}

\begin{proof} We postpone the full proof of this to \S\ref{subsec:min}. We only
give here the argument in the case when the manifold $L$ as above admits a perfect
Morse function. Obviously, all the critical points of such a perfect Morse function $f$
different from its minimum $m$ and from its maximum are of odd index. 
It follows that the cluster complex $\Cl (L,J,f)$ does not have any high free terms and, as
in the corollary above, this means that $dm\not=0$ in $\Cl(L,J,f)$ and implies that there
is a $J$-disk passing through $m$. By conjugating the given Morse function by a diffeomorphism
we may construct prefect Morse functions with the minimum in any desired point of $L$ which
shows the claim in this case.
\end{proof}  

\begin{remark}{\rm The statement in Corollary \ref{cor:cond_uni_ruled} 
is most relevant in the case when the minimal Maslov index of $J$-disks is 
non-positive. If this minimal Maslov index is greater
or equal than $2$, we shall see that the Corollary is valid without the orientability,
relative spin and homological assumptions by a simpler argument - see our Remark \ref{rem:sseq_comparison} c., as well as the end of \S\ref{subsec:bded_disks}.
In the monotone case a different proof is due to Peter Albers \cite{Al}.} 
\end{remark}

\begin{example}\label{ex:uni_rule}{\rm The homological restriction in the corollary above
 is serious but, still, there are many examples of such manifolds:
$S^{1}\times S^{n-1}$ and its connected sums with itself provides maybe the simplest 
examples. 
}\end{example} 

In \S\ref{subsec:bded_disks} we will improve these results in the displaceable case 
by also bounding from above the area of the disks detected in terms of the displacing energy. 
This upper bound and Gromov's compactness theorem then imply:

\begin{corollary}\label{cor:disks} Suppose that the relatively spin, orientable 
Lagrangian submanifold $L$ is displaceable by a Hamiltonian isotopy and let $E(L)$ be its Hofer displacement energy. Any $\omega$-tame almost complex structure $J$
has the property that one of the following is true:
\begin{itemize}
\item[i.]  for any point $x\in L$  there exists a
$J$-pseudo-holomorphic disk of symplectic area at most $E(L)$
whose boundary rests on $L$ and which passes through $x$.
\item[ii.] there exists a $J$-disk of Maslov index at most 
$$2-\min\{2k\in\N \backslash \{0, \dim(L) \} : H_{2k}(L;\Q)\not=0\}$$
and of symplectic area at most $E(L)$.
\end{itemize}
\end{corollary}

\

If the set $\{2k\in\N \backslash \{0, \dim(L) \} : H_{2k}(L;\Q)\not=0\}$ is empty, the statement should be understood as saying that (i) automatically  holds. 
In \S\ref{subsec:bded_disks} we 
also deduce an interesting geometric consequence of this.

\
\subsubsection{Constraints on Maslov indices.}
The cluster complex setting provides straighforward proofs of various
constraints regarding Maslov indices of Lagrangian submanifolds. 
For instance, we give in \S~\ref{subsec:calc} a rapid proof of Fukaya's 
recent result in \cite{Fu} that we may state in the following general 
form: if $S^1 \times S^{n-1}$ admits a Lagrangian embedding in a 
symplectic manifold in such a way that it is displaceable by a 
Hamiltonian isotopy, then $\{2,n\} \cap {\rm Im} (\mu) \neq \emptyset$ for 
$n$ even and $\{2, 3-n\} \cap {\rm Im} (\mu) \neq \emptyset$ for $n$ odd.

\

\subsubsection{Detection of periodic orbits.}
In \S\ref{subsec:orbits} we discuss another application. In the presence of an orientable
relatively spin pair of Lagrangian submanifolds $L$ $L'$, we show that, by replacing
in the definition of the clustered moduli spaces anyone (one and only one) of the $J$-disks by a pseudo-holomorphic
cylinder with one boundary on $L$ and the other boundary on $L'$,
one can construct a chain map:

$$\cyl: \Cl C(L) \to \Cl C(L) \otimes \Cl C(L') \otimes \La_{\Phi_0}$$
where  $\La_{\Phi_{0}}$ is an appropriate Novikov ring.

This map induces a morphism in homology whose non-triviality is used to detect
the existence of periodic orbits of Hamiltonian diffeomorphisms  that {\em separate} 
$L$ from $L'$. These results generalize considerably those in \cite{GL}.

\section{The cluster complex and its homology.}\label{sec:cluster_cplx_hom}

We now discuss in more detail the construction of the cluster complex - in particular that
of the underlying moduli spaces. We then pursue this section with the proof of the
invariance of the cluster homology and with some other of its properties.

\subsection{Clustered moduli spaces} \label{sec:cluster_moduli}

As we mentionned in \S~\ref{subsec:cluster_definition}, the definition of the cluster differential
$$d:\Cl(L,J;(f,g))\to \Cl(L,J;(f,g))$$
depends on certain
moduli spaces which combine  Morse theory with the main features of  the standard
construction of moduli spaces of stable maps. We gave an idea of what these spaces
look like in \S~\ref{subsec:cluster_definition}.
We now describe them in detail.

\

We fix a compact oriented Lagrangian submanifold $L\subset M$ as well as an almost
complex structure $J$ compatible with $\omega$ and we
also fix $f:L\to \R$ a Morse function and $g$ a Riemannian metric
on $L$ so that the pair $(f,g)$ be Morse-Smale. Finally, we fix as well a relative
spin structure on $L$.

\

\subsubsection{Trees.}We start with the notion of $\n$-labeled, oriented,
metric trees where $\n=(n_{1},n_{2})\in \N\times\N$.

These are couples $(\mathcal{T},\psi)$ where:
$\mathcal{T}=(T,V_{D}, \phi_{T},\mu_{T})$ with $T$ a (connected)
tree with oriented edges and with exactly one edge entering each
vertex except for one - the root of the tree $v_{0}\in Vertex(T)$
- in which no edge enters; $V_{D}\subset Vertex(T)$ is a subset of
the set of vertices of $T$ which generates a connected subtree of
$T$ and $v_{0}\in V_{D}$; $\phi_{T}:\{0,1,\ldots, n_{1}\}\to
V_{D}$, $\phi_{T}:\{n_{1}+1,\ldots , n_{1}+n_{2}\}\to
Vertex(T)\backslash V_{D}$ assigns markers to the vertices of $T$
so that $0\in \phi_{T}^{-1}(v_{0})$; $\mu_{T}:Vertex(T)\to \La$
assigns classes in $\pi_{2}(M,L)/\sim$ to each vertex of $T$;
$\psi :Edges(T)\to [0,+\infty)$ assigns lengths to the edges of
$T$ so that if $e\in Edges(T)$ relates two vertices so that one of which (at least)
does not belong to $V_{D}$, then $\psi_{T}(e)=0$.
Below we denote $V_{S}=Vertex(T)\backslash V_{D}$.

\

We denote the fact that $\alpha$ is related by an edge to $\beta$
by $\alpha E \beta$. By a slight abuse of notation we shall denote
the edge joining $\alpha$ to $\beta$ also by $\alpha E \beta$. For
further use, we let for each $\alpha\in Vertex(T)$
$$n_{\alpha}=\#(\phi^{-1}_{T}(\alpha))+\#\{\ \beta\in Vertex(T) :
\alpha E \beta \ \rm{or} \  \beta E \alpha\ \}.~$$

\subsubsection{Pseudo-holomorphic disks and flow lines}\label{subsec:disks_flow}
For a Morse-Smale pair $(f,g)$, we denote by $\gamma_{t}(-)$ the
negative gradient flow induced by $-\nabla f$. We understand
$\gamma_{\infty}(x)$ as $\lim_{t\to\infty}\gamma_{t}(x)$.

For our fixed Lagrangian submanifold $L$, almost complex structure $J$ and
Morse-Smale pair $(f,g)$, we define the clustered moduli space
modeled on the metric tree $(\mathcal{T},\psi)$:
$$\mathcal{M}_{\mathcal{T},\psi}(L,J,f)=
\{[(\mathbf{u},\mathbf{z})] :(\mathbf{u},\mathbf{z})=
(\{u_{\alpha}\}_{\alpha\in Vertex(T)},
(\{z_{\alpha\beta}\}_{\alpha E \beta\ \rm{or}\ \beta E\alpha},
\{z_{i}\}_{1\leq i\leq n_{1}+n_{2}}))\}$$ where the pairs
$(\mathbf{u},\mathbf{z})$ are so that:
\begin{itemize}
\item[i.]$u_{\alpha}:D^{2}\to M$ is a $J$-holomorphic disk with
boundary on $L$ and with Maslov class $\mu(\alpha)$ when
$\alpha\in V_{D}$. \item[ii.] when $\alpha\in V_{S}$,
$u_{\alpha}:S^{2}\to M$ is a pseudoholomorphic sphere with Chern
class $2\mu(\alpha)$ .
 \item[iii.] If $\alpha \in V_{D}$, then
$z_{\alpha\beta}$ are points in $D^{2}$ with the property that if
$\beta\in V_{D}$, $\alpha E \beta$, then $z_{\alpha\beta}\in
S^{1}$, $z_{\beta\alpha}\in S^{1}$ and $\gamma_{\psi(\alpha
E\beta)}(u_{\alpha}(z_{\alpha\beta}))=u_{\beta}(z_{\beta\alpha})~.~$
If $\alpha\in V_{S}$, then $z_{\alpha\beta}\in S^{2}$ and if
$\beta\in Vertex(T)$ so that $\alpha E\beta$ or $\beta E\alpha$,
then $u_{\alpha}(z_{\alpha\beta})=u_{\beta}(z_{\beta\alpha})~.~$

\item[iv.] We have $z_{i}\in S^{1}=\partial D^{2}$ if $\phi_{T}(i)\in V_{D}$
and $z_{i}\in S^{2}$ if $\phi_{T}(i)\in V_{S}$. For $\alpha\in
Vertex(T)$ fixed, the points $z_{\alpha\beta}$, $z_{\delta\alpha}$
and $z_{i}, \ i\in \phi_{T}^{-1}(\alpha)$ are pairwise distinct.
The union of all these points is the set of {\em special points}
on the disk (or sphere) $u_{\alpha}$. Their number is equal to
$n_{\alpha}$.

 \item[v.] If $u_{\alpha}$ is the constant map,
 then $n_{\alpha}\geq 3$ and there exists a chain of edges of  length $0$
 connecting $\alpha$ to a vertex $\beta$ with $u_{\beta}$ non-constant.
 \end{itemize}

The points $z_{\alpha\beta}$ will be called the {\em incidence
points} corresponding to the disk (or sphere) $\alpha$ and the
points $z\in \phi_{T}^{-1}(\alpha)$ are the {\em
marked points}.

\

The obvious reparametrization group acts on these objects and
 $[(\mathbf{u},\mathbf{z})]$ are the equivalence
 classes obtained by taking the quotient with respect to this
 action.  We call such an equivalence class a {\em cluster configuration}
or, sometimes, a cluster element or cluster tree.

We now define other associated moduli spaces
$$\mathcal{M}_{\mathcal{T}}(L,J,f)=
\bigcup_{\psi}\mathcal{M}_{\mathcal{T},\psi}(L,J,f)$$ and
$$\mathcal{M}(L,J,f;\la,\n)=\bigcup_{\mathcal{T}\in\mathcal{T}(\n,\la)}\mathcal{M}_{\mathcal{T}}(L,J,f)~.~$$
 where the set $\mathcal{T}(\n,\la)$ consists of all trees $\mathcal{T}$ so that: $$\#(Dom(\phi_{T}))=n_{1}+n_{2} \ , \ \sum_{\alpha\in Vertex(T)} \mu_{T}(\alpha)=\la~.~$$
The moduli spaces $\mathcal{M}_{\mathcal{T};\psi}(L,J,f)$ are
endowed with the
 Gromov topology and this induces topologies on the
 other moduli spaces introduced above (except for the minor modifications that we
 introduced, this is by now a standard construction;
 for an exposition see, for example,
 McDuff-Salamon \cite{McSa}). Using this topology,
 notice also that the last union above is certainly not disjoint.
 Indeed, if the length of an edge joining the vertices $a$ and $b$
 in a tree shrinks to zero,  the corresponding cluster
 element is also obtained by bubbling from
 an element for which the underlying tree has a single vertex
 instead of $a$, $b$ (and the edge $aEb$ is absent).

\subsubsection{Compactification} \label{subsec:compactification}
 We now discuss the compactification of
 $\mathcal{M}(L,J,f;\la,\n)$, $\bar\Mm(L,J,f;\la,\n)$.
 We claim that
 \begin{equation}\label{eq:pre_bdry}  \begin{array}{l}
 \bar{\mathcal{M}}(L,J,f;\la,\n)\backslash
 \mathcal{M}(L,J,f;\la,\n)=\\[1ex]
 =\bigcup_{x\in \Crit(f)}\bar{\mathcal{M}}_{x}(L,J,f;\la',\n')
 \times \bar{\mathcal{M}}^{x}(L,J,f;\la'',\n'')
 \end{array}
 \end{equation}
  where
 $\la'+\la''=\la$, $\n'+\n''=\n$ and the space
 $\mathcal{M}_{x}(L,J,f;\la',\n')$ is defined similarly to
 $\mathcal{M}(L,J,f;\la,\n)$ except that we modify the definition of
 $\mathcal{M}_{T,\psi}(L,J,f)$ in the following way: there exists
 $\alpha\in V_{D}$ and one special point $z_{x}\in
 S^{1}$ so that
 $\lim_{t\to\infty}\gamma_{t}(u_{\alpha}(z_{x}))=x$. The
 definition of $\mathcal{M}^{x}(L,J,f;\la'',\n'')$ is analogous
 with the only change that $\lim_{t\to
 -\infty}\gamma_{t}(u_{v_{0}}(z_{0}))=x$ (recall that $v_{0}$ is the root of the
 tree).

 This formula is immediate as the only possible reason for non-compactness of
 the space $\mathcal{M}(L,J,f;\la,\n)$ is due to the fact that the length of
 some gradient flow line of $f$ may tend towards $\infty$.

 \begin{remark}{\rm
The role of the flow lines that appear in the definition of our
moduli spaces  becomes apparent from this product formula as it
implies that the bubbling of disks is an {\em internal}
``codimension one" phenomenon here. For example, a configuration
of two disks joined in a point may occur as the limit of bubbling.
But then it is also the limit of a configuration formed by two
disks connected by a flow line whose length decreases to $0$.}
 \end{remark}

Some related moduli spaces play a key role in the construction of
the differential of the cluster complex.

\

For $x,x_{1},\ldots,x_{k}\in \Crit(f)$ we let
$\mathcal{M}^{x}_{x_{1},\ldots,x_{k}}(L,J,f;\la,\n)$ be defined in
the same way as above except that
$\mathcal{M}_{\mathcal{T},\psi}(L,J,f)$ is replaced with
$(\mathcal{M}^{x}_{x_{1},\ldots,x_{k}})_{\mathcal{T},\psi}(L,J,f)$
where the definition of this last moduli space is the same as that
of the moduli space $\mathcal{M}_{\mathcal{T},\psi}(L,J,f)$ with
one additional condition:

\

\begin{itemize}
\item[vi.] We have $\lim_{t\to
-\infty}\gamma_{t}u_{v_{0}}(z_{0})=x$ (we recall that $v_{0}$ is
the root of the tree $T$ and that $0\in\phi_{T}^{-1}(v_{0})$) and
there is an ordered set of marked points $(z_{1},\ldots,z_{k})$,
$\phi_{T}(i)\in V_{D}$ called {\em  terminations} so that
$$\lim_{t\to +\infty}\gamma_{t}(u_{\phi_{T}(j)}(z_{j}))=x_{j}~.~$$
\end{itemize}

This is clearly an extension of the moduli spaces
$\mathcal{M}_{x}(L,J,f;\la,\n)$, $\mathcal{M}^{x}(L,J,f;\la,\n)$
that appeared above.

When $\mathcal{T}$ is the void graph we let
$(\Mm^{x}_{x_{1}})_{\emptyset}(L,J,f)$ be the set of negative
gradient flow lines joining the critical points $x$ and $x_{1}$.

\

Notice that, for each permutation $\sigma\in \Sigma_{k}$ there is
an obvious homeomorphism:

$$\xi_{\sigma}:\Mm^{x}_{x_{1},\ldots,x_{k}}(\la,\mathbf{n})\to
\Mm^{x}_{x_{\sigma_{1}},\ldots,x_{\sigma_{k}}}(\la,\mathbf{n})$$
which is induced by permuting the indexes of the terminations
$z_{i}$ (in vi. above) according to $\sigma$.

\begin{remark} {\rm For the existence of these homeomorphisms, it is
significant that the trees used in our construction are not required
to be planar (for a planar tree, the terminations
inherit a natural order and not each permutation of them can be extended
to a homeomorphism of the given planar tree to another planar tree)}.
\end{remark}

\subsubsection{Justification for equation (\ref{eq:dsquare1})}
As mentioned in \S\ref{subsec:cluster_definition}, formula (\ref{eq:dsquare1})
together with the analysis of orientations imply the vanishing of the
square of the cluster differential. In this subsection, we justify this
formula and, in the next one, we discuss orientations.

We shall assume from now on that {\em the critical points of $f$ are
(strictly) ordered} (this is useful because later on we will have to
 take orientations into account and the maps $\xi$ do
not necessarily preserve orientation).

\

Let $S=\{x_{1},\ldots, x_{k}\}$ be a set consisting of critical
points of $\Crit(f)$ (repetitions being allowed) so that $S$ is a
partially ordered set whose elements respect the order fixed on
$\Crit(f)$ as in \S~\ref{subsec:cluster_definition}.

We now define:
$$\Mm^{x}_{S}(\la;\n)=\Mm^{x}_{x_{1},\ldots,x_{k}}(\la,\n)~.~$$

\

These moduli spaces can again be topologized by using the Gromov
topology and their compactifications verify the following product
formula :

\begin{equation}\label{eq:bdry}\begin{array}{l}
\bar{\Mm}^{x}_{S}(\la,\n)\backslash \Mm^{x}_{S}(\la,\n)= \\[1ex]
 \bigcup_{S=<S'\cup S''>,y, \la'+\la''=\la, \n'+\n''=\n}
 \bar\Mm^x_{<S',y>}(\la',\n')\times \bar\Mm^y_{S''}(\la'',\n'')
\end{array}\end{equation}

\

To see that this formula is correct, first notice that bubbling in a family
of elements in $\Mm^{x}_{S}$ leads to an element which is still
in $\Mm^{x}_{S}$ (this would not necessarily be so if, to define our
moduli spaces, we would have used planar graphs together with the
resulting ordering of the terminations). Moreover, use
(\ref{eq:pre_bdry}) and notice that if one flow line corresponding
to an edge $\alpha E \beta$ breaks, then the resulting element
is represented by two trees $T'$ and $T''$ so that $T''$ is the
full subtree generated by $\beta$ and $T'$ is the tree obtained by
eliminating from $T$ all the vertices of $T''$. There only remains to see
that the orientations in this formula agree -- see the next subsection for this.

Finally, we may also define unparametrized moduli spaces
denoted by
$$\Mm^{x}_{x_{1},\ldots,x_{k}}(L,J;f,\la)$$ (and similarly for the
other moduli spaces involved - the difference in notation is that
we do not indicate the marked points anymore). These are subspaces
in $\cup_{\n}\Mm^{x}_{x_{1},\ldots,x_{k}}(\la,\n)$ defined by the
additional condition that $\n=(n_{1},n_{2})$ verifies $n_{1}=k+1$,
$n_{2}=0$.

Our moduli spaces may be viewed as intersections  of moduli spaces of
disks and appropriate moduli spaces of flow lines by means of evaluation maps
at the special points. Under ideal transversality conditions, one expects that these
moduli spaces have the structure of singular manifolds with
boundary and that their dimension is given by formula (\ref{eq:dim}).
Here a singular closed manifold of dimension $n$ is a stratified
space $A$ with two strata $A=A_{1}\cup A_{2}$ where $A_{1}$ is
an open, dense  manifold of dimension $n$ while $A_{2}$ is a closed
stratified space of dimension at most $n-2$. The
stratification condition means that there is a neighbourhood of
$A_{2}$ in $A$ which admits $A_{2}$ as a deformation retract. This
implies, in particular, that $A$ admits a fundamental class.
Similarly, a singular manifold with boundary $A$ is a stratified space
with three strata, $A=A_{1}\cup A_{2}\cup A_{3}$, where $A_{1}$,
$A_{2}$ are as above and $A_{3}$ is a closed singular manifold of
dimension $n-1$ so that $\partial A\backslash A_{2}=A_{3}$. We
will denote $\partial A=A_{3}$.

It is well-known that the  idea transversality conditions
mentioned above are not in fact satisfied even for generic choices
of $J$ (this is seen in our case in particular when two terminations cross - this
situation also indicates the need for ghost bubbles in our model).
 However, following the methods recently
developped by Hofer, Wisocki and Zehnder in \cite{HoWiZe}, it is possible to
construct a system of perturbations $\nu$ so that for generic $J$ and
$(f,g)$, enough transversality is achieved to
obtain, at least,  a moduli space carrying  a rational cycle. More precisely we
have:

\

 {\bf Basic analytic result.} \label{thm:basic-analytic-intro}
\emph{For  generic choices of $(f,g)$ and $J$ which is compatible with $\omega$ and for each
 $\la \in \pi_2(M,L)/ \sim$, $k \ge 0$ and $x, x_1, \ldots x_k \in \Crit(f)$ there is
 a system of perturbations:
$$
\nu = \nu(f,g,J, x,x_1, \ldots,x_k, \la)
$$
of the equation $\bar{\p}_{J} = 0 $ on each pseudoholomorphic disk of the tree $\Tt'$, such that
each moduli space $\Nu \Mm^x_{x_1, \ldots, x_k} (\la)$, consisting by definition
of the configurations satisfying the perturbed equations $\bar\partial_{J}(u)=\nu(u)$
admits a natural compactification $\Nu \bar{\Mm}^x_{x_1, \ldots, x_k} (\la)$ that carries
the structure of an oriented rational singular cycle with boundary
of dimension
$$
|x|-\sum_{i=1}^{k}|x_{i}|+\mu(\la)-1 ~.~$$}

\

As mentioned in the introduction, the proof of this regularity
property of the perturbed, clustered moduli spaces, is postponed
to a subsequent paper.
Variants of this result will be used below in the description of
other moduli spaces similar to the clustered ones.

\

If $|x|-\sum
|x_{i}|+\mu(\la)=2$, then, as a consequence of (\ref{eq:bdry})
 we deduce formula (\ref{eq:dsquare1}).

\subsubsection{Orientation.}\label{subsec:orientation}

We will show in this section that, although the orientation is not the one that one would expect a priori -- see Definition~\ref{def:orientations} -- the fact that $L$ is orientable and relatively spin implies that we may assign
coherent orientations to our clustered moduli spaces. This will be done in three steps:

\

1)  It follows from \cite{FOOO} that the moduli space of
perturbed pseudo-holomorphic disks in a given class $\la$ inherits a natural orientation associated to a choice of (i) an orientation of $L$, 
(ii) the choice of an extension $st \in H^2(M; \Z)$ of the second Stiefel-Withney class of $L$, and (iii) the choice of a relative spin structure on $(M, L)$ (which is made possible by the last two conditions on $L$). 

2)  This gives an orientation on the ``root'' cells of our moduli spaces $$\Nu \Mm(L,J,(f,g); \la, (n_1,n_2) = (0,0))$$ without marked points, corresponding to the trees consisting of a single vertex; we then show how this can be used to orient, in a coherent way, all other cells of any codimension in $\Nu \Mm(L,J,(f,g); \la, (n_1, n_2)=(0,0))$ (i.e correponding to any other tree).

3) Given coherent orientations on $\Nu \Mm(L,J,(f,g); \la, {\bf n} = (0,0))$, we get orientations on $\Nu \Mm(L,J,(f,g); \la, {\bf n})$ for general $\bf n$, and finally on the moduli spaces $\Nu \Mm^x_{x_1, \ldots, x_k} (L,J,(f,g); \la)$ (see Definition~\ref{def:orientations})

\

   We will then look at the compactifications of these moduli spaces and prove that, with these orientations, the formula $d^2=0$ for the cluster differential is valid.

     For the convenience of the reader, we are going to briefly review the relevant arguments from \cite{FOOO} in (1) above. So let us assume that an orientation of $L$ as well as an extension $st$ of the second Stiefel-Withney class have been chosen, and let us define the orientation on the root cell of $\Nu \Mm(L,J,(f,g); \la, {\bf n} = (0,0))$, which are just maps $u$ satisfying the equation $\bar{\partial}_J = \nu$ for generic choices of $J$ and $\nu$. More precisely, denote by $\Bb$ the space of all $W^{1,p}$-maps $u: (D^2, \partial D^2) \to (M,L)$ with fixed integer $p > 2$ and by $\Ee$ the infinite dimensional vector space over $\Bb$ whose fiber $F_u$ at $u$ is the space of $(0,p)$-sections  $\Ga^{0,p}(\La^{0,1} D^2 \otimes u^*TM )$. The operator $\bar{\partial}_J$ is a section of this bundle. For a small generic section $\nu$ of this bundle, our moduli space is the inverse image of $\nu$ by the  $\bar{\partial}_J$-section. In order to linearize this operator at a solution $u$, one needs a Hermitian connection on this bundle. Indeed, if $\Hh$ is such a connection, viewed as a horizontal distribution, the linearisation is given by definition by:
$$
D^{\Hh}_u (\bar{\partial}_J) = \Pi_{\Hh} \circ d_u (\bar{\partial}_J): \Ga^{1,p}(D^2, \partial D^2; u^*TM, u^*TL) \to \Ga^{0,p}(\La^{0,1} D^2 \otimes u^*TM )
$$
where
$$
d_u : T_u \Bb = \Ga^{1,p}(D^2, \partial D^2; u^*TM, u^*TL) \to T_{\bar{\partial}_J(u)} \Ee
$$
is the ordinary differential and $\Pi_{\Hh}$ is the linear projection, induced 
by the connection, on the tangent space to the fiber at $ \bar{\partial}_J(u)$, identified 
with the fiber itself. Such a connection is most naturally induced by a connection on the ambient 
manifold $M$. Let $ \nabla$ be the Levi-Civita connection associated to the 
metric $g( \cdot, \cdot) = \om (\cdot, J \cdot)$ and define, as in \cite{McSa}, the following connection  
$$
 \tilde{\nabla}_v X = \nabla_v X - \frac{1}{2} J (\nabla_v J) X.
 $$
 The point is that, with respect to this new connection, $J$ is parallel (but not necessarily $g$ anymore). 
 (Recall that if $J$ is not integrable, the Levi-Civita connection need not preserve $J$). 
 Then, as showed in \cite{McSa}, the explicit expression for $D^{\Hh}_u \bar{\partial}_J$ is:
 $$
D^{\Hh}_u \bar{\partial}_J (\xi) = (\tilde{\nabla} \xi)^{0,1}  + \frac{1}{4} N_J ( \xi, \partial_J (u)) 
$$
where the first term is the complex-anti-linear part of $\tilde{\nabla} \xi$, equal by definition to 
$$
\frac{1}{2} ( \tilde{\nabla} \xi + J(u) \circ (\tilde{\nabla} \xi ) \circ j_{D^2}
$$
(here  $ j_{D^2}$ is the standard complex structure on $D^2$ equal to the product by $\sqrt{-1}$), and where in the second term, $N$ is the integrability tensor of $J$: $N(X,Y) = [JX,JY] [X,Y] - J[JX,Y] - J[X,JY]$. It is easy to check that the first term, of order $1$, of $D^{\Hh}_u \bar{\partial}_J (\xi)$ is the $(j_{D^2}, J(u))$-complex linear part while the second term, of order $0$, is the $(j_{D^2}, J(u))$-anti-complex linear part.  By hypothesis, we may assume that  the differential $D_u = D^{\Hh}_u  \bar{\partial}_J - D^{\Hh}_u \nu$ is onto, where of course the new term $D^{\Hh}_u \nu$ is of order $0$.

   The argument for the point (1) above breaks in two parts: one first shows that the choice of a trivialization of the real oriented bundle $u^*(TL) \to \p D^2$ naturally induces an orientation on the tangent space of $\Nu \Mm(\la) =    \Nu \Mm(L,J,(f,g); \la, {\bf n} = (0,0))$ at $u$, i.e an orientation on the kernel of $D_u$, and one then shows that the choices made of an orientation and a relative spin structure on $L$ enables one to define the trivializations coherently along any loop in the moduli space. Let us explain the first part. Pick a $SO(n)$-trivialisation 
$$
\Phi:  \p D^2 \times \R^n \to u^*(T_*L)|_{\p D^2} 
$$
and consider the $U(n)$-trivialization 
$$
\Phi^{\C}:  \p D^2 \times \C^n \to u^*(T_*M)|_{\p D^2} 
$$ 
that it induces by complexification (using $J(u)$).  
Let us denote by $E \to D^2$ the bundle obtained by attaching 
 $u^*(TM) \to D^2$ to the bundle $\p D^2 \times \C^n$ via the map $\Phi^{\C}$. 
 Because the above Hermitian connection is flat on the boundary, one may assume 
 that both the bundle $E$ and the connection have been trivialized on a closed 
 collar neighbourhoud $N$ of the boundary of the disk whose inner circle will be 
 denoted by $C$. By shrinking $C$ to a point, and using standard gluing techniques 
 (see Appendix A in \cite{McSa}), one gets an isomorphism
$$
\ker D_u \simeq \ker(Hol_{\nu,J}(D^2, \p D^2: \C^n, \R^n) \times Hol_{\nu,J}(S^2; E'_{tr}) \stackrel{\to}{\ev} \C^n )
$$
where the first term is the space of $(\nu,J)$-holomorphic maps from   
$(D^2, \p D^2)$ to $(\C^n, \R^n)$, the second term is the space of $(\nu,J)$-holomorphic 
sections of the linearized operator in the absolute case (that is to say defined on a 
closed Riemann surface) of the resulting Hermitian fiber bundle $E'$ in which 
the fiber at the point $p \in S^2$ is  identified with $\C^n$ though $\Phi^{\C}$, 
and where the evaluation map sends $(v,w)$ to   $v(0) - (\Phi^{\C})^{-1}(w(p))$. 
But the first space consists only of constant maps, and can be identified 
with $\R^n$ while the second space carries a canonical orientation (see \cite{McSa} for this). 
Therefore the right hand side is oriented, and the left hand side inherits this orientation. 

\begin{remark}{\rm For the convenience of the reader, let us mention that one can see this canonical orientation of the second space by considering the above equation $D_u = D^{\Hh}_u \bar{\partial}_J - D^{\Hh}_u \nu$ which applies as well to the case of the sphere. Since the order $1$ term is $(j_{D^2}, J(u))$-complex linear, one may deform the operator, within the space of Fredholm operators, to the operator consisting of the first term only, which then carries a natural complex structure (being the kernel of a complex map) and therefore a natural orientation. But this orientation can then be carried back, along the one-parameter family of Fredholm operators, to an orientation of the original kernel. Note that along this generic deformation, one may assume that only real dimension one cokernels may appear; in any case, it is an easy exercise to check that the orientation on the determinant bundle $\ker \otimes {\rm coker}$ can be carried over the family.}
\end{remark}

  This shows how to orient the tangent space of the root cell of 
  $$
  \Nu \Mm(\la) =: \Nu \Mm(L,J,(f,g); \la, {\bf n} = (0,0))
  $$ 
  at a given point $u_0$. To prove that this defines an orientation of the the root cell, we now come back to the argument in \cite{FOOO}. Starting from a given map $u_0$, we may transport the orientation at $u_0$ along any path to an orientation at some other map $u_1$ by noting that a trivialisation over $\p D^2$ extends canonically to a trivialisation over $\p D^2 \times [0,1]$.  To show that the orientation defined above at two different maps $u$ and $u'$ does not depend on the choice of the path, we must show that this orientation carries back to itself along any loop $u_{t \in [0,1]}$ based at some map $u = u_0 = u_1$. 
 Assuming  (i) and (ii) of the above point (1), one first triangulate $L$ and extend this triangulation to $M$. Let $V$ be an oriented real fiber bundle over the $3$-skeleton $M^{(3)}$ of $M$ such that $w_2(V) = st |_{M^{(3)}}$. By the choice of $V$, the bundle $TL \oplus V \to L^{(3)}$ has both of his first and second Stiefel-Withney classes equal to zero, and is therefore spin. We can thus make the choice of  a spin structure on $TL \oplus V \to L^{(3)}$, i.e the choice of a trivialization of $TL \oplus V$ over $L^{(2)}$. If $u_t$ is a loop based at $u$, consider the map
 $$
 U: (D^2, \p D^2) \times S^1 \to (M,L)
 $$
 sending $(x,t)$ to $u_t(x)$, which can be homotoped to $(M{(3)}, L{(2)})$. 
 The spin structure gives a trivialization of $TL \oplus V$ over $L{(2)}$, i.e a 
 trivialization of $U^*(TL \oplus V) \to \p D^2 \times S^1$. But   $U^*(TL \oplus V) = U^*(TL) \oplus U^*(V)$ 
 and the second term extends to $D^2 \times S^1 \simeq S^1$ and 
 is therefore trivialized (because $V$ is oriented). Thus, this gives
 a stable trivialization of $u_t^*(TL)$ that varies continuously over time. 
 But over $S^1$, the $SO(n)$-trivializations coincide with the stable $SO(n)$-trivialisations. 
 This completes the proof of 1) above.

 (2)   The orientation on the root cell of $\Nu \Mm(\la) = \Nu \Mm(L,J,(f,g); \la, {\bf n} = (0,0))$ gives rise in the obvious way to an orientation on the same moduli space with $m$ distinct marked points that we will denote $\Nu \Mm_k(\la)$.   Given a partition of $\la$ into a sum $\la_1 + \ldots + \la_m$, and a tree $T$ with $m$ vertices and  incidence points $z_{\al \be} \in \p D^2$, we consider, for each $1 \le \al  \le m$,  the moduli space $\Nu \Mm_{\al}(\la_{\al})$ equal to $\Nu \Mm_{m_{\al}}(\la_{\al})$ where $m_{\al}$ is the valence of the vertex $\al$ and $\la_{\al}$ the class corresponding to the vertex $\al$. Consider the evaluation map:
 $$
\Phi: ( \Pi_{\al}  \;  \Nu \Mm_{\al}(\la_{\al}) ) \times (\R^+ \times \ldots \times \R^+) \to L^2 \times \ldots \times L^2
$$
where the number of copies of $\R^+ = [0, \infty)$, as well as the number of copies of $L^2$, is equal to the number of edges in $T$ and where, for each edge $\al E \be$ with $\la < \be$, $\Phi$ maps 
$$
((u_{\al}, z_{\al \be}), (u_{\be}, z_{\be \al}), t) \in \Nu  \Mm_{\al}(\la_{\al}) \times \Nu \Mm_{\be}(\la_{\be}) \times \R^+ 
$$
 to $(\phi_t(u_{\al}(z_{\al \be})), u_{\be}(  z_{ \be \al} )) \in L^2$. Clearly, 
 the parametrised moduli space $$\Nu \Mm_{T, \psi} (L,J,(f,g))$$ is equal to the 
 inverse image by $\Phi$ of the products of diagonals $\De \times \ldots \De \subset  L^2 \times \ldots \times L^2$.  
 Since $L$ is oriented, $\De$ is co-oriented in $L^2$ and this induces a co-orientation, 
 and therefore an orientation, on $\Nu \Mm_{T, \psi} (L,J,(f,g))$. 

    We must now show that these orientations are compatible at the codimension 0 strata of 
    the boundaries of the compactifications of $\Nu \Mm_{T, \psi} (L,J,(f,g))$. Generically, 
    the only case to consider corresponds to the elementary surgery on a tree $T$ in 
    which an edge $\al E \be$ is contracted to a point while the two vertices $\al, \be$ are 
    identified. Let's denote by $T'$ the tree after this surgery. This corresponds to 
    the gluing of two $(\nu, J)$-holomorphic disks, while of course the reverse 
    surgery corresponds to the bubbling-off of a given disc. Denote by $\Nu \Mm_{T, \psi, \al \be} (L,J,(f,g))$ 
     the codimension $1$ stratum of $\Nu \Mm_{T, \psi} (L,J,(f,g))$ formed of 
     configurations in which the time-flow $\psi(\al E \be)$ is zero.
    Consider the corresponding attaching map
$$
\Nu \Mm_{T, \psi, \al \be} (L,J,(f,g))   \to \Nu \bar{\Mm}_{T', \psi'} (L,J,(f,g))
$$
 given by gluing. Equipping the domain with the orientation induced as the boundary of $\Nu \Mm_{T, \psi} (L,J,(f,g)) $ (which of course is simply the one given by the inverse image of the above map $\Phi$ with the modification that $\Phi$ now maps the factor $((u_{\al}, z_{\al \be}), (u_{\be}, z_{\be \al}), t) \in \Nu \Mm_{\al}(\la_{\al}) \times \Nu \Mm_{\be}(\la_{\be})  $
 to $(u_{\al}(z_{\al \be}), u_{\be}(  z_{ \be \al} )) \in L^2$), and equipping the codomain by the natural orientation, we must show that this attaching map preserves the orientations. It is evident that this reduces to prove the following statement.

 \begin{lemma} Let $ \Mm(\la_1, \la_2)$ denote the space of  pairs $((u_1, z_1), (u_2, z_2))$ made of  discs $u_i$ with boundary on $L$ satisfying the equation $\bar{\p}_J u_i= \nu(u_i)$ in non-trivial classes $\la_1, \la_2$, and with marked points $z_i \in \p D^2$  that satisfy $u_1(z_1) = u_2(z_2)$. Let $\Oo (\la_1, \la_2)$ be the orientation on $ \Mm(\la_1, \la_2)$ defined above through $\Phi$, and $\Oo(\la)$ be the orientation of $\Nu \Mm(\la)$, the space of $(\nu,J)$-holorphic discs in class $\la = \la_1 + \la_2$. Then, near any such Fredholm regular pair, the gluing at  $p = u_1(z_1) = u_2(z_2)$ maps the orientation  $\Oo (\la_1, \la_2)$ to the boundary-orientation of $\Oo(\la)$.
 \end{lemma} 

   This lemma is, by now, a standard result in the theory of $J$--holomorphic disks with Lagrangian boundary conditions: see for instance \cite{FOOO} for a proof.

  Note that the full moduli space $\Nu \Mm(L,J,(f,g);\la,\bf n)$ is made, by definition, of the 
  above cells  $\Nu \Mm_{T, \psi} (L,J,(f,g))$ glued together; the above argument and lemma show 
  that they glue coherently (i.e their orientations match). Moreover,  it is evident that the 
  boundary of $\Nu \Mm(L,J,(f,g);\la,\bf n)$ is given by 
formula~ (\ref{eq:pre_bdry}) understood with signs: indeed, the boundary only consists of breaking of 
gradient trajectories of the ordinary Morse function $f$ and it is an elementary exercise to check 
that, near any breaking point, the orientation of the right hand side in formula~ (\ref{eq:pre_bdry})
 coincides with the boundary-orientation of the left hand side. 

(3) We must now show how the above orientations on the moduli spaces $$\Nu \Mm(L,J,(f,g); \la, {\bf n})$$ for general $\bf n$, induces orientations on the moduli spaces $\Nu \Mm^x_{x_1, \ldots, x_k} (L,J,(f,g); \la)$.

As before, we assume
that $L$ is oriented and we pick orientations of the stable manifolds  so that
the unstable and stable manifolds (in that order)  of the same critical point have
intersection number equal to $1$ (with respect to the orientation
of $L$).

 Consider now the following
commutative diagram:

$$
\begin{array}{c}
\Mm^x_{x_1, \ldots, x_k} (\la)  \\
\downarrow \iota \\
 \Mm_{k+1} (\la)  \\
 \downarrow \ev    \\
  L^{k+1}  \\
\uparrow \iota_1  \\
W^u_{x} \times W^s_{x_1} \times \ldots W^s_{x_k}
\end{array}
$$
where $\Mm_{k+1}(\la)$ denotes the moduli spaces of
$J$-holomorphic discs with boundary on L in class $\la$ with $k+1$
ordered marked points $z_0$ (the root), $z_1, \ldots, z_k$, where
$\iota, \iota_1$ denote the various inclusions and $\ev$ the
evaluation map. Up to perturbations of $(f,g)$,  the map $ev$ is
transverse to the product of the unstable manifold of $x$ and the
stable manifolds of the $x_{i}$'s in the target space $L^{k+1}$.

\begin{definition} \label{def:orientations} When $n$ is odd, we give the moduli space  $\Mm^x_{x_1, \ldots,  x_k} (\la)$ the orientation induced by  $\ev$; more precisely, we set the orientation of $\Mm^x_{x_1,  \ldots, x_k} (\la) $ so that
$$
\Oo (\Mm^x_{x_1, \ldots, x_k} (\la) ) \; + \;  {\rm coorientation} (\Mm^x_{x_1, \ldots, x_k} (\la) ) \; = \;  \Oo( \Mm_{k+1} (\la))
$$
where the co-orientation is induced by the pull-back by $\ev$ of the co-orientation of $W^u_{x} \times W^s_{x_1} \times  \ldots W^s_{x_k} $ in $L^{k+1} $.

   When $n$ is even we change this orientation by multiplying it by
$$
(-1)^{\Pi_{1 \le i \le k} (i+1) |x_i |} \; .
$$
\end{definition}

With this definition, our orientations verify:

\begin{lemma}
$\Oo ( \Mm^x_{x_1, \ldots, x_{i+1}, x_i, \ldots, x_k} (\la)) = (-1)^{|x_i| |x_{i+1}|} \Oo (\Mm^x_{x_1, \ldots, x_i, x_{i+1}, \ldots, x_k} (\la)).$
\end{lemma}

\proof{} Let $$W_1=W^u_{x} \times W^s_{x_1} \times \ldots \times W^s_{x_i} \times W^s_{x_{i+1}}
\times \ldots W^s_{x_k}~,~$$
$$W_2=W^u_{x} \times W^s_{x_1} \times \ldots \times W^s_{x_{i+1}} \times W^s_{x_i}
\times \ldots W^s_{x_k}$$ and consider the following commutative diagram:
$$
\begin{array}{ccc}
\Mm^x_{x_1, \ldots, x_i, x_{i+1}, \ldots, x_k} (\la) & & \Mm^x_{x_1, \ldots, x_{i+1}, x_i, \ldots, x_k} (\la) \\
\downarrow \iota & & \downarrow \iota \\
 \Mm_{k+1} (\la) & \stackrel{\tau}{\longrightarrow} & \Mm_{k+1} (\la) \\
 \downarrow \ev & &  \downarrow \ev   \\
  L^{k+1} & \stackrel{\bar{\tau}}{\longrightarrow} & L^{k+1} \\
\uparrow \iota_1 & & \uparrow \iota_2 \\
W_1 & &
W_2
\end{array}
$$
 Here $\iota, \iota_1, \iota_2$ are the various inclusions,
 $\tau$ permutes $z_i$ and $z_{i+1}$ and $\bar{\tau}$ permutes
 the two corresponding factors in $L^{k+1}$. Endow  the moduli
 space $\Mm^x_{x_1, \ldots, x_i, x_{i+1}, \ldots, x_k} (\la)$ with the
 orientation induced by $\ev^{-1}(W_1)$  (ie so that the co-orientation of $W_1$ in $L^{k+1}$
 lifts to the co-orientation of $\Mm^x_{x_1, \ldots, x_i, x_{i+1},
 \ldots, x_k} (\la)$ in $\Mm_{k+1}(\la))$ and similarly  endow $\Mm^x_{x_1,
 \ldots, x_{i+1}, x_i, \ldots, x_k} (\la)$ with the orientation induced
 by $\ev^{-1}(W_2)$ .   Obviously, $\tau$ identifies $\Mm^x_{x_1, \ldots,
 x_i, x_{i+1}, \ldots, x_k} (\la)$ with $\Mm^x_{x_1, \ldots, x_{i+1}, x_i,
 \ldots, x_k} (\la)$. We claim that it preserves or reverses the
 orientation if $\dim L + \dim W^s_{x_i}   \dim W^s_{x_{i+1}} +1$
 is even or odd respectively.  Denote by $\Oo_1$ the co-orientation
 of $\Mm^x_{x_1, \ldots, x_{i}, x_{i+1}, \ldots, x_k} (\la)$ in $\Mm^x_{k+1}$
 and by $\Oo_2$ the co-orientation of $\Mm^x_{x_1, \ldots, x_{i+1}, x_i, \ldots, x_k} (\la)$ in $\Mm^x_{k+1}$. To prove the claim, and since the map $\tau$ clearly reverses the orientation of $\Mm^x_{k+1}$ (because it permutes two successive $S^1$-factors), it suffices to show that the co-orientation of $\Mm^x_{x_1, \ldots, x_{i+i}, x_{i}, \ldots, x_k} (\la)$ in $\Mm^x_{k+1}$, defined as the image $\tau(\Oo_1)$, agrees with $\Oo_2$ iff  $\dim L +  \dim W^s_{x_i})   \dim W^s_{x_{i+1}} $ is even. But, because the diagram is commutative, the co-orientation $\tau(\Oo_1)$ can as well be obtained throught the map $\bar{\tau}^{-1} \circ ev$. Thus $\tau(\Oo_1)$ agrees with $\Oo_2$ if and only if the image by $\bar{\tau}$ of the co-orientation of $W_1$ inside $L^{k+1}$ agrees with the co-orientation of $W_2$
inside $L^{k+1}$. But this is the case exactly when the parity of
$\dim L + \dim W^s_{x_i} \dim W^s_{x_{i+1}}$ is even since the first
term takes into account the change in orientation of the total space
$L^{k+1}$ while the second takes into account the discrepancy  between
the orientation of the $\bar{\tau}$-image of $W_1$ and the orientation of $W_2$.
Summarizing, $\tau$ identifies $\Mm^x_{x_1, \ldots, x_i, x_{i+1},
\ldots, x_k} (\la)$ with $\Mm^x_{x_1, \ldots, x_{i+1}, x_i, \ldots, x_k} (\la)$
 and the change in orientation is given by the parity of
$$
\begin{array}{lll}
 & n +  (n-\ind (x_i)) (n-\ind (x_{i+1})) +1 &    \\
 =  &    n + n^2 - n(\ind (x_i) + \ind (x_{i+1})) + \ind (x_i)  \ind (x_{i+1}) +1 &  \\
 = & n(\ind (x_i) + \ind (x_{i+1})) + \ind (x_i)  \ind (x_{i+1} ) +1 & \\
 = & n(|x_i| + |x_{i+1}|) + |x_i|  |x_{i+1}| + |x_i| +  |x_{i+1}| &.
\end{array}
$$
When $n$ is odd, this is equal to $|x_i|  |x_{i+1}|$  as claimed in the
statement of this lemma, and when $n$ is even, the formula of the lemma
remains valid by our particular choice of orientations.
\QED

With these conventions, we can now identify the relevant orientations
in formula (\ref{eq:dsquare1}).
First, formula (\ref{eq:pre_bdry}) holds without any additional
signs. It then follows that in (\ref{eq:dsquare1}) the sign of the term
$(^{\nu}\bar\Mm^x_{<S',y>}(\la'))\times
(^{\nu}\bar\Mm^y_{S''}(\la''))$ is: $\epsilon(\{y\}, <S',y>)\epsilon(S'',S)$.
Here, $<S',y>$ is, as before, the ordered set obtained by the rearrangement of the
elements of $S'\cup \{y\}$ by following the order in $\Crit(f)$
and $\epsilon(T',T)$ is the sign necessary to move the ordered subset $T'\subset T\subset\Crit(f)$
so that it precedes $T\backslash T'$ (in counting these signs we always use the
convention $ab=(-1)^{|a||b|}ba$).

Using these conventions, it is an easy exercise to verify
$d^{2}=0$ where $d$ is defined by formula (\ref{eq:diff}).

\subsection{Invariance of the cluster homology.} \label{sec:cluster_invariance}

  Recall that for our choices of almost complex
structure $J$, Morse-Smale pair $(f,g)$ and perturbations $\nu$,
the cluster homology of $L$ is defined by:
$$\Cl H_{\ast}(L, J, (f,g),\nu)=H_{\ast}((S\Q<s^{-1}\ \Crit(f)>\otimes\ \Lambda)^{\wedge},d^{\nu})~.~$$

\begin{theorem}\label{thm:cluster_invariance}
This homology is independent up to isomorphism of these choices
and thus we define the cluster homology of $L$ by:
$$\Cl H_{\ast}(L)= \Cl H_{\ast}(L, J, (f,g),\nu)$$
where $J, (f,g), \nu$ is any generic choice.
\end{theorem}

\proof{} The proof is based on the usual argument
in Floer theory in which two complexes associated to two choices
of data are compared by relating the moduli spaces defining the
differentials of the two complexes through cobordisms consisting
of moduli spaces associated to ``homotopies" between the two
selections of data. For this, we need the following notions.

\subsubsection{Clustered moduli spaces for Morse cobordisms}\label{subsubsec:cluster_cob}
We shall use the notion of Morse cobordism between two Morse-Smale
pairs $(f,g)$ and $(f',g')$. This is a homotopy $$F:L\times
[0,1]\to \R$$ together with a metric $G$ on $L\times [0,1]$ with
the following properties:
\begin{itemize}
\item[-]$(F,G)|_{L\times\{0\}}=(f,g)$,
$(F,G)|_{L\times\{1\}}=(f'+k,g')$ for some constant $k$.
\item[-]$\frac{\partial F}{\partial t}(x,0)=0=\frac{\partial
F}{\partial t}(x,1)$ for all $x\in L$, $\frac{\partial F}{\partial
t}(x,t)>0$ for all $x\in L$, $t\in (0,1)$. \item[-] $(F,G)$ is a
Morse-Smale pair so that $$\Crit_{i}(F)= \Crit_{i}(f)\times\{0\}\cup
\Crit_{i-1}(f')\times\{1\}~.~$$
\end{itemize}

As described in Cornea-Ranicki \cite{CoRa}, given that $L$ is compact, Morse
cobordisms exist and are easy to contruct between any two
$(f,g)$, $(f',g')$.

Consider two sets of data, $J, (f,g)$ and $J', (f',g')$, so that
are defined the moduli spaces $$\Mm^{x}_{x_{1},\ldots,
x_{k}}(\la)=\Mm^{x}_{x_{1},\ldots,x_{k}}(L,J,(f,g);\la)\ , \
x,x_{1},\ldots, x_{k}\in \Crit(f)$$ and
$$(\Mm')^{x'}_{x'_{1},\ldots,
x'_{k'}}(\la)=\Mm^{x'}_{x'_{1},\ldots, x'_{k'}}(L,J',(f',g');\la)\
, \ x',x'_{1},\ldots, x'_{k'}\in \Crit(f')~.~$$

We also fix a Morse cobordism $(F,G)$ between $(f,g)$ and
$(f',g')$ as well as a smooth one-parameter family of almost
complex structures $J_{t}$ with $J_{0}=J$ and $J_{1}=J'$. We
denote by $\bar{\ga}:(L\times [0,1])\times\R\to L\times [0,1]$ the
negative gradient flow induced by $(H,G)$.

Obviously, at the heart of the construction will be certain moduli
spaces that we shall denote by $\Nn^{x}_{y_{1},\ldots,y_{k}}(\la)$
(or by $\Nn^{x}_{y_{1},\ldots,y_{k}}(\la,F)$ in case the Morse
cobordism $F$ needs to be explicitly mentioned) where $x\in
\Crit{f'}$, $y_{i}\in \Crit(f)$.

\

The definition of $\Nn^{x}_{y_{1},\ldots,y_{k}}(\la)$ is perfectly
analogous to that of $\Mm^{x}_{x_{1},\ldots,x_{k}}(\la)$ from
\S~\ref{subsec:disks_flow} with the modification that the
moduli spaces $\Mm_{\mathcal{T},\psi}(L,J,f)$ are replaced with
moduli spaces $\Nn_{\mathcal{T},\psi}$ in whose definition
conditions i. ii. and iii. from \S \ref{subsec:disks_flow} are
replaced by i'. ii'. iii'. below (of course, iv. and v. are still required):

\begin{itemize}
\item[i'.]$u_{\alpha}:D^{2}\to M\times\{d(\alpha)\}\subset M\times
[0,1]$ is a $J_{d(\alpha)}$-holomorphic disk with boundary on
$L\times\{(d\alpha)\}$ and with Maslov class $\mu(\alpha)$ when
$\alpha\in V_{D}$. \item[ii'.] when $\alpha\in V_{S}$,
$u_{\alpha}:S^{2}\to M\times\{d(\alpha)\}$ is a
$J_{d(\alpha)}$-pseudo-holomorphic sphere with Chern class
$2\mu(\alpha)$ .
 \item[iii'.] If $\alpha \in V_{D}$, then
$z_{\alpha\beta}$ are points in $D^{2}$ with the property that if
$\beta\in V_{D}$, $\alpha E \beta$, then $z_{\alpha\beta}\in
S^{1}$, $z_{\beta\alpha}\in S^{1}$ and $\bar{\gamma}_{\psi(\alpha
E\beta)}u_{\alpha}(z_{\alpha\beta})=u_{\beta}(z_{\beta\alpha})~.~$
If $\alpha\in V_{S}$, then $z_{\alpha\beta}\in S^{2}$ and if
$\beta\in Vertex(T)$ so that $\alpha E\beta$ or $\beta E\alpha$,
then $u_{\alpha}(z_{\alpha\beta})=u_{\beta}(z_{\beta\alpha})~.~$
 \end{itemize}

We proceed again as in \S~\ref{sec:cluster_moduli} to construct
$\Nn^{x}_{y_{1},\ldots,y_{k}}(\la)$ out of
$\Nn_{\mathcal{T},\psi}$ except that the analogue of condition vi.
uses $\bar{\gamma}$ - the negative gradient flow of the Morse
cobordism $F$ - instead of the flow $\ga$.

\

The expected dimension of these moduli spaces is:
$$\dim (\Nn^{x}_{y_{1},\ldots, y_{k}}(\la))= \ind_{f'}(x) - (\sum_{i}(\ind_{f}(y_{i})-1)+\mu(\la)-1~.~$$

As in the case of the moduli spaces
$\Mm^{x}_{x_{1},\ldots,x_{k}}(\la)$, there are perturbed moduli
spaces $^{\nu}\Nn^{x}_{y_{1},\ldots,y_{k}}(\la)$ so that at least
when $\ind_{f'}(x) - (\sum_{i}(\ind_{f}(y_{i})-1)+\mu(\la)-1=1$,
$^{\nu}\Nn^{x}_{y_{1},\ldots y_{k}}(\la)$ is a rational
$1$-dimensional cycle. Moreover these perturbations may be chosen
in a way compatible with the perturbations associated to the data
$(J,(f,g))$ and  $(J', (f',g'))$. The resulting
perturbed moduli space admits a compactification whose boundary is
described by the formula below (to avoid carrying a notation that would be too heavy, we omit
the perturbations from the writing of the moduli spaces):

\begin{equation}\label{eq:cob_moduli}
\begin{array}{l}
\partial \bar{\Nn}^{x}_{S}(\la)=
\bigcup_{S'\cup S''=S,y, \la'+\la''=\la}\bar\Nn^x_{<S',y>}(\la')\times \bar\Mm^y_{S''}(\la'')\ \bigcup\\[1ex]
\bigcup_{\la'+\la''=\la,x_{1},\ldots,x_{s}}
[\bar{(\Mm')}^x_{x_{1},\ldots x_{s}}(\la')
  \times_{\cup S_{i}=S,\sum \la_{i}=\la''}\  (\bar\Nn^{x_{i}}_{S_{i}}(\la_{i}))]~.~
\end{array}
\end{equation}

Here $y, y_{i}\in \Crit(f)$ and $x,x_{i}\in \Crit(f')$.

Using the moduli spaces above and, as before, for $x\in \Crit(f')$, $y_{i}\in
\Crit(f)$ with $\ind_{f'}(x) -
\sum_{i}(\ind_{f}(y_{i})-1)+\mu(\la)-1=0$, we let
$$b^{x}_{y_{1},\ldots,y_{k}}(\la)=\#(^{\nu}\Nn^{x}_{y_{1},\ldots,y_{k}}(\la))$$
where the elements of $^{\nu}\Nn^{x}_{y_{1},\ldots,y_{k}}(\la)$
are counted with signs (the orientations are defined in a way similar to that
used in \S\ref{subsec:orientation}).

We now define:
 $$\phi^{F}: \Cl(L,J';(f',g'))\to \Cl(L,J;(f,g))$$
by the formula
\begin{equation}\label{eq:comparison_morph}\phi^{F}(x)=\sum b^{x}_{y_{1},\ldots,
y_{k}}(\la)y_{1}\ldots y_{k}e^{\la}
\end{equation} and from
(\ref{eq:cob_moduli}), we immediately obtain (again, up to a sign verification) that $\phi^{F}$ is a
morphism of chain complexes.

\begin{remark}{\rm In order to get the relation $d
\phi^F = \phi^F d$ that makes $\phi^F$ into a complex morphism, it
is required that the top dimensional stratum of the boundary of
the moduli space $\bar{\Nn}^{x}_{S}(\la)$ be asymmetric: indeed,
applying first
 $d$ to some critical point $x$ and then $\phi^F$ means
 that there might be many breaking points on the $M \times \{1\}$ side,
 whereas applying first $\phi^F$ to some critical point $x$ and
 then $d$ makes appear only one breaking point of the same
 moduli space $\bar{\Nn}^{x}_{S}(\la)$ (this asymmetry is of
 course due to the fact that there is only one root but many ends).
 The fact that, indeed, the top dimensional stratum
 of $\bar{\Nn}^{x}_{S}(\la)$  behaves in this way can be seen as follows.
 Notice first that one may regard an
element of $\Nn^{x}_{y_{1},\ldots, y_{k}}$ as a ``tree'' in
$M\times [0,1]$ having disks (or spheres) instead of vertices and
with its edges represented by flow lines of $\bar{\ga}$. The
origin of the tree is at $y$ and it has terminal vertices in
$x_{1},\ldots x_{k}$. The boundary of the compactification of this
moduli space consists of broken trees, the breaking points
occuring in critical points of $F$. Therefore, these breaking
points appear on $L\times \{1\}$ or on $L\times\{0\}$ and if a
broken tree appears in the top dimensional stratum of the boundary
of $\Nn^{x}_{y_{1},\ldots,y_{k}}$, then it has breaking points
only on one of these. The interesting point is that such a top
dimensional stratum tree may be broken in {\em multiple} points on
$L\times\{1\}$ (but brakes at only {\em one} critical point on $L \times \{0\}$). In essence, this happens because if a flow line
inside the tree breaks on $L\times \{1\}$ then the disk (or
sphere) at the origin of this flow line is contained in
$M\times\{1\}$ and the only way to leave $M\times\{1\}$ via a flow
line of $\bar{\ga}$ is through a critical point of $f'$. A simple
dimension computation shows that indeed this multiple breaking points
phenomenon appears in the top dimensional stratum of the boundary.
}\end{remark}

\

To show that $\phi^{F}$ induces an isomorphism in homology, one
could apply a variant of the standard approach in Morse-Floer
theory which is to show that a cobordism of Morse cobordisms
induces a chain homotopy. However, things are more delicate here
because, while this approach provides the definition of the
chain-homotopy on the generators of $\Cl(-)$, it is not clear how
to directly extend this chain homotopy to the words of length
longer than one in the commutative, graded algebra
$S\Q<Crit(f)>\otimes \Lambda$ (as an interesting side comment,
this is easy to do for non-commutative free algebras \cite{An}).
In view of this, our proof of the existence of an isomorphism
between $\Cl H(L,J;(f,g))$ and $\Cl H(L,J';(f',g'))$ will be based
(besides the construction of $\phi^{F})$ on the properties of a
differential filtration which is of independent interest.

\subsubsection{The word-area filtration.} \label{subsec:area_filtration} Denote by
$\epsilon_{D}$ the infimum of $\int_{D^{2}}u^{\ast}\omega$ over the set of 
maps $ u : (D^{2},S^{1})\to (M,L)$ which are $J_t$-holomorphic for some $t \in [0,1]$ and non-constant. This number is
strictly positive and we use it to define a filtration of the
cluster complex $\Cl(L,J;(f,g))$:
$$F^{k}\Cl(f)=\Q<\ x_{1}x_{2}\ldots x_{m}e^{\la}\
|\ \ m\frac{\epsilon_{D}}{2}+\omega(\la)\geq k
 \frac{\epsilon_{D}}{2}\ \ >~.~$$

Notice that the cluster differential preserves this filtration.
Indeed, if a term in the differential of $x_{1}\ldots
x_{m}e^{\la}$ is of the form $y_{1}\ldots y_{l}e^{\la'}$ and
$l<m$, then $l=m-1$ and $e^{\la'}=e^{\la}e^{\la''}$ with $\la''$
represented by a sum of pseudoholomorphic disks. Therefore,
$\omega(\la'')\geq \epsilon_{D}$ and our filtration is
differential.  This remark implies in fact a little more: if we
set the {\em weight} of a monomial $x_{1}\ldots x_{k}e^{\la}$ to
be
$$w(x_{1}\ldots x_{k}e^{\la})=k+2\frac{\omega(\la)}{\epsilon_{D}}$$
  and if $m\in \Cl(L,J;(f,g))$ is a monomial, then we may write
$dm=d_{0}m+\sum_{i}m_{i}$ with $d_{0}$ the Morse differential and
the $m_{i}$ are monomials with $w(m_{i})\geq w(m)+1$.

\

We now consider the spectral sequence $E^{r}(f)$ associated to the
filtration $F^{k}\Cl(f)$ which we shall call further the {\em word-area
spectral sequence}. The remark above implies that the total vector
space of the term $E^{1}(f)$ is isomorphic to
$(S(s^{-1}H_{\ast}(L;\Q))\otimes\La)^{\wedge}$ because the $0$-order differential in
the spectral sequence coincides with the Morse differential.

We obviously have a similar filtration $F^{k}\Cl(f')$ and an
associated spectral sequence $E^{r}(f')$. It is easy to notice (by
an argument similar to the one applied above to the differential)
that the chain morphism $\phi^{F}$ preserves the word-area filtration
and thus induces a morphism of spectral sequences
$E^{r}(\phi^{F})$. Moreover, $E^{0}(\phi^{F})$ coincides with the
morphism induced by the Morse comparison morphism associated to
the Morse cobordism $F$ and, as this morphism induces an
isomorphism in homology, we deduce that $E^{1}(\phi^{F})$ is an
isomorphism. This means that $E^{r}(\phi^{F})$ is an isomorphism
for all $r\geq 1$.

The cluster differential of a term $x\in \Crit(f)$ may contain
infinitely many terms (even if only finitely many terms of weight
bounded by some $k$ can appear). Consequently, knowing that
$E^{\infty}(\phi^{F})$ is an isomorphism does not imply directly
that $H_{\ast}(\phi^{F})$ is an isomorphism. We will now prove
this last fact under the additional assumption that $\phi_{0}$,
the Morse-comparison map relating the Morse complexes of $f'$ and
of $f$ and induced by $F$, is surjective.

\

We start by noting that if $\phi_{0}$ is surjective, then
$\phi^{F}$ is surjective too. Indeed, if we denote by
$x_{1},\ldots , x_{n}$ the critical points of $f$, we may
 write $x_{i}=\phi_{0}(z_{i})$ with $z_{i}\in \Q<\Crit(f)>$. We also have
$\phi_{0}(z_{i})=\phi^{F}(z_{i})-\sum_{\omega(\la)>0}m_{\la}e^{\la}$
where $m_{\la}$ are monomials in the generators $x_{k}$'s.  In all
the monomials $m_{\la}$ for which $\omega(\la)$ is minimal, we
replace the $x_{k}$'s by the expression $\phi^{F}(z_{k})-\sum
m_{\la}e^{\la}$. This process replaces these monomials with
monomials written in the $\phi^{F}(z_{i})$'s summed with monomials
(in the $x_{i}$'s) which are multiplied with some $e^{\la'}$ with
$\omega(\la')>\omega(\la)$. Recursively, this shows that $x_{i}$
can be written (as formal series) only in terms of the
$\phi^{F}(z_{i})$'s and thus $\phi^{F}$ is surjective. This
argument does also imply that the restriction of $\phi^{F}$ to
$F^{k}\Cl(f')$ is onto $F^{k}\Cl(f)$.

Now let $K=\ker(\phi^{F})$. We intend to show that $K$ is acyclic.
Let $$A^{n}=\Cl(f')/(F^{n}\Cl (f'))\ \ {\rm  and\ let}\ \
B^{n}=\Cl(f)/(F^{n}\Cl(f))~.~$$ Clearly, our filtrations induce
filtrations on both $A^{n}$ and $B^{n}$ and the map $\phi^{h}$
induces chain morphisms $\phi^{n}:A^{n}\to B^{n}$ which respect
the filtrations. As before, we see that at the level of the induced
spectral sequences, $E^{1}(\phi^{n})$ is an isomorphism. Thus
$E^{\infty}(\phi^{n})$ is an isomorphism. The differential of the
algebra generators in $A^{n}$ and $B^{n}$ are {\em finite} sums of
monomials and so we get that $H_{\ast}(\phi^{n})$ is an
isomorphism for all $n\geq 1$. Thus the kernel $K_{n}=Ker
(\phi^{n})$ is acyclic. This kernel fits into an exact sequence
$$0\to K\cap F_{n}\Cl(f')\longrightarrow K\stackrel{p_{n}}{\longrightarrow} K_{n}\to 0$$
and the projections $p_{n}$ commute with the obvious projections
$K_{n+1} \stackrel{q_{n}}{\longrightarrow} K_{n}$ which are thus
surjective. We have $$K=\lim_{\leftarrow}K_{n}~.~$$

Let $\sigma\in K$ so that $d\sigma =0$. As $K_{n}$ is acyclic,
there are elements $\xi_{n}\in K_{n}$ so that
$d\xi_{n}=p_{n}\sigma$, $\forall n\geq 1$. Therefore,
$d(q_{n-1}\xi_{n}-\xi_{n-1})=0$. Using the acyclicity of $K_{n-1}$
again there is $\tau\in K_{n-1}$ so that $d\tau=
q_{n-1}\xi_{n}-\xi_{n-1}$. But $q_{n-1}$ is surjective, so there is
$\tau'\in K_{n}$ with $q_{n-1}\tau'=\tau$. Now define a
preturbation of $\xi_{n}$ by $\xi'_{n}=\xi_{n}-d\tau'$. Notice
that $d\xi'_{n}=p_{n}\sigma$ and $q_{n-1}\xi'_{n}=\xi_{n-1}$. By
applying this process recursively, we get a sequence of elements
$\xi'_{n}$ so that $q_{n-1}\xi'_{n}=\xi'_{n-1}$,
$d\xi'_{n}=p_{n}\sigma$, $\forall n$. This implies that we may
define $\xi=\lim_{\leftarrow}\xi'_{n}$ and that we have
$d\xi=\sigma$, which shows that $K$ is acyclic and that $\phi^{F}$
induces an isomorphism in homology.

\

To prove that $\Cl H(L,J;(f,g))$ and $\Cl H(L',J'; (f',g'))$ are
isomorphic without any surjectivity condition for $\phi_{0}$, it
is sufficient to notice that for any two Morse- Smale pairs
$(f,g)$ and $(f',g')$, there is a chain of Morse cobordisms:
$$(f,g)\leftarrow \odot\rightarrow \odot \leftarrow \odot \ldots \odot\rightarrow (f',g')$$
so that each Morse cobordism in the chain induces a surjective
chain map. This is quite easy to show either as a consequence of
bifurcation analysis or by the rigidity results in \cite{CoRa}.

This completes the proof of Theorem~\ref{thm:cluster_invariance}.

\begin{remark}\label{rem:canonical_comp}{\rm
We expect the comparison isomorphism $H_{\ast}(\phi^{F})$ to
be canonical, as it is in Morse theory or in
Floer theory. Unfortunately, the argument above does not provide
this stronger result. However, it does imply that the morphisms
induced by $\phi^{F}$ on the pages of order $1$ and more of the
area spectral sequence are canonical. Indeed, $E^{1}(\phi^{F})$ is
induced by the morphism in Morse homology which is canonical. So
$E^{1}(\phi^{F})$ does not depends of the choice of $h$. But this
means that $E^{r}(\phi^{F})$ does not depends of this choice for
$r\geq 1$. }\end{remark}

\subsection{Various special cases.}

We start by recalling from Example \ref{rem:no_bubbbling} that in the absence of bubling
we have: $$\Cl H_{\ast}(L,J;(f,g)\simeq S(s^{-1}H_{\ast}(L;\Q))\otimes\La)^{\wedge}~.~$$

\subsubsection{Higher dimensional cluster moduli spaces}
It is an interesting question to see, in the spirit of \cite{BaCo}, how to encode the information, in moduli spaces of clusters of any dimension (not only $0$-dimensional), that has to do with the loop structure of the ``boundary'' of each cluster tree. By this, we mean that we could try to assign, to each cluster  configuration in $\Mm^{x}_{x_{1},\ldots,x_{k}}(\la)$, a loop in $L$ that starts at $x$ and goes round the tree (i.e round each $J$-holomorphic disk once, and twice (back and forth) on each gradient flowline) and comes back to $x$. Because special marked points may cross each other (not directly, but through ghost bubbles), this loop structure is not  well-defined, but there is certainly a tree-structure. 
For now just notice that, in fact, these
higher dimensional moduli spaces appear implicitly in the cluster differential because
we may view $\Mm^{x}_{x_{1},\ldots,x_{k}}$ as the intersection of
$\Mm^{x}_{x_{1},\ldots, x_{k-1}}$ with the unstable manifold of $x_{k}$ and thus
high dimensional moduli spaces become visible in the cluster differential once
they are reduced by intersection to $0$-dimensional spaces. This is, of course, similar
to the process leading to Gromov-Witten invariants if the moduli spaces of
disks are replaced by moduli spaces of pseudoholomorphic spheres.

\subsubsection{Relation to Floer
homology.}\label{subsec:relation_floer}
In this subsection we assume
that the minimal Maslov index of a pseudo-holomorphic disk  with
boundary on $L$ is at least equal to $2$. We shall denote this
minimal Maslov index of some pseudo-holomorphic disk by $\mu_{min}$. For $x,y\in \Crit(f)$ and
$\la\in \La$, let $\Ff^{x}_{y}(\la)$ be the moduli space
constructed in the same way as the moduli spaces
$\Mm^{x}_{x_{1},\ldots, x_{k}}(\la)$ except that all graphs
$\mathcal{T}$ that are used in the construction are {\em linear}
chains (they consist simply of chains of edges $xEy_{1},
y_{1}Ey_{2},\ldots, y_{s}Ey$). We shall call further these cluster
configurations {\em linear cluster elements}. The key point is
that, because $\mu_{min}\geq 2$, the compactification of
$\Ff^{x}_{y}(\la)$, $\bar{\Ff}^{x}_{y}(\la)$, has the property
that (after appropriate perturbations) if
$|x|-|y|+\mu(\la)-1=1$, then
\begin{equation}\label{eq:linear_tree}\partial\bar{\Ff}^{x}_{y}(\la)=
\sum_{z,\la'+\la''=\la}\bar{\Ff}^{x}_{z}(\la')\times\bar{\Ff}^{z}_{y}(\la'')~.~
\end{equation}
This happens because bubbling in $\Ff^{x}_{y}(\la)$ can be
transformed into an interior point just by using linear trees.
Indeed, fix one disk $D$ inside a linear chain $\kappa$ of flow lines
and disks joining $x$ to $y$. As $\kappa$ is a linear chain, there are
precisely two double points on $D$. Assume now that, by bubbling
off inside $\bar{\Ff}^{x}_{y}(\la)$, the disk $D$ splits into two
other disks $D'$ or $D''$. In this case, for dimensional reasons,
each of $D'$ and $D''$ will carry one of these double points. More precisely,
suppose that a linear tree from $x \in \Crit(f)$ to $y \in \Crit(f)$ is the starting point
of a $1$-dimensional cobordism of configuartions, and suppose that at some time $t_0 > 0$,
it reaches a non-linear tree by bubbling off from some disk $D_0$ at a point which is disctinct from the two special points of $D_0$. The class $D_0$ would then split in $D_0' + D_0''$ where say $D_0'$ would still be part of a new linear tree $\kappa'$ (but not $D_0''$). But this is impossible since the new linear tree $\kappa'$ would then have dimension at least $2$ less that the original tree $\kappa$ (by our assumption on the minimal Maslov index). Since $\kappa$ was assumed to be of dimension $0$, this cannot happen generically in a one-parameter family.

 This
immediately implies (\ref{eq:linear_tree}).

Clearly, formula (\ref{eq:linear_tree}) may be used to define a
chain complex $$\Cc(L,J;(f,g))=(\Q<\Crit(f)>\otimes (\La)^{\wedge}, d)$$ with
$$d(x)=\sum_{y,\la}(\#\Ff^{x}_{y})ye^{\la}$$
where $x\in \Crit(f)$, the sum is over all $y\in \Crit(f)$, $\la\in
\La$ so that $\dim(\Ff^{x}_{y}(\la))=0$ and the elements in
$\Ff^{x}_{y}(\la)$ are counted with signs. The completion $^{\wedge}$
only applies now to the group ring $\La$ which reduces here to the usual Novikov
ring.

Essentially, this same complex has appeared before in an annoucement by
Oh (see the very end of the paper \cite{Oh}) who has indicated that its homology should coincide
with the standard Floer homology. Indeed, it is easy to show first
that the homology of this complex is independent of $J, (f,g)$
(for example, by the ``chain homotopy" method typical in Floer
theory) and then it is not difficult to find a chain map defined
on the Floer complex $$\psi:FC(L,J,H)\to \Cc(L,J;(f,g))~.~$$ Here
we use the Floer complex associated to the action functional
$$\Aa:\widetilde{\Pp_{0}(L,L)}\to \R$$ where $\Pp_{0}(L,L)$ is the
space of paths $$\gamma: [0,1]\to M, \ \gamma(0),\gamma(1)\in L$$
which represent the trivial element in $\pi_{1}(M,L)$,
$$\Aa(\gamma)=-\int u^{\ast}\omega +\int_{0}^{1}H(\gamma(t)dt$$
and $\ \widetilde{\ \ }\ $ represents the universal cover. The
Floer trajectories associated to this setting are semi-cylinders
with their boundary on $L$ which join orbits of the Hamiltonian
flow which start from $L$ and arrive in $L$ in time $1$. The
comparison map $\psi$ may be defined by means of moduli spaces of
objects that ressemble semi-disks connected by a flow line to a
chain of flow lines and disks (this is the Piunikin, Salamon, Schwarz method \cite{PSS}).
The same argument with one additional parameter shows that any two
comparison maps obtained in this way are chain homotopic. Comparison
maps $\psi':\Cc(L,J;(f,g))\to FC(L,J,H)$
 can also be constructed and they verify similar properties. It is then easy to
see that this implies that $\psi$ induces an isomorphism in
homology.

\

Finally, we discuss the relation between the complex
$\Cc(L,J;(f,g))$ and the cluster complex.

\begin{proposition}
The natural projection $$l:S(s^{-1}\Q<Crit(f)>)\to s^{-1}\Q<Crit(f)>~,~$$ which
sends to $0$ all the words of length longer or equal than $2$ as
well as the unit, induces a morphism of chain complexes: $$l:
\Cl(L,J;(f,g))\to s^{-1}\Cc(L,J;(f,g))~.~$$
\end{proposition} \begin{proof}
The only reason that could prevent the map $l$ to be a morphism of chain
complexes is that a priori there might exist some element in
a $0$-dimesional moduli space $\Mm^{x}_{y}(\la)$ whose
corresponding tree $\mathcal{T}$ is not a linear tree but rather a
tree with some branches. As this tree has one root and one single
end, we may consider the linear subtree $\mathcal{T}'$ which
connects the root to the end. Let $\la'$ be the total Maslov class
of the disks corresponding to the vertices which are not on
$\mathcal{T}'$. Then the cluster configuration corresponding to
$\mathcal{T}'$ provides an element of $\Mm^{x}_{y}(\la-\la')$. But
the dimension of this last moduli space is
$|x|-|y|+\mu(\la)-\mu(\la')-1$. As $|x|-|y|+\mu(\la)-1=0$ and
$\mu(\la')\geq 2$, this is impossible (for a generic $J$).
\end{proof}

\begin{remark}\label{rem:sseq_comparison}{\rm
a. There is obviously an area filtration also on the complex
$\mathcal{C}(L,J;(f,g))$. The spectral sequence associated to this
filtration is essentially a variant of Oh's spectral sequence from
\cite{Oh2}. It is useful to notice also that the map $l$ defined
above preserves this filtration and we have therefore an induced
morphism of spectral sequences which compares our area spectral
sequence to Oh's spectral sequence.

b. One may view the definition of cluster homology as the result of
an attempt to define a quantum product for Lagrangian submanifolds to be defined on
Floer homology. Without the condition
$\mu_{min}\geq 2$, such a product cannot be defined. However,
if this condition is satisfied,  it is easy to define a bilinear map
which descends to Floer
homology and which is just the dual of the linear map
$$s\Cc(L,J;(f,g))\to s\Cc(L,J;(f,g) \otimes s\Cc(L,J;(f,g))$$
given by the quadratic part of the cluster differential (for a
generator $x\in \Crit(f)$, this is the part of $dx$ which is written
as a sum of words of length two). Indeed, the fact that the linear
part of the cluster differential is itself a differential implies
that this bilinear map is a chain map. Here $s$ indicates raising
the degrees of the generators in $\Cc$ by one unit. This bilinear
map is associative but is not a product due to the lack of a unit
(as Paul Biran pointed out to us). Indeed, our cluster
construction does
not ``integrate" the usual cup-product (this is done on purpose).

c. In the case when $\mu_{min}\geq 2$, it is easy to use directly
the complex $\mathcal{C}(L,f)$ described above to
show that the conclusion (i) of Corollary \ref{cor:uni_ruled} holds: the argument is exactly the one
in Corollary \ref{cor:uni_ruled}. Moreover, in this case the result remains true
without any relatively spin or orientability restrictions because it can be seen that
 the complex $\mathcal{C}(L,f)$ is  well defined over $\Z/2$. }\end{remark}

\section{The fine Floer complex and its homology}

\subsection{Moduli spaces for the fine Floer Homology.} \label{sec:fine_moduli}

The purpose of this section is to give a precise definition of the moduli spaces
used in the definition of the differential of the fine Floer complex.
They were very briefly described in \S~\ref{subsubsec:fine_def_first}. As in that
subsection, we assume that we have a transversal pair of orientable,
relatively spin Lagrangian submanifolds $L_0$ and $L_1$,
as well as an almost complex structure $J$, and Morse-Smale pairs  $(f_i,g_i)$ on $L_i$.
We denote by $\gamma^{i}$ the negative gradient flow of $f_{i}$. We also assume that
$\Crit(f_{j})\cap L_{0}\cap L_{1}=\emptyset$ for $j=0,1$ and that the points in 
$L_{0}\cap L_{1}$ belong to the unstable manifolds of the minima of $f_{j}$.

Th readers should have in mind the large Novikov ring $\bar\La$ as well as the desired coefficient
ring for the fine complex $\mathcal{R}$ from equation (\ref{eq:coefficients_fine}).

Recall also that the moduli spaces to be described are associated to
$a,b\in I_{\eta}$, $\la\in \bar\La$, $r_{i}\in \Crit(f_{0})\cup\Crit(f_{1})$ and
are denoted by
$\Ww^{a}_{r_{1},\ldots,r_{k};b}(\la)$.
They are defined in a way similar to the moduli spaces $\Mm^{a}_{x_{1},\ldots,x_{i}}(\la)$
from \S\ref{sec:cluster_moduli}, but there are a number of modifications
that we now describe.

As in \S\ref{sec:cluster_moduli}, we start with trees and consider 
couples $(\mathbf{u},\mathbf{z})$ with the difference that, in this case, there are at least two marked points
$z_{0},z_{-1}$ on the root of the tree (in other words $0,-1\in \phi^{-1}_{T}(v_{0})$). Besides this,
 the properties (i.-vi.) become:

\begin{itemize}
\item[i''.] The map $u_{v_{0}}$ verifies:
\begin{equation}\label{eq:Flo_eq}
u_{v_{0}}:\R\times [0,1] \to M \ , \ \frac{\partial u}{\partial s}+J\frac{\partial u}{\partial t}=0~,~
\end{equation}
$\lim_{s\to -\infty}u(s,t)=a$, $\lim_{s\to +\infty}u(s,t)=b$; $u_{v_{0}}(\R\times\{i\})\subset L_{i}$.
We view $u_{v_{0}}$ as a continuous map defined
on $D^{2}$ and which sends the marked point $z_{0}$ to $a$ and
 sends the second marked point, $z_{-1}$, to $b$.
\item[ii''.] For $\alpha\not=v_{0}$, $\alpha\in V_{D}$, $u_{\alpha}$ are disks
so that $u_{\alpha}(\partial D^{2})\subset L_{i}$ (for some $i\in\{0,1\}$) and
$u_{\alpha}$ is a $J$-pseudoholomorphic disk.
In both cases, these disks are of
Maslov class $\mu(\alpha)$.
\item[iii''.] For $\alpha\in V_{S}$, $u_{\alpha}$ is a $J$-holomorphic
sphere of Chern class $2\mu(\alpha)$.
\item[iv''.] Properties iii., iv. from \S\ref{sec:cluster_moduli} are verified with the
only modification that we require $\gamma^{j}_{\psi(\alpha
E\beta)}(u_{\alpha}(z_{\alpha\beta}))=u_{\beta}(z_{\beta\alpha})$ in case
$u_{\alpha}(z_{\alpha\beta})\in L_{j}$.
\item[v''.] If $u_{\alpha}$ is constant, then $n_{\alpha}\geq 3$. Moreover, if $\alpha\not=v_{0}$,
there is a chain of edges of zero length that connects $\alpha$ to a vertex $\beta$
with $u_{\beta}$ non-constant.
\item[vi''.] Instead of vi, we have the following property:  there exists an ordered subset
 $(z_{1},\ldots,z_{k})$
 of marked points  (again called terminations) where
 $i\in \phi^{-1}_{T}(V_{D})$ so that $u_{\phi_{T}(i)}(z_{i})\in L_{j}$ for $j=0$ or $j=1$,
and  $$\lim_{t\to
+\infty}\gamma^{j}_{t}(u_{\phi_{T}{i}}(z_{i}))=r_{i}~.~$$
\end{itemize}

An element $(\mathbf{u},\mathbf{z})$ as defined above may be viewed roughly
as a pseudoholomorphic strip (of ends $a$ and $b$)
corresponding to the root of the tree $v_{0}$ together with a number
of cluster (sub)trees attached at some incidence points $z_{v_{0}\eta_{i}}$ (and maybe some
additional pseudoholomorphic spheres). Clearly, according to whether
$u_{\eta_{i}}(S^{1})\subset L_{0}$ or $u_{\eta_{i}}(S^{1})\subset L_{1}$,
the corresponding cluster subtree is included in $L_{0}$ or in $L_{1}$.
It is also clear in the point vi'' of the definition above that, as $r_{i}$ belongs to a single
$L_{j}$, all of the cluster subtrees starting in $v_{0}$ and having $r_{i}$ as termination
must lie on $L_{j}$. As in the construction of the moduli spaces $\Mm$ in $\S\ref{sec:cluster_moduli}$,
the number of marked points for a moduli space $\Ww^{a}_{r_{1},\ldots,r_{k};b}(\la)$ is the minimal
possible and is equal to $k+2$ (because, in this case, the root carries two special marked
points). This means, in particular, that a constant strip by itself does not belong to
our moduli spaces.

Clearly, a permutation $\sigma$ of the indexes of the marked points $z_{i}$ induces
an identification (up to orientation matters which will be discussed below):
$$\Ww^{a}_{r_{1},\ldots, r_{k};b}(\la)\to \Ww^{a}_{r_{\sigma(1)},\ldots, r_{\sigma(k)};b}(\la)~.~$$
Assume now that we fix an order on $\Crit(f_{j})$ for $j\in 0,1$ and we then
order $\Crit_(f_{0})\cup\Crit(f_{1})$ so that the critical points of $f_{1}$ follow
those of $f_{0}$.  For an ordered subset $S=(x_{1},\ldots,x_{k},y_{1},\ldots,y_{s})$, $x_{i}\in\Crit(f_{0})$,
$y_{i}\in \Crit(f_{1})$, we set
$$\Ww^{a}_{S,b}(\la)=\Ww^{a}_{x_{1},\ldots,x_{k},y_{1},\ldots y_{s};b}(\la)~.~$$

\
\subsubsection{Compactification}
The moduli spaces above admit compactifications
$$\bar{\Ww}^{a}_{S,b}(\la)$$
and after appropriate perturbations these are singular cycles
(with boundary) of dimension $$|b|-|a|-\sum_{i}|x_{i}|+\mu(\la)-1$$
and, at least if they are of dimension $1$, their boundaries
satisfy a formula analogous to  (\ref{eq:dsquare1}).

To write this formula, we recall that
$S$ denotes the ordered set $(x_{1},\ldots,x_{k},y_{1},\ldots,y_{s})$. We set
$S_{0}=(x_{1},\ldots,x_{k})$, $S_{1}=(y_{1},\ldots,y_{s})$. We
need also to consider the clustered moduli spaces $\Mm$ for each
of $(L_{i},(f_{i},g_{i}),J)$. To simplify notation, is it useful to
define $\Mm^{r}_{R}(\la)$ for any $r\in \Crit(f_{i})$ and for $R$
any ordered subset of $ \Crit(f_{0})\cup \Crit(f_{1})$ as follows:
if $R\not\subset \Crit(f_{i})$, then $\Mm^{r}_{R}(\la)=\emptyset$
and, if $R\subset \Crit(f_{i})$, then
$\Mm^{r}_{R}(\la)=\Mm^{r}_{R}(L_{i},(f_{i},g_{i}), J;\la)$. With
this notation and when $a\not= b$, it is easy to see that the
formula is:

\begin{equation}\label{eq:dddsquare}
\begin{array}{l}\partial\bar{\Ww}^{a}_{S,b}(\la)=\bigcup_{y,\la=\la'+\la'', S=<S'\cup S''>}
\bar{\Ww}^{a}_{<S',y>, b}\times
\bar{\Mm}^{y}_{S''}(\la'')\ \cup \\[1ex] \hspace{0.75in}\bigcup_{c, \la=\la'+\la'', S=<S'\cup S''>}
\bar{\Ww}^{a}_{S',c}(\la')\times \bar{\Ww}^{c}_{S'',b}
\end{array}
\end{equation}

Notice that, in a one-dimensional moduli space as above, one and
only one cluster tree attached to a strip might slide to one of
the ends of this strip thus creating a boundary point. This type
of boundary point is taken into account above by the terms
$\bar{\Ww}^{a}_{S',a}(\la')\times \bar{\Ww}^{a}_{S'',b}$ and
$\bar{\Ww}^{a}_{S',b}(\la')\times \bar{\Ww}^{b}_{S'',b}$. Thus,
to take into account these terms, the only moduli spaces that need
to be considered are those which satisfy the auxiliary condition
that:

\begin{itemize}
\item[vii''.] if $u_{v_{0}}$ is constant then $n_{v_{0}}=3$.
\end{itemize}

As we only need to deal with one parametric moduli spaces here, we
will assume from now on that this condition is verified. This
condition means that a constant strip can only contain the
attaching point of a single cluster tree. Thus the moduli spaces
in which $u_{v_{0}}$ is the constant strip equal to $a\in L_{0}\cap L_{1}$
are of two types: $T^{a}_{S_{0}}(\la)\subset
\Ww^{a}_{S_{0},a}(\la)$ with $S_{0}\subset \Crit (f_{0})$ and
$R^{a}_{S_{1}}\subset \Ww^{a}_{S_{1},a}(\la)$ where $S_{1}\subset
\Crit (f_{1})$ (an element in these spaces corresponds to one
cluster tree attached to $a$ contained in $L_{0}$ in the case of
$T^{a}_{\ldots}$ and in $L_{1}$ in the case of $R^{a}_{\ldots}$).

\

We now show that the boundary formula
(\ref{eq:dddsquare}) remains true even in the case when $a=b$.
Indeed, there are only two possibilities which need to be
explicited in this case as they are different from the case of the
clustered moduli spaces. The first one is the case when $z_{0}$ and
$z_{-1}$ become identified. This can be viewed as the bubbling of
a disk with boundary on just one $L_{i}$, say $L_{0}$, which is
attached to the constant strip $a$ (of course some cluster trees
might be attached to this disk). This means that this bubbling can
only happen in a space $\Ww^{a}_{S_{0},a}(\la)$ with
$S_{0}\subset\Crit (f_{0})$. Moreover, this element is not a
boundary in $\Ww^{a}_{S_{0},b}(\la)$ because this configuration
also appears in the compactification of $T^{a}_{S_{0}}(\la)$ (when
the flowline joining $a$ to the root of some cluster tree tends to
$0$). The second case already appears in the discussion above:
inside $T^{a}_{S_{0}}(\la)$ (or similarly, inside
$R^{a}_{S_{1}}(\la)$) the length of the flow line joining the root
of some cluster tree to $a$ tends to $0$. Thus the previous
discussion together with the gluing result described in
\S\ref{subsubsec:diff_fine1} takes care of both of these cases and
hence provides the formula.

\

As in \S \ref{subsec:orientation}, the moduli spaces constructed
above also admit orientations. The only significant remark here is
that for the clusters on $L_{1}$, we must take the inverse orientations
compared of those described in  \S \ref{subsec:orientation}. We
take the orientation described there for the clusters on $L_{0}$.
With these choices, as in \S \ref{subsubsec:diff_fine1}, it is
easy to check that $d_{F}^{2}=0$.

 \

\

\subsection{Invariance of the fine Floer Homology.} \label{sec:fine_invariance}

  The goal of this section is to prove:

\begin{theorem} \label{thm:fine_invariance}
The homology $\FF H_{\ast}(L_{0},L_{1},\eta;J,(f_i,g_i))$ is independent up
to isomorphism of $J$, $(f_i,g_i)$ and if $\psi:M \to M$ is a
Hamiltonian diffeomorphism, then $\FF
H_{\ast}(L_{0},\psi (L_{1}),\eta')$ is isomorphic to $\FF
H_{\ast}(L_{0},L_{1},\eta)$ (here $\eta'=\psi (\eta)$).
\end{theorem}

\begin{proof}  First note that, since loops of Hamiltonian diffeomorphisms 
act trivially on the homology of $M$ (see for instance Lalonde, McDuff, Poletorvich \cite{LMP}), 
two different Hamiltonian isotopies from the identity to $\psi$ lead to the 
same class $\eta'$ given by concatenation of $\eta$ and the Hamiltonian 
isotopy applied to the endpoint of $\eta$. Thus the class $\eta'$ in the statement 
of the Theorem is well-defined. So
fix $L_{i}, f_{i},g_{i}$ as before and consider also any Hamiltonian
isotopy $\psi^{t}$, $t\in [0,1]$ from the identity to $\psi$. We denote $L'_{1}=\psi^{1}(L_{1}), L'_{0}=L_{0}$
and we also fix $(f'_{i},g'_{i})$ Morse-Smale pairs on $L'_{i}$. We also
need almost complex structures $J$ and $J'$. We assume that we are in a generic
case and our purpose is to construct a comparison morphism:
$$\Phi:\FF C(L_{0},L_{1},\eta; J,(f_{i},g_{i}))\to \FF C(L'_{0},L'_{1},\eta'; J',(f'_{i},g'_{i}))~.~$$
The construction of this morphism is quite similar to the construction of the comparison
morphism $\phi^{F}$ in the case of the cluster complex (see \S\ref{subsubsec:cluster_cob}).
The second step will be to use such morphisms to show the desired isomorphism.

\

To construct $\Phi$, we fix Morse cobordisms $(F_0:L_{0}\times [0,1]\to \R, G_0)$
between $(f_{0},g_{0})$ and $(f'_{0},g'_{0})$ 
and $(F_1:L_{1}\times [0,1]\to \R, G_1)$ between $(f_{1},g_{1})$ and $(f'_{1} \circ \psi^1,(\psi^1)^{-1}_*(g'_{1}))$
(see \ref{sec:cluster_invariance}). It is preferable here to assume that the Morse cobordisms
are in fact defined on $L_{i}\times\R$ but this is easy to arrange.
We also extend $\psi$ to $M\times \R$ so that $\psi^{s}=id$ for $s\leq 0$ and $\psi^{s}=\phi^{1}$
for $s\geq 1$. We fix a time-dependent almost complex structure $\bar J$
which interpolates between $J$ and $J'$ again in infinite time. As before we denote the generators of the
fine Floer complex associated to $(L_{0},L_{1})$ by $I_{\eta}$ and we denote by $I'_{\eta'}$
the similar generators of the complex associated to $(L'_{0},L'_{1})$.
We also fix the notation $\mathcal{R}$ for the coefficient ring of $\FF C(L_{0},L_{1})$
and $\mathcal{R}'$ for the coefficient ring of $\FF C(L'_{0},L'_{1})$.
The morphism $\Phi$ is a module morphism whose restriction to the coefficient
ring $\mathcal{R}$ is the morphism
 $\Phi:\mathcal{R}\to \mathcal{R}'$ induced
by $\phi^{F_0}\otimes\bar{\phi}^{F_1}$. Here $\phi^{F_0}$ is induced by $F_0$ and $\bar{J}$
 and is defined as in (\ref{eq:comparison_morph}). We need to be
a bit more precise about the morphism $\bar{\phi}^{F_1}$.  We consider the map
$$\bar{F_1}=F_1\circ ((\psi^{t})^{-1}, t):\bar L_{1}=\{(\psi^{t}(L_{1}),t)\}\subset M\times\R\to \R$$
 and define
the comparison morphism $\bar{\phi}^{F_1}$ with respect to this Morse cobordism and
the almost complex structure  $\bar{J}$.

\

On the generators of the fine Floer complex, $I_{\eta}$, $\Phi$ has the form:
$$\Phi(a)=\sum v^{a}_{x_{1},\ldots,x_{k},y_{1},\ldots, y_{s};b}(\la)x_{1}\ldots x_{k}y_{1}\ldots y_{s}be^{\la}$$
where $x_{i}\in \Crit(f'_{0})$, $y_{j}\in \Crit(f'_{1})$, $a\in I_{\eta}$, $b\in I'_{\eta'}$, $\la\in \bar\La$.

The coefficients $v^{a}_{\ldots;b}(\la)$ are again obtained by counting the elements
of some appropriate moduli spaces. For $r_{i}\in \Crit(f'_{0})\cup \Crit(f'_{1})$ and with the
rest of the notation as before we denote these moduli spaces
by: $$\Vv^{a}_{r_{1},\ldots,r_{k};b}(\la)~.~$$
They are obtained by an obvious adaptation of the definition
of the moduli spaces  $\Nn$ (from \S\ref{sec:cluster_invariance}) following the exact same
procedure that gave the moduli spaces $\Ww$  (in \S \ref{sec:fine_moduli})
when starting with the moduli spaces $\Mm$: in short, the root disk is replaced
by a pseudoholomorphic $\bar J$-strip $u$ which verifies instead of (\ref{eq:Flo_eq})
the moving boundary condition $$u_{v_{0}}(\R,0)\subset L_{0} \ , \ u_{v_{0}}(s,1)\in \psi^{s}(L_{1})$$
and, in the analogue of property iv'', (of \S\ref{sec:fine_moduli}) we use the gradient flow of
$F_0$ on $L_{0}\times \R$ and the gradient flow of $F_1$ on $\{(\psi^{t}(L_{1}),t)\}\subset M\times\R$.
Verifying that $\Phi$ is a chain morphism is now a simple exercise by the same arguments as in
\S\ref{sec:cluster_invariance}.

\

We will now use this construction to show the existence of an isomorphism between
the two fine Floer complexes.

First we assume that the Hamiltonian diffeomorphism $\psi$ is constant equal to the identity.
Therefore $I_{\eta}=I'_{\eta'}$.
We introduce a word-area filtration as in \S\ref{subsec:area_filtration} (with the obvious
adaptation that the area of the strips involved is also taken into account).
In the associated spectral sequence the term $E^{1}$ is an isomorphism.
This happens because the only connecting orbits that are taken into account at the level
$E^{0}$ are the Morse orbits and the constant strips joining a point in $I_{\eta}\times -\infty$
with the same point in $I_{\eta}\times +\infty$. Applying the same type of algebraic argument
as in \S\ref{subsec:area_filtration} we deduce in this case that $H_{\ast}(\Phi)$ is an isomorphism.

Now we fix $F_0$ and $F_1$ (of course, this means that also
$(f_{i},g_{i})$ as well as $(f'_{i},g'_{i})$ are fixed) and we also pick $\bar{J}$ (which is
in general time-dependent) so that $\psi^{s}_{\ast}(\bar{J}_{s,0})=\bar{J}_{s,1}$. This means that
$\bar{\phi}^{F_1}=\phi^{F_1}$.  We also assume that $L_{i},L'_{i}$, $i\in \{0,1\}$ are fixed.
Clearly, the resulting chain morphism $\Phi$ continues
to depend on the Hamiltonian diffeomorphims $\psi$, on $\bar{J}$ and on various other
data which we summarize by $\nu$. We denote this chain morphism
by $\Phi=\Phi(\psi,\bar{J},\nu)$. We claim that any two such chain morphisms are chain homotopic.
The key reason for this is that if $\Phi'$ is a second such morphism,
then $\Phi|_{\mathcal{R}}=\Phi'|_{\mathcal{R}}$. To simplify notation,
we denote $k=\Phi|_{\mathcal{R}}$.
Now, as usual in Floer theory, we may repeat the construction of the moduli spaces $\Vv$
with one additional parameter: the construction is simply the analogue in our setting (with our
big coefficient rings) of the usual proof that the comparison morphisms
in Floer theory induce a canonical map
in homology. This provides a linear map
$$\xi: \Q_{\ast}<I_{\eta}>\to (\mathcal{R}'\otimes\Q<I'_{\eta'}>)_{\ast+1}$$
so that if we now define $\xi :\FF C(L_{0},L_{1})\to \FF C_{\ast+1}(L'_{0},L'_{1})$
by the formula $\xi (r\otimes a)=k(r)\xi(a)$ this provides a chain homotopy between
$\Phi$ and $\Phi'$.

\

Finally, by using the inverse Hamiltonian
diffeomorphism we obtain that the composition $\Phi(\psi^{-1},\bar{J}',\nu')\circ \Phi(\psi,\bar J,\nu)$
is chain homotopic to a chain morphism of the form $\Phi(id, \bar{J''},\nu'')$
but we have seen before that such a chain morphism induces an isomorphism in homology and this ends the proof.

\end{proof}

We formulate the corresponding statement for the symmetric fine Floer homology as a Corollary
as its proof follows exactly the same scheme.

\begin{corollary}\label{cor:sym_fine} The isomorphism type of the
 symmetric fine Floer homology $$(\hat\FF H(L,H,J, (f,g)), d_{\hat F})$$
 defined in \S\ref{ss:cluster_fine_relation} does
not depend on the generic choice of $H, J, (f,g)$.
\end{corollary}

\begin{remark}\label{rem:inv_symm_fine}{\rm
In the symmetric case, it is be useful to make explicit the
following situation. Assume that we have two Hamiltonians $H,H':M\times
S^{1}\to \R$. Fix the almost complex structure $\bar{J}$ which
interpolates between $J$ and $J'$ (we might even take $J=J'$ and
$\bar{J}_{s}=J$, $\forall s$) as well as the Morse-Smale pair
$(f,g)$. Generically, we have two symmetric fine Floer complexes:
$$\hat\FF C(L,H,J, (f,g)) \ , \ \hat\FF C(L,H',J', (f,g))~.~$$
Clearly, the coefficient rings $\hat{\mathcal{R}}$,
$\hat{\mathcal{R}}'$ are just the cluster algebras $\Cl(L,J,
(f,g))$, $\Cl(L,J', (f,g))$.  Fix also a trivial Morse cobordism
between $(f,g)$ and itself. As in \S\ref{sec:cluster_invariance}
this data provides a chain morphism $\phi:\Cl(L,J, (f,g))\to
\Cl(L,J', (f,g))$ (it is easy to see that this is even an
isomorphism).  The argument in the proof above shows that, in
conjunction with this fixed data, any two homotopies $G$, $G'$ between $H$ and
$H'$ provide a chain morphisms:
$$\Phi(G), \Phi(G'):\hat\FF C(L,H,J, (f,g)) \to  \hat\FF C(L,H',J', (f,g))$$
 which are chain homotopic and, thus, induce the same morphism at the level 
 of the symmetric, fine Floer homology.
 Moreover, this morphism is an isomorphism.}\end{remark}

\section{Some applications and comments.}\label{sec:applic}

It is natural to expect that the symmetric fine Floer homology has
a simple expression when the Hamiltonian $H$ is a very small
time-independent Morse function. We start in this section by
discussing this point.

\subsection{The symmetric fine Floer
complex and the cluster complex} \label{se:relation_cluster_fine}

Assume that, as before, we have a fixed Morse-Smale pair on $L$
$(f,g)$ and an almost complex structure $J$ so that the cluster
complex $\Cl(L,(f,g),J)$ is defined.

Let $h:L\to \R$ be another Morse function so that the pair $(h,g)$
is Morse-Smale. By using $h$ we may construct a complex
$\hat{\mathcal{C}}(L,J,f,h,g)$ similar to the fine Floer complex
except that the intersection points $I_{\eta}$ are replaced by the
critical points of $h$ and the strips at i''. in equation
(\ref{eq:Flo_eq}) are replaced by linear cluster elements
corresponding to $L,J,(h,g)$, see
\S~\ref{subsec:relation_floer}. While this construction is quite
similar to many of our other constructions discussed before, one  
significant point is worth explicit notice: in the moduli spaces
which provide the differential in this new complex the cluster
elements of $f$ can only be anchored at points on the disks
appearing in the linear cluster element of $h$ (and not on points
appearing on the flow lines of $h$). Clearly, all bubbling is
dealt with by using cluster elements of $f$. Given these remarks,
the definition of the differential of this complex is 
transparent as well as the facts  that its square vanishes 
and that its homology is independent of
the choice of $J,h,f,g$.

Moreover, it is not difficult to define a comparison morphism
relating the symmetric fine Floer complex to
$\hat{\mathcal{C}}(L,J,f   ,h,g)$ as well as a morphism in the other
direction. The construction of these morphisms is based on the PSS
\cite{PSS} idea which was already mentioned in
\S\ref{subsec:relation_floer}. 
More precisely, to define these morphisms, we fix a 
bump function $\beta:\R\to [0,1]$ (as, for example, in Schwarz \cite{Schw}) so that
$\beta(s)=0$ for $s\leq 1/2$, $\beta(s)=1$ for $s\geq 1$ and $\beta$ is increasing.
We consider moduli spaces which, up to codimension two elements, are made
of objects consisting of: a chain 
of $J$-disks $u_{1},\ldots, u_{k}$ so that a critical point $x$ of $h$
is connected by a negative gradient flow line of $h$ to $u_{1}$ and each disk $u_{i}$ is related
to the disk $u_{i+1}$ by a negative gradient flow line of $h$; 
a finite energy solution $u$ of the equation 
$$\partial_{s}u+\partial_{t}u+\beta(s)\nabla H=0$$
where 
$$
u: (\Si, \p \Si) \to (M, L)
$$
is defined on $\Sigma=\R\times [0,1]$; we ask that 
a point on the disk $u_{k}$ be related to $u(-\infty)$ by a negative gradient flow line
of $h$ and that $\lim_{x \to \infty} u = \ga$ where $\ga$ is a trajectory of $H$; 
$f$-clusters possibly attached to the disks $u_{i}$ as well as
to $u$. Under the usual genericity conditions and by the same arguments
used in describing the comparison morphisms for the fine Floer homology,
 it is easy to see that by counting the elements of such $0$-dimensional 
 moduli spaces, we get a morphism:  
$$\hat{\mathcal{C}}(L,J,f   ,h,g)\to \hat\FF C(L,H,J,(f,g)~.~$$
The morphism in the other direction is constructed in a similar way.
Both morphisms are $\Cl(L,(f,g),J)$-morphisms. By treating both of these maps as
module morphisms over the ring $\hat{\mathcal{R}}=\Cl(L,(f,g),J)$, it is not difficult
to show, along the lines in
(\ref{rem:inv_symm_fine}) and by making use of the PSS method, that the composition 
$$\hat{\mathcal{C}}(L,J,f   ,h,g)\to \hat\FF C(L,H,J,(f,g)\to \hat{\mathcal{C}}(L,J,f   ,h,g)$$
is chain
homotopic to the identity.

Thus we obtain:

\begin{proposition}\label{prop:equivalence}
With the notation above, we have:
$$\hat{\FF} H_{\ast}(L)\simeq H_{\ast}\hat{\mathcal{C}}(L,J, f,h,g)~.~$$
\end{proposition}

\

The construction of $\hat{\mathcal{C}}(L,J,f,h,g)$ is possible
even if $h=f$. A simple exercise shows that, in this case,
$\hat{\mathcal{C}}(L,J,f,f,g)$ coincides with
$s\widetilde{\Cl}(L,J;(f,g))$, which provides the proof of
Proposition \ref{prop:cluster_fine} in
\S\ref{subsubsec:relation_cluster_fine} (the suspension $s$
appears here because, in $\mathcal{C}(L,J,f,h,g)$, the degrees of the
generators represented by the critical points of $h$ are the same
as their Morse indices).

\subsection{The cluster complex in minimal form}\label{subsec:min}
As we shall see further in this section, it is useful to reduce algebraically
the cluster complex to a minimal form that we now describe.

\begin{proposition} \label{prop:min} Assume that $L$, $(f,g)$, $J$ are so that the
cluster complex  $$\Cl(L,(f,g),J)$$ is defined. 
There exists a commutative differential graded algebra 
$\Cl_{min}(L, (f,g),J)$ and a surjective differential algebra morphism
$$P: \Cl(L,(f,g),J)\to \Cl_{min}(L, (f,g),J)$$ so that, 
as an algebra, $\Cl_{min}(L, (f,g),J)=(S (s^{-1}\ H_\ast(L;\Q))
\otimes\ \Lambda)^{\wedge}$, the differential $d_{min}$ of $\Cl_{min}(L,(f,g),J)$
increases strictly the word-area filtration and $P$ induces an isomorphism in homology. 
\end{proposition}

\begin{remark} {\rm If $f$ happens to be a perfect Morse function and thus its
Morse differential vanishes, then $\Cl(L, (f,g),J)$ itself may be taken
as $\Cl_{min}(L,f,J)$.}
\end{remark}

\begin{proof} Let $d$ be the differential of the cluster complex $\Cl(L,f,J)$.
We decompose $d=d_{0}+d'$ where $d_{0}$ is induced by the Morse differential
(in other words, this is the part that does not increase the word-area filtration).
By a  change of basis inside $\Q<\Crit (f)>$ we may assume that
each generator $x$ verifies either $d_{0}x=0$ or $d_{0}x=y$ with $y$ another
generator of the complex. We fix such a pair of generators $x, y=d_{0}x$
and consider the element $y'=dx$. We want to replace, in the cluster
complex,  the generator $y$ by the element $y'$. For this
we only need to show that $y$ can be expressed in terms of the other generators
of $\Cl(L,f)$ together with $y'$. Write $y'=y+b$. Certainly, $b$ is of higher
word-area filtration than $0$. We now express $y=y'-b$. For further use
we put $u_{0}=y'$ and $-b=r_{0}$ so that $y=u_{0}+r_{0}$. The term $b$ might still
contain some $y$'s and we write $b=a+y^{i}b'+c$ where $a$ and $b'$ do not contain $y$ and
$c$ might contain $y$ but is of higher word-area filtration than $y^{i}b'$. 
Replace $y$ by $y'-b$ in the expression of $b$. This gives $y=y'-a+(y'-b)^{i}b'+c$.
This means that we have now written $y$ as a sum $y=u_{1}+r_{1}$ 
so that $u_{1}$ is only expressed in $y'$ and in the generators of $\Cl$ which are different
from $y$ and $r_{1}$ has higher word-area filtration than $r_{0}$. We apply iteratively 
this process and, with respect to our completion, this does express $y$ as the limit of the $u_{i}$'s. 

Now let $I(x,y')$ be the differential ideal generated by $x$ and $y'=dx$. 
We consider the projection $$P_{1}:\Cl(L,f,J)\to \Cl(L,f,J)/I=\Cl_{1}~.~$$
We notice that $\Cl_{1}$ is a free commutative differential algebra generated 
by the same generators as $\Cl(L,f,J)$ with the exception of $x$ and $y'$.
Moreover, it is a simple exercise to see that the ideal $I(x,y')$ is acyclic 
so that $P_{1}$ induces an isomorphism in homology.

This process has reduced the number of pairs of generators $x,y$ of our algebra
which are related by the formula $d_{0}x=y$  by one. We iterate it till no
such pairs remain. The resulting algebra is $\Cl_{min}$. 
\end{proof}

\begin{remark}{\rm 
We will always apply this proposition to functions $f$
with a single local minimum and a single local maximum. These two generators are always $d_{0}$-cycles
and they are not pertubed by the algebraic reduction decribed above.}
\end{remark}

As an application of this proposition we complete the proof of 
Corollary \ref{cor:cond_uni_ruled}. Indeed, in the construction of the symmetric fine Floer
homology we may very well use $\Cl_{min}(L,f,J)$ instead of $\Cl(L,f,J)$: it is simply a
matter of changing coefficients by using the map $P$. The advantage is that
if $H_{\ast}(L;\Q)$ vanishes in even degrees except in dimension $0$ and $\dim(L)$,
then $\Cl_{min}(L,f,J)$ can not have high free terms. Moreover, for any function $f$
with a unique minimum $m$, if $d_{min}m\not=0$, then $dm\not=0$. The argument
used to prove the Corollary in the special case when a perfect Morse function exists can therefore be
applied also in this case and this concludes the proof of \ref{cor:cond_uni_ruled}.

\subsection{Detection of disks of bounded
area}\label{subsec:bded_disks}

Recall that $L$ is said to be displaceable if there exists a
Hamiltonian $H:M\times S^{1}\to \R$ whose associated Hamiltonian
diffeomorphism $\psi_{H}$ verifies $\psi_{H}^{1}(L)\cap
L=\emptyset$. The displacement energy of $L$, $E(L)$, is then the
infimum of the Hofer energies of the Hamiltonians with this
property. The Hofer energy of a Hamiltonian $H$ will be taken here
simply as $\max_{x,t}H-\min_{x,t}H$.

\

The purpose of this subsection is to prove Corollary~\ref{cor:disks} from \S~\ref{subsubsec:GS-problem}.
For the convenience of the reader we recall that this corollary claims that 
for a relatively spin, orientable  Lagrangian submanifold $L$ which is displaceable, 
any $\omega$-tame almost complex structure $J$
has the property that one of the following is true:
\begin{itemize}
\item[i.]  for any point $x\in L$  there exists a
$J$-pseudo-holomorphic disk of symplectic area at most $E(L)$
whose boundary rests on $L$ and which passes through $x$.
\item[ii.] there exists a $J$-disk of Maslov index at most 
$$2-\min\{2k\in\N^{\ast}\backslash \dim(L) : H_{2k}(L;\Q)\not=0\}$$
and of symplectic area at most $E(L)$.
\end{itemize}

\begin{proof} 
We start with a number of more general constructions.

Recall from Remark \ref{rem:action}
that it is possible to define an action in the setting of the fine
Floer complex. Of course, the definition given in that remark also
applies to the symmetric fine complex and it is this version that
will be of use here. 
Clearly, there is also an action in the complex $\hat{\mathcal{C}}(L,J,f,h)$:
 the action of a generator $\bar{x}e^{\la}$ in this complex
is $h(x)-\omega(\la)$. Both functions $f$ and $h$ used here are  assumed to have
a single local minimum and a single local maximum. We fix a small disk $D$ in $L$ and we assume
that both the minimum of $f$ and that of $h$ are included in this disk. We indend
to show that either there is a $J$-disk passing through $D$ and of area $\leq E(L)+\epsilon$
or there exists a $J$-disk of non-positive Maslov index with area again bounded by $E(L)+\epsilon$
where $\epsilon$ is an arbitrarily small constant and whose Maslov index satisfies the condition
ii. above. As the disk $D$ is arbitrary this 
implies the claim by Gromov compactness.

The proof is based on a comparison between
three symmetric, fine Floer complexes. For all three of them we
will use the same data $((f,g),J)$. The three Hamiltonians will be
$h:M\to \R$, $H:M\times S^{1}\to \R$ and $h':M\to \R$ so that: $h$
is a Morse function with a single minimum $\bar{m}$,
$\epsilon> \max(h)-\min(h)$, $\min_{x,t}H>0$ and
$K\geq\max(h)>\min(h)>\max_{x,t}(H)$, $h'= h-\max(h)$ where $K$ is some
positive constant and $\epsilon$ is small and $h'$ is very $C^{2}$-small.
The key remark that we will use is that, as in the usual Floer case,
if we compare the different Hamiltonians by using monotone homotopies,
then the action can not increase between the domain and the image of a map.
In short, if one moduli spaces $\Ww^{a}_{\ldots;b}(\la)$ is not null, then
$\mathcal{A}(a)-\mathcal{A}(be^{\la})\geq 0$ and similarly for the moduli spaces
of type $\Vv$ which are associated to the homotopies relating the different
Hamiltonians in case these homotopies are monotone. Our choice
of Hamiltonians imply the existence of such monotone homotopies and as a result
we obtain morphisms:

$$\hat{\FF}C(L,h)\stackrel{\Phi}{\longrightarrow}
\hat{\FF}C(L,H)\stackrel{\Phi'}{\longrightarrow}\hat{\FF}C(L,h')$$
so that both preserve the action level.

As discussed in \S\ref{se:relation_cluster_fine} we may replace
$\hat{\FF}C(L,h)$ with $\hat{\mathcal{C}}(L,h)=\hat{\mathcal{C}}(L,J,h,f,g)$
because it is easy to see that the comparison morphisms described
there are also action preserving. We have a similar identification between
$\hat{\FF}C(L,h')$ and $\hat{\mathcal{C}}(L,h')$. The composition
$\Phi'\circ\Phi$ is chain homotopic - see Remark \ref{rem:inv_symm_fine}
to a morphism induced by a trivial Morse cobordism $G$
between $h$ and $h'$ (which has the form $G(x,t)=h(x)+k(t)$). Moreover, this chain homotopy
is also action preserving.

For the complex $\hat{\FF}C(L,H)$ we may define its action filtration as follows:\\
$\hat{\FF}C^{s}(L,H)$ is the (appropriate completion of the) rational vector
space generated by the  monomials $x_{1}\ldots x_{k}ae^{\la}$ where
$\mathcal{A}(\bar{a}e^{\la})\leq s$
(here, as before, $x_{i}\in \Crit(f)$, $a$ is a generator of the fine symmetric complex,
and $\la\in\bar{\La}$). We have similar filtrations also for $\hat{\mathcal{C}}(L,h)$
and $\hat{\mathcal{C}}(L,h')$ (clearly, the stages in these filtrations are not
anymore $\bar{\La}$-modules but this will not affect in any way the arguments
that follow).

These filtrations are differential and for $s'< s$
we may define the chain complex: $$\hat{\FF}C^{s,s'}(L,h)=\hat{\FF}C^{s}(L,h)/\hat{\FF}C^{s'}(L,h)$$
and similarly for the two other complexes.
Due to the fact that $\Phi$ and $\Phi'$ are action preserving the following
composition is well defined:
$$\Psi=\Phi'\circ\Phi: \hat{\mathcal{C}}^{K,-\epsilon}(L,h)\longrightarrow
\hat{\FF}C^{K,-\epsilon}(L,H)\longrightarrow\hat{\mathcal{C}}^{K,-\epsilon}(L,h')$$
and as the chain homotopy between $\Phi'\circ\Phi$ is action preserving,
$\Psi$ is chain homotopic with the chain morphism $\Psi'$ induced by the trivial
Morse cobordims $G$. In the argument below we will only use the differential of these
truncated complexes.

We now specialize our construction: we notice that we
may pick the constant $K$, the Hamiltonian $H$ and the Morse
function $h$ so that $\phi_{1}^{H}(L)\cap L=\emptyset$ and $K=E(L)+\epsilon$.
In this case $\Psi=0$.

We also need to fix the monotone homotopies $L:h\simeq H$ and $L':H\simeq h'$ which induce
respectively the $\phi$ and $\phi'$ as well as a monotone homotopy of homotopies $T:L\# L'\simeq G$
(here $L\# L'$ is the composed homotopy which joins $h$ to $h'$). The chain homotopy $\eta^{T}$ between
$\Psi$ and $\Psi'$ is provided by this $T$ - see Remark \ref{rem:inv_symm_fine}.

Let $m'$ be the unique minimum of $h'$. We will add bars over the critical points
of $h$ or $h'$ to denote the respective generators of the fine Floer complexes.
Clearly, if we have $d\bar{m}\not=0$ or $d\bar{m}'\not=0$
by the same argument as in Corollary \ref{cor:uni_ruled} we obtain that there is a $J$-disk passing
through $m$ whose area is at most $E(L)+\epsilon$ which consludes the proof. So from now on
we shall assume that $d\bar{m}=0=d\bar{m}'$. 

Notice also that we may assume that $\Psi'(\bar{m})=\bar{m}'$. Indeed, we have that $\Psi'(\bar{m})=\bar{m}'+\ldots$
because there is a unique flow line of $G$ exiting $\bar{m}$ and this flow line ends in $\bar{m}'$. Geometrically,
this flow line is constant equal to the critical point $m$. Therefore, if there is any cluster tree
originating at $\bar{m}$ there will be a $J$-disk which crosses this line - but this means
that this disk passes through $m$ which again implies the claim.

In view of this, we have $d\eta^{T}(\bar{m})=\bar{m}'$. Given Remark \ref{rem:cl_moduli} the only
possibility for $\bar{m}'$ to be a boundary is that there exists some $x\in Crit(f)$
so that $dx=a_{0}e^{\la}+\ldots$ (we may assume $|x|\geq 0$ otherwise, again, our claim follows).

But this means precisely that the truncated cluster complex has high free terms which shows
that there exist $J$-disks of area bounded by $E(L)+\epsilon$
and whose Maslov index is at most $2-\min\{2k\in\N^{\ast}\backslash \dim(L): \Crit_{2k}(f)\not=0\}$.
The last step in the proof is to note that, in this argument, instead of the cluster 
complex associated to $f$ we could have used as well the minimal form of it, $\Cl_{min}(L,f,J)$, as described in 
\S\ref{subsec:min}. This allows us to replace the condition $\Crit_{2k}(f)\not=0$ by the one 
in the statement.
\end{proof}

\begin{remark}\label{rem:method}{\rm
The method of proof used here is an adaptation of one used in \cite{BaCo} and Barraud-Cornea \cite{BaCo2}
to detect pseudoholomorphic strips of bounded area. 
}\end{remark}

A particular case is worth stating separately.

\begin{corollary} \label{cor:hom_cond}If a relatively spin, orientable Lagrangian $L$ verifies
$H_{2k}(L;\Q)=0$ if $2k\not\in\{0,\dim(L)\}$ and is displaceable, then
for any $J$, through each point of $L$ passes a $J$-disk of symplectic
area at most $E(L)$.
\end{corollary}

\ 

For any Lagrangian submanifold $L\subset (M,\omega)$ define (as in \cite{BaCo2}, see also \cite{BaCo}) its real, or relative, Gromov
radius $r(L)$ as the infimum of the positive constants $r$ so that
there is a symplectic embedding of the standard standard sphere
$(B(r),\omega_{0})\stackrel{e}{\hookrightarrow} (M,\omega)$ with the property
that $e^{-1}(L)=\R^{n}\cap B(r)$.

\begin{corollary}\label{cor:energy_radius}
If an orientable, relatively spin Lagrangian with $H_{2k}(L;\Q)=0$ for $2k\not\in\{0,\dim(L)\}$ 
is displaceable, then
$$E(L)\geq \pi r(L)^{2}/2~.~$$
\end{corollary}

This follows immediately from Corollaries~\ref{cor:disks} and \ref{cor:hom_cond} by a typical Gromov's ``minimal surface" argument
as in \cite{BaCo}. If $\mu_{min}\geq 2$, the result in Corollary~\ref{cor:energy_radius} remains true 
without the homological,
orientability or relatively spin conditions because we may apply the method in the Corollaries
above directly to the complex $\mathcal{C}(L,J,f)$ described in \S\ref{subsec:relation_floer}
and this complex is defined over $\Z/2$. 
Corollary~\ref{cor:energy_radius} is a particular
case of a conjecture formulated in \cite{BaCo2}.

\subsection{Some calculations.}\label{subsec:calc}

Our machinery easily leads to restrictions on the topology of Lagrangian submanifolds in $\C^{n}$.
We will only discuss here one such immediate example. In this case, we recover results that agree with those
obtained in \cite{Fu} (see also \cite{FOOO}). We will pursue with other examples and applications elsewhere.

\begin{example}\label{ex:s1-sn}{\rm
We assume that our relatively spin Lagrangian submanifold $L\subset M$ is diffeomorphic to $S^{1}\times S^{n-1}$ and 
$\hat\FF H_{\ast}(L)=0$ (for example, as happens in the case considered in \cite{Fu} where $M=\C^{n}$).

We indend to deduce restrictions on the image of the Maslov
index homomorphism $$\mu:\pi_{1}(L)\approx\pi_{2}(\C^{n},L)\to \Z~.~$$

\

We construct the cluster complex  and the associated complex $\widetilde{\Cl}(L,J;f)$
which computes the symmetric fine Floer homology (as in \S\ref{se:relation_cluster_fine})
by using a perfect Morse function $f$. We denote by $m$ its minimum, by $M$ its maximum, by $a$
its critical point of index $1$ and by $b$ its critical point of index $n-1$. We denote by 
$\bar{x}$ the corresponding generators of $\widetilde{\Cl}(L,J;f)$ (viewed as a module
over $\Cl(L,J,f)$).

\

The first key question in this type of calculation is whether the cluster complex has free terms or not
(see Proposition \ref{prop:free_terms}). Suppose first that it has free terms. This means
that some generator $x$ verifies $dx=a_{0}e^{\la}+\ldots$. As $\mu(\la)$ is even, it follows
that $x$ may be $m$ (and, in this case $\mu(\la)=2$) but cannot be $a$. 
It might also be $b$ in which case $n-1$ is even and
 $\mu(\la)=3-n$. In all cases, it cannot be $M$ by Remark \ref{rem:free_Max}.

\

The other possibility to be considered is if the cluster differential does not have any free terms.
In this case we will use the fact that the complex $\widetilde{\Cl}(L,J;f)$ has a trivial homology
and also has a differential that cannot decrease the word length. 

\

The remark above concerning the maximum applies as well to $\bar{M}$ so that $d\bar{M}=0$. This means
that there has to be another generator $y$ of $\widetilde{\Cl}(L,J;f)$ verifying $dy=a_{1}\bar{M}e^{\la}+\ldots$,
$a_{1}\in\Q^{\ast}$.
First, $y\not=\bar{M}$ because $\mu(\la)$ is even. In case  $y=\bar{b}$, then $\mu(\la)=2$.
Another possible choice for $y$ is  $y=\bar{m}$ and then $\mu(\la)=n+1$ which has to be even. 
Similarly $y=\bar{a}$ implies $\mu(\la)=n$ which needs to be even. 

We now summarize our results for $n$ even:
\begin{itemize}
 \item[-]if $n$ is even, then $\{2, n\}\cap Im(\mu)\not=\emptyset$. 
\end{itemize}

Clearly, if $2\in Im(\mu)$, then $n\in Im(\mu)$ also but there exists an example due
to Polterovich in which $2\not\in Im(\mu)$.

\

We now pursue our discussion in case $n$ is odd. As $2\in Im(\mu)$ implies $3-n\in Im(\mu)$, 
our discussion till now shows that
if $3-n\not\in Im(\mu)$, then the cluster complex has no free terms.
In this case we see that $d\bar{m}=a_{1}\bar{M}e^{\la}+\ldots$ and so $\mu(\la)=n+1$. 
We now write $d\bar{m}=\bar{M}e^{\la}(a_{1}+a)+ c$ where $a\in \Cl(L,f)$ and  is a term
of word-area filtration greater than $0$, $a_{1}\in \Q^{\ast}$ and $c$ does not contain $\bar{M}$. We define
a $\Cl(L,f)$-differential module epimorphism $$q:\widetilde{\Cl}(L,f)\to (\Cl(L,f)\otimes \Q<\bar{a},\bar{b}>,d')=
\mathcal{K}$$ by sending $\bar{m}\to 0$, $\bar{M}\to -c(a_{1}+a)^{-1}e^{-\la}$ where the 
differential $d'$ in the target module is the one induced by $q$ (notice that $(a_{1}+a)$ is 
invertible in $\Cl(-)$ due to our completion). It is easy to see that this map $q$ induces an isomorphism in
homology and thus $H_{\ast}(\mathcal{K})$ is trivial.  This means that either $d'\bar{a}=a_{2}be^{\la'}+\ldots$
and in this case $\mu(\la')=n-1$ and/or  $d'\bar{b}=a_{3}ae^{\la''}+\ldots$  ($a_{i}\in\Q^{\ast}$) 
which means $\mu(\la'')=3-n$.
But as we also have $\mu(\la)=n+1$ both these possibilities contradict $3-n\not\in Im(\mu)$. 
To conclude the argument in this case:

\begin{itemize}
\item[-]if $n$ is odd, then $\{2, 3-n\}\cap Im(\mu)\not=\emptyset$.
\end{itemize}
}\end{example}

\subsection{Detection of closed orbits of Hamiltonian systems}\label{subsec:orbits}

   We show here that, given any two compact oriented relatively spin Lagrangian 
   submanifolds $L,L'$ in a closed or geometrically bounded  
   symplectic manifold $(M, \om)$, there is a natural morphism of chain complexes
 $$
 \cyl^{\Phi_0}: \Cl (L)  \to   \Cl (L)  \otimes \Cl (L')  \otimes \La_{\Phi_0}
 $$
 defined using holomorphic cylinders (of varying conformal structure) with one 
 boundary in $L$ and the other in $L'$ (here $\La_{\Phi_0}$ is a coefficient 
 ring defined below). When this morphism is not trivial in homology, and if 
 some Hamiltonian in $M$ {\em separates}  $L$ and $L'$, one can easily deduce 
 the existence of closed orbits of $H$ with bounded period. Let us explain this more precisely.

 \subsubsection{Reviewing the Maslov index}
 Let $L,L'$ be two oriented relatively spin Lagrangian submanifolds in $(M, \om)$. 

To define the Maslov index of a cylinder $u \in [ S^1 \times [0,1], S^1 \times \{0\}, S^1 \times \{1\}; M, L, L']$, first choose a $SO(n)$-trivialization of $(u |_{S^1 \times \{0\}})^*(T_*L)$, extend it in the natural way to a $U(n)$-trivialization of $(u |_{S^1 \times \{0\}})^*(T_*M)$, and then to a well-defined, up to homotopy,  $U(n)$-trivialization of $u^*(T_*M)$. The Maslov index of $u$ is set to be equal to the Maslov index of the Lagrangian subbundle $(u |_{S^1 \times \{1\}})^*(T_*L')$ with respect to this trivialization.

It is easy to see that this definition is independent of the choice of the $SO(n)$-trivialization of $(u |_{S^1 \times \{0\}})^*(T_*L)$ and that it is a variant of the definition
of the Maslov morphism as defined in \S\ref{subsec:fine_definition}. Similarly, one easily sees that, in the case when the homotopy class of $u |_{S^1 \times \{0\}}$ vanishes, it is equal to the Maslov index of 
$(u |_{S^1 \times \{0\}})^*(T_*L)$ with respect to the trivialization of any capping disc $v$, minus the  Maslov index of 
$(u |_{S^1 \times \{1\}})^*(T_*L')$ with respect to the trivialization of the capping disc $u \# v$. Note finally that the index of the space of parametrized holomorphic cylinders (with unprescribed conformal structure) in a class $\Phi$ is equal to $\mu(\Phi) + 1$, so that the unparametrized index is $\mu(\Phi)$.

 Now fix any homotopy class $\Phi_0$ of continuous maps
 $$
 S^1 \times [0,1] \to M
 $$
 with the boundary component $S^1 \times \{0\}$ in $L$ and the boundary 
 component $S^1 \times \{1\}$ in $L'$. Consider  the obvious left 
 action of $\pi_2(M,L)$ on the space $[ S^1 \times [0,1], S^1 \times \{0\}, S^1 \times \{1\}; M, L, L']$ 
 and the right action of $\pi_2(M,L')$ on the same space. Obviously these two 
 actions commute. Let us denote by $K$ the subset of  $[ S^1 \times [0,1], 
 S^1 \times \{0\}, S^1 \times \{1\}; M, L, L']$ equal to the orbit of these two actions on the element $\Phi_0$.

 \begin{remark} {\rm We could as well take the full space $[ S^1 \times [0,1], S^1 \times \{0\}, S^1 \times \{1\}; M, L, L']$ instead of the orbit $K$. Note that these two sets are in general distinct (for instance for Lagrangian submanifolds in a cotangent given as sections, there in no left and right action, so that $K$ contain no multiple covering of $\Phi_0$). The latter being smaller than the former, the applications of this theory will be simpler to compute.} 
 \end{remark}

 As usual,  denote by $ \sim$ the quotient of the objects in $[ S^1 \times [0,1], S^1 \times \{0\}, S^1 \times \{1\}; M, L, L']$, $\pi_2(M,L)$ and $\pi_2(M,L')$  by the relation requiring equality of both the Maslov index and the symplectic area. Now let $\bar K$ be the abelian subgroup of $\Z \times \R$ generated by the elements in $\pi_2(M,L)/ \sim, \pi_2(M,L')/ \sim$ and $K/  \sim$, which of course is the same as taking the subgroup generated by the first two and by the singleton $(\mu(\Phi_0), \om(\Phi_0))$. Let finally $\La_{\Phi_0}$ be the usual completion of the group ring $\Q <\bar K>$.

\subsubsection{Appropriate moduli spaces} We define the morphism $ \cyl^{\Phi_0}$  by considering the following moduli spaces. Fix Morse functions and metrics $(f,g), (f',g')$ on $L$ and $L'$, a compatible complex structure $J$ and a generic perturbation $\nu$ (for simplicity, we will omit the $\nu$ in the symbols denoting our moduli spaces). For elements $x, x_1, \ldots, x_k \in \Crit(f)$, elements  $y_1, \ldots, y_m \in \Crit(f')$, and a class $\bar{\la} \in \bar K $, denote by $\Cc^x_{x_1, \ldots, x_k,y_1,\ldots, y_m} (\bar{\la})$ the space of configurations (after quotient by the automorphism group) made from

1) a cluster on $L$ rooted at $x$, with terminations $x_1, \ldots, x_k$ and where all vertices of the graph defining the cluster correspond to boundaries of $(J, \nu)$-holomorphic disks in $(M,L)$ except precisely one which is replaced by a loop $\ga_L$ in $L$;

2) a $(J, \nu)$-holomorphic cylinder in class $K$
$$
u: (S^1 \times [0,\si], S^1 \times \{\si\}, S^1 \times \{0\})  \to  (M, L, L')
$$
satisfying
$$
u |_{S^1 \times \{\si\}} = \ga_L, 
$$
where $\si > 0$ is the height of the cylinder defining its conformal structure (the cylinder is endowed with the  standard complex structure coming from the quotient of the strip $[0,1] \times [0, \si] \subset \C$ by the translation by $1$ in the $x$-direction); 

3) a cluster on $L'$ rooted at a vertex replaced by the loop $\ga_{L'} :=  u |_{S^1 \times \{0\}}$, with terminations $y_1, \ldots, y_m$ in which all other vertices of the graph defining the cluster correspond to boundaries of $(J, \nu)$-holomorphic disks in $(M,L')$.

   We define the Maslov index of such a configuration $C$ and its area by:
$$
\sum_{\Dd} (\mu(D), \om(D)) + (\mu(u), \om (u)) 
$$
where $\Dd$ is the set of $(J, \nu)$-holomorphic discs appearing in the cluster elements of the configuration (in $(M,L)$ as well as in $(M,L')$) and $u$ is the above cylinder. We finally require that

4)  $(\mu(C), \om(C)) = \bar{\la}$.

\begin{lemma} a) The dimension of $\Cc^x_{x_1, \ldots, x_k,y_1,\ldots, y_m} (\bar{\la})$ is equal to
$$
|x| - \sum_i |x_i| - \sum_j |y_j| + \mu(\bar{\la}) -n+2.
$$

   b) The cell of maximal dimension of the boundary of the compactification
   $$
\bar{\Cc}^x_{x_1, \ldots, x_k,y_1,\ldots, y_m} (\bar{\la}) - {\Cc}^x_{x_1, \ldots, x_k,y_1,\ldots, y_m} (\bar{\la})   
   $$ 
   is made of all configurations obtained from the ones described before by admitting one broken flowline of one of the flows $- \nabla_g f, - \nabla_{g'} f'$.  The break point belongs to $\Crit(f) \cup \Crit(f')$. The full boundary is obtained in a similar way by adding $d_1$ broken flowlines and $d_2$ $(J, \nu)$-holomorphic $2$-spheres attached to points belonging to discs or to the cylinder, where $d_1 + 2d_2 \le d$ (here $d$ is the dimension of the moduli space). 
   \end{lemma}

     The proof of this lemma is left to the reader. Two points are worth mentioning:  the dimension formula in (a) is an easy consequence of the formula for the real dimension of parametrized pseudoholomorphic cylinders $u$ with one boundary on $L$ and the other on $L'$ (of non-specified conformal structure): it is equal to $\mu(u) + 1$, see \cite{GL} (thus the unparametrized  moduli space has dimension $\mu(u)$).  The second remark concerns (b):  a cylinder in class $\Phi$ may degenerate to a cylinder in class $\Phi - \la - \la'$ where $\la \in \pi_2(M,L)/ \sim$ and $\la' \in \pi_2(M,L')/ \sim$ together with
$J$-disks with boundaries in $L$ and $L'$ and of classes, respectively, $\la$ and $\la'$ . This shows that the correct Novikov ring to consider is indeed the larger ring $\La_{\Phi_0}$ which contains both $\La(L)$ and $\La(L')$.  \QED

\

Set:
$$
\cyl^{\Phi_0}: \Cl (L) \to \Cl (L) \otimes \Cl (L') \otimes \La_{\Phi_0}
$$
defined by 
$$
\cyl^{\Phi_0} (x) = (-1)^{|x|} \sum_{x_1, \ldots, x_k, y_1, \ldots, y_m, \bar{\la}} \# \Cc^x_{x_1, \ldots, x_k,y_1,\ldots, y_m} (\bar{\la}) x_1, \ldots, x_k  y_1, \ldots, y_m e^{\bar{\la}}
$$
when the dimension of this configuration space is zero, and zero otherwise. Extend this definition by the Leibniz rule 
$$
\cyl^{\Phi_0} (xy) = \cyl^{\Phi_0} (x) \, y + (-1)^{|x|} x \, \cyl^{\Phi_0} (y)
$$ 
and then by linearity.  
We define the degree of an element in $\Cl (L) \otimes \Cl (L') \otimes \La_{\Phi_0}$ by the usual formula.

\begin{lemma}
$\cyl^{\Phi_0}$ is a morphism of chain complexes of degree $2-n$.
\end{lemma}

\proof{}  Let $x \in \Crit(f)$. We must show that $\cyl^{\Phi_0} \circ d (x) = d \circ \cyl^{\Phi_0} (x)$. A  term on the left hand side corresponds to a configuration made from the one defining a term $a=  x_1 \ldots x_k e^{\la} $ in the differential of $x$ attached at one of its terminations $x_i$ to the root $x_i$ of a configuration defining a term $b = x_1'  \ldots x_{k'}' y_1 \ldots y_m e^{\bar{\la}}$ in $\cyl^{\Phi_0} (x_i)$. Smoothing the resulting configuration at $x_i$, one gets a real one-parameter family of configurations in 
$\Cc^x_{x_1, \ldots, \hat{x}_i, \ldots, x_k, x_1', \ldots , x_{k'}',y_1,\ldots, y_m} (\la +\bar{\la})$. By the above compactification formula, the other end of this bordism must be obtained from an element of  
$\Cc^x_{x_1, \ldots, \hat{x}_i, \ldots, x_k, x_1', \ldots , x_{k'}',y_1,\ldots, y_m} (\la +\bar{\la})$ by adding one broken flowline either on $f$ or on $f'$. If the broken edge of the tree appears below the vertex corresponding to the cylinder $u$ (on $L$ or on $L'$, it does not matter), then the resulting configuration belongs to $d \circ \cyl^{\Phi_0} (x)$. Otherwise, it belongs to $\cyl^{\Phi_0} \circ d (x)$. One checks easily that the same holds when one starts from a term on the right hand side of $\cyl^{\Phi_0} \circ d (x) = d \circ \cyl^{\Phi_0} (x)$.  The statement concerning the degree is an immediate consequence of the formula for the dimension of the moduli spaces in the last lemma. \QED

 \subsubsection{Applications to the existence of closed orbits}

   We will say that a morphism  $\Cl H_* (L) \to H_*(\Cl  (L) \otimes \Cl (L') \otimes \La_{\Phi_0})
$ is {\em trivial} if it vanishes on all of the $\La$-generators of $\Cl H_* (L)$ 
(note that it need not vanish on the unit).

    \begin{corollary} Let $L,L'$ be oriented relatively spin submanifolds in some 
    geometrically bounded manifold $(M, \om)$. Let $H_{t \in [0,1]}$ be a Hamiltonian on $M$ 
    that separates $L$ and $L'$ in the sense that the gap
    $$
    \int_0^1 (\min_{L} H_t - \max_{L'} H_t) dt
    $$
    is strictly positive. Then, if the morphism $\cyl^{\Phi_0}$  is not trivial in homology, 
    there must be a closed orbit $\ga(t)$ of the system $T H_t$ of period $1$, for 
    some $T \le A/gap$ where $A$ is the symplectic area of $\Phi_0$. In particular, 
    if $H$ is autonomous, this gives the existence of a periodic orbit of $H$ of period bounded above by $A/gap$.
    \end{corollary}

    \proof{} We proceed in three simple steps: the first one gives an a priori estimate, 
    the second one is the core of the proof, and the third gives the definition of a homotopy 
    operator used in the second step. 

    \

    1) {\em First step}. \, From \cite{GL} (or Hofer-Viterbo \cite{HV} in  a 
    slightly different context), recall that there is an a priori estimate that gives an upper 
    bound of the $s$-energy of the cylinder $u$ in class $K$ satisfying the 
    equation $\bar{\p}_J u= - T \nabla_{g_M} H_t$, where $g_M$ is 
    the metric $\om(J \cdot, \cdot)$ on $M$ and $T$ is allowed to vary in the 
    semi-infinite interval $[0, \infty)$ (here $s$ is the vertical coordinate, 
    whereas $t$ will denote the $S^1$-horizontal coordinate). This estimate is
    $$
    \int_{C_{\si}} \| \frac{\p u}{\p s} \|^2 ds dt  \le A - T \, gap.
    $$
   To see this, calculate:
   \smallskip
    $$
    \begin{array}{lcc}
    \int_{C_{\si}} \| \frac{\p u}{\p s} \|^2 ds dt  & = & \\
     \int_{C_{\si}} <\frac{\p u}{\p s}, J \frac{\p u}{\p t}  - T \nabla_{g_M} H_t>  ds dt & = & \\
      \int_{C_{\si}} <\frac{\p u}{\p s}, J \frac{\p u}{\p t}> ds dt   -  
      \int_{C_{\si}} <\frac{\p u}{\p s}, T \nabla_{g_M} H_t> ds dt.
      \end{array}
    $$ 
    \smallskip
    But the first term is the symplectic area of the cylinder in class $K$, 
    bounded above by $A= \om(\Phi_0)$. The second term is the the average 
    over $t$ of the difference of $H_t (t, \si) - H_t(t,0)$, which immedialty yields the above estimate.

   \

     2)  {\em Second step}. \,  Now assume that, in homology, $\cyl_{\Phi_0}$ is 
     not trivial and that $T H_t$ has no closed orbit of period $1$ for any $T \le A/gap$. 
     Then there is  $\eps >0$ sufficiently small so that $T H_t$ has no closed orbit 
     of period $1$ for any $T \le (A/gap) + \eps$.  Consider the one-parameter 
     family of morphisms $\cyl^{\Phi_0}_T$ in which the defining equation of the cylinder is replaced by
       $$
   \bar{\p}_J u= - T \nabla_{g_M} H_t.
   $$       
For $ A/gap \le T \le (A/gap) + \eps$, the right hand side of our estimate becomes 
non-positive so that no solution exists and therefore $\cyl^{\Phi_0}_T$ obviously 
vanishes on all generators for these values of $T$. On the other hand, we prove 
in the third step below that $\cyl^{\Phi_0}$  is chain homotopic to 
$\cyl^{\Phi_0}_{T}$ for any generic $T \le (A/gap) + \eps$.  This is a contradiction.

\

3)  {\em Third step}. \,  There remains to construct the homotopy operator between 
the two morphisms $\cyl^{\Phi_0}$   and $\cyl^{\Phi_0}_{T}$ 
where $T \in [A/gap, A/gap + \eps]$ is generic. For this, consider the 
morphism of degree one more than the degree of $\cyl^{\Phi_0}$:
$$
\Psi:  \Cl (L) \to \Cl (L) \otimes \Cl (L') \otimes \La_{\Phi_0}
$$
defined on a generator $x$ by 
$$
\Psi(x) =  \sum_{\tau \in [0,T]} \cyl^{\Phi_0}_{\tau}(x)
$$
when, as usual,  the corresponding obvious moduli space $\Cc^x_{x_1, 
\ldots, x_k,y_1,\ldots, y_m} (\bar{\la}, \tau \in [0,T] )$ is of 
dimension $0$. Extend it by the Leibniz rule and by linearity. In order to show the formula
$$
\p \Psi - \Psi \p =   \cyl^{\Phi_0}_{T} -  \cyl^{\Phi_0}
$$
i.e to see that this is indeed a homotopy operator between the 
two morphisms on the right hand side, we consider a one-parameter 
moduli space of the form $\Cc^x_{x_1', \ldots, x_k',y_1',\ldots, y_m'} 
(\bar{\la}', \tau \in [0,T] )$. This space can only degenerate to:

  i) a configuration that reaches one of the two obvious boundaries corresponding to $\tau = 0, T$;

  ii) a configuration that breaks in the ordinary Morse way, for an 
  interior time $\tau \in (0, T)$, into two trees,  at a critical 
  point $\bar x  \in \Crit(f)$ or $ \bar y \in \Crit(f')$. In the 
  second case, this belongs to $\p \Psi (x)$ while in the first case,  
  this belongs to  either $\p \Psi (x)$ or  $ \Psi \p (x)$ according to 
  whether the tree rooted at $\bar x$  contains or not the cylinder.

  Note that the eventual bubbling off, from the cylinder, of a pseudoholomorphic 
  disk with boundary on $L$ or $L'$ is not a boundary point of 
  our $1$-dimensional family since the cluster setting sees this phenomenon as an internal point. 

  Finally, there is no other case than the ones in (i) and (ii) above 
  because obviously  the  degeneracy corresponding to the conformal 
  parameter of the cylinder going to infinity cannot happen since 
  there is, by assumption, no closed orbit of period $1$ of $T H_t$ for $\tau \le T$.
 \QED

   There are many examples where the hypotheses of our theorem are satisfied.

\MS   \NI
 {\bf Examples}. 

 \MS \NI
A. {\em In cotangent spaces}

 \smallskip \NI
 Note first that if one considers sections of a cotangent bundle, 
 there is no bubbling off  and therefore one does not need to orient 
 moduli spaces of holomorpic discs, so the above theorem applies as 
 well without the hypothesis that $L$ and $L'$ be orientable or 
 relatively spin. In general, if $(V,g)$ is a Riemannian manifold 
 and $\al$ is a closed $1$-form whose $g$-dual $X$ is $g$-parallel, 
 the holomorphic cylinders, with respect to the natural almost 
 complex structure on $TV$, from the zero section $L$ to the 
 graph of $\al$,  $L'$, correspond to $1$-jets of integral 
 curves $\be$ of $X$ which are also geodesics of $V$: the cylinder 
 is then given by $\{ (\be(t), s \frac{d \be}{dt}) : t \in [o,a], s \in [0,1] \}$ 
 (actually, one reparametrizes the cylinder so that its circle-base is of length $1$). 

\medskip \NI
1)  The simplest example is a generalisation from the case treated in \cite{GL} : 
one considers $V= K \times W$ where $K$ is the Klein bottle and $W$ any closed 
manifold admitting a perfect Morse function. Take $L$ the zero section 
and $L'$ the graph of the pull-back to $K \times W$ of the constant 
closed $1$-form $\al$ on $K$ which corresponds, through the 
euclidean paring, to the constant vector field $\p/\p x$ in $([0,1] \times [0,1]/R$ 
where $R$ identifies the top and bottom boundaries, and reverse the 
orientation of the left and right boundaries. On $K$ the geodesic 
whose image is $[0,1] \times \{1/2\}$ is an integral curve of 
the $g$-dual of $\al$, and it is moreover unique in its 
homotopy class. Therefore, on $L = K \times W$, there is 
a $W$-family of such geodesics. One can consider a product 
Morse function on $L=K \times W$ and the configuration 
in    $\Cc^x (L,L')(\bar{\la})$ made of a flowline from a point $x$ 
of index $\dim W$ to one element in the $W$-family of geodesics, 
with the corresponding holomorphic cylinder from $L$ to $L'$, 
and with no flowline in $L'$. These configurations are unique and give a map:
$$
 \cyl^{\Phi_0} :  \Cl H (L) \to \Cl H(L) \otimes \Cl H(L') \otimes \La (\Phi_0)
 $$
 which sends $x$ to $1 \otimes e^{(\mu(\Phi_0), \om(\Phi_0))}$. 
 Because the Morse function on $W$ is perfect, the ordinary flowlines in $L$ 
 do not contribute in the morphism. Thus $ \cyl^{\Phi_0}$ does not vanish since 
 $1 \otimes e^{(\mu(\Phi_0, \om(\Phi_0))}$ is not zero in  
 $$
 \Cl H(L) \otimes \Cl H(L') \otimes \La(\Phi_0) \simeq (S(s^{-1}H_*(K \times W)) 
 \otimes S(s^{-1}H_*(K \times W))) \otimes \La(\Phi_0).
 $$

 \MS
 \NI
 2)  The same arguments apply modus vivendi to $T^n$. Indeed, if $L$ is the zero 
 section and $L'$ the graph of a constant one-form in one of the $S^1$ direction, 
 the geodesics in that same class foliate the torus. Thus once again, the 
 map $ \cyl^{\Phi_0} $ sends the minimum $x$ of a Morse function on $T^n$ 
 to a configuration made of $x$, sitting on a unique geodesic, and of the 
 corresponding holomorphic cylinder. Hence, once again, $ \cyl^{\Phi_0} (x) = 
 1 \otimes e^{(\mu(\Phi_0), \om(\Phi_0))}$ and the rest of the argument is 
 the same. One may here also generalise this to a product of $T^n \otimes W$ where $W$ admits a perfect Morse function.

 \MS
 \NI
 3)  Let $(V,g)$ be a hyperbolic manifold, consider a closed $1$-form $\al$ and 
 set as before $L$ equal to the zero section and $L'$ the graph of $\al$.  Clearly, 
 if there is a geodesic which is an integral curve $C$ of the $g$-dual of $\al$, 
 it is unique in its homotopy class. Assume that there is a Morse function on $V$ 
 with a critical point $p$ of index $n-1$ such that its Morse differential vanishes 
 and that its stable manifold has non-zero linking number with $C$. Then the map
 $$
 \cyl^{\Phi_0} :  \Cl H (L) \to \Cl H(L) \otimes \Cl H(L') \otimes \La (\Phi_0)
 $$
sends $p$ to $1 \otimes e^{(\mu(\Phi_0), \om(\Phi_0))}$ which is non-zero in homology.

\medskip
 \NI
 B. {\em In manifolds}

 \medskip \NI
Take any of the above examples imbedded as a Lagrangian submanifold of a geometrically 
bounded symplectic manifold $(M, \om)$. Assume that neither $L$ or $L'$ has 
free terms (i.e their cluster homology does not vanish). Fix a class $\Phi_0$ 
whose area is no greater than the smallest non-constant $J$-holomorphic disc 
in $(M,L)$ or $(M,L')$. Then the same maps as in (A), (B) and (C) send some element  to a cycle of the form
$$
b=1 \otimes e^{(\mu(\Phi_0, \om(\Phi_0))} + \ldots.
$$
But by hypothesis, the class $\Phi_0$  cannot degenerate so that, in the above 
expression, the only free term is $1 \otimes e^{(\mu(\Phi_0, \om(\Phi_0))}$. Thus $b$ does not vanish in homology.

\

{\bf Concluding remark.} \, The morphism $\cyl$ is the first term of a series that 
defines a differential on the tensor product $\Cl (L) \otimes \Cl (L') \otimes \La (\Phi_0)$.  
To explain this, first denote by $S_{k+1}$ the standard sphere with $k+1$ 
disjoint open disks $D_0, D_1, \ldots, D_k$ removed,  and by $\Jj_{k+1}$ 
the space of conformal structures on $S_{k+1}$. For each integer $\ell > 0$, define:
$$
d^{\ell \Phi_0}_\ell :  \Cl (L) \to \Cl (L) \otimes \Cl (L') \otimes \La (\Phi_0)
$$
by the formula
$$
d^{\ell \Phi_0}_\ell (x) = \sum_{x_1, \ldots, x_k, y_1, \ldots, y_m, \bar{\la}} 
\# \Cc^x_{x_1, \ldots, x_k,y_1,\ldots, y_m} (\bar{\la}, \ell) x_1, \ldots, x_k  y_1, \ldots, y_m e^{\bar{\la}}
$$
when the dimension of this configuration space is zero, and zero otherwise. 
Extend this definition by the Leibniz rule. Here the moduli space of 
configurations $\Cc^x_{x_1, \ldots, x_k,y_1,\ldots, y_m} (\bar{\la}, \ell)$ is 
defined like $\Cc^x_{x_1, \ldots, x_k,y_1,\ldots, y_m} (\bar{\la})$ except 
that the unique cylinder in the definition of the latter moduli space is 
replaced by $q$ surfaces $(\Si_1, j_1), \ldots , (\Si_q, j_q)$ with each $\Si_i$ 
equal to $S_{k_i +1}$ and $j_i \in  \Jj_{k_i+1}$ so that $\sum_i k_i = \ell$ and 
that each $k_i$ is larger or equal to $1$. More precisely, for 
elements $x, x_1, \ldots, x_k \in \Crit(f)$, elements 
 $y_1, \ldots, y_m \in \Crit(f')$, and a class $\bar{\la} \in \bar K $, 
 consider the space $\Cc^x_{x_1, \ldots, x_k,y_1,\ldots, y_m} (\bar{\la}, \ell)$ of configurations made from

\

1) a cluster on $L$ rooted at $x$, with terminations $x_1, \ldots, x_k$ and where all 
vertices of the tree defining the cluster correspond to 
boundaries of $(J, \nu)$-holomorphic disks in $(M,L)$ except precisely $q$ of 
them which are replaced by  loops $\ga_i, 1 \le i \le q,$ in $L$;

2) $q$ maps
$$
u_i: (\Si_i, \p D_0,\p D_1 \cup \ldots \cup \partial D_k)   \to  (M, L, L'), \;  \; 1 \le i \le q,
$$
satisfying
$$
\bar{\p}_{(j_i, J)} u_i = \nu
$$
and
$$
(u_i) |_{\partial D_0} = \ga_i;
$$ 

3) $\ell$ clusters on $L'$, each one being rooted at a vertex corresponding 
to a loop $\ga'_i, 1 \le i \le \ell$, on $L'$ so that these $\ga'_i$'s are equal to
the restrictions of the $u_i$'s 
to the various boundary circles $\p D_j, 1 \le j \le k_i$ via a one-to-one correspondence.  
In each such cluster,  
all other vertices (i.e all but the root) of the tree defining the cluster
correspond to boundaries of $(J, \nu)$-holomorphic disks in $(M,L')$. The set 
of all terminations of the these $\ell$ clusters consists of the ordered set $y_1, \ldots, y_m$.

   We finally require that the total Malov index (whose computation can be reduced 
   to the computation of a finite (connected) sum of cylinders) be equal to $\bar{\la}$, 
   and that, topologically, the $\ell$ cylinders, be in class $\ell \Phi_0$ up to the left 
   and right actions of $\pi_2(M,L)$ and $\pi_2(M,L')$.

   \

   Thus, for $\ell = 0$, this is simply the ordinary differential of the 
   complex $\Cl (L)$; for $\ell = 1$, it is, up to sign, equal to the 
   morphism $\cyl^{\Phi_0}$ above. For $\ell = 2$, the definition takes 
   into account pairs of pants relating $L$ to $L'$. Now extend each $d^{\Phi_0}_\ell$ to
$ \Cl (L) \otimes \Cl (L') \otimes \La (\Phi_0)$ by defining it as the usual 
cluster differential on the $\Cl(L')$-term,  and denote the resulting map 
by the same symbol. Finally consider the {\em big differential}
   $$
   d_{\infty}^{\Phi_0} = \sum_\ell d^{\Phi_0}_\ell : \Cl (L) 
   \otimes \Cl (L') \otimes \La (\Phi_0) \to \Cl (L) \otimes \Cl (L') \otimes \La (\Phi_0).
   $$
   At least up to a verification of signs, it is not difficult to see that $d_{\infty}^2$ vanishes 
   (note that this vanishing represents an infinite sequence of equations of which the first 
   two are: $d^2 = 0$ where $d$ is the usual cluster differential, and $\cyl^{\Phi_0} d - d \, \cyl^{\Phi_0} = 0$ 
   which is another way of proving that $\cyl^{\Phi_0}$ is a morphism of chain complexes). One could of course 
   define the big differential $d_{\infty}$ in which one would not limit to the class $K$ the choices of the 
   cylinders (but would instead consider all of them).

      The computation of the corresponding big homology seems an interesting problem. By analogy with the 
      differential of the rational model of the total space of a fibration, which is defined on the tensor 
      product of the complexes of the base and fiber and is ``twisted'' only on the fiber, it seems that 
      what this big differential computes is a quantum version of a fibration whose base would be $L'$ and 
      fiber $L$. But the geometric significance of this quantum fibration is not yet clear.

\end{document}